\documentclass{article}
\usepackage{indentfirst}
\usepackage{amsmath}
\usepackage[utf8]{inputenc}
\usepackage[english]{babel}
\usepackage{amsthm}
\usepackage{bbm}
\usepackage{bm}
\usepackage{dsfont}
\usepackage{blindtext}
\usepackage{amssymb}
\usepackage{extarrows}
\usepackage{authblk}

\usepackage{amsmath,amsfonts,stmaryrd,amssymb} 

\usepackage{enumerate} 

\usepackage[ruled]{algorithm2e} 

\usepackage[framemethod=tikz]{mdframed} 

\usepackage{listings} 
\usepackage{geometry} 

\usepackage[utf8]{inputenc} 
\usepackage[T1]{fontenc} 

\usepackage{XCharter} 

\mdfdefinestyle{commandline}{
	leftmargin=10pt,
	rightmargin=10pt,
	innerleftmargin=15pt,
	middlelinecolor=black!50!white,
	middlelinewidth=3pt,
	frametitlerule=false,
	backgroundcolor=black!5!white,
	frametitle={Command Line},
	frametitlefont={\normalfont\sffamily\color{white}\hspace{-1em}},
	frametitlebackgroundcolor=black!50!white,
	nobreak,
}

\mdfdefinestyle{file}{
	innertopmargin=4.0\baselineskip,
	innerbottommargin=3.5\baselineskip,
	topline=false, bottomline=false,
	leftline=false, rightline=false,
	leftmargin=2cm,
	rightmargin=2cm,
	singleextra={%
		\draw[fill=black!10!white](P)++(0,-1.2em)rectangle(P-|O);
		\node[anchor=north west]
		at(P-|O){\ttfamily\mdfilename};
		\def\l{3em}
		\draw(O-|P)++(-\l,0)--++(\l,\l)--(P)--(P-|O)--(O)--cycle;
		\draw(O-|P)++(-\l,0)--++(0,\l)--++(\l,0);
	},
	nobreak,
}

\mdfdefinestyle{question}{
	innertopmargin=3.5\baselineskip,
	innerbottommargin=3.0\baselineskip,
	roundcorner=5pt,
	nobreak,
	singleextra={%
		\draw(P-|O)node[xshift=1em,anchor=west,fill=white,draw,rounded corners=5pt]{%
		Question \theQuestion\questionTitle};
	},
}

\newcounter{Question} 

\mdfdefinestyle{warning}{
	topline=false, bottomline=false,
	leftline=false, rightline=false,
	nobreak,
	singleextra={%
		\draw(P-|O)++(-0.5em,0)node(tmp1){};
		\draw(P-|O)++(0.5em,0)node(tmp2){};
		\fill[black,rotate around={45:(P-|O)}](tmp1)rectangle(tmp2);
		\node at(P-|O){\color{white}\scriptsize\bf !};
		\draw[very thick](P-|O)++(0,-1em)--(O);
	}
}

\mdfdefinestyle{info}{%
	topline=false, bottomline=false,
	leftline=false, rightline=false,
	nobreak,
	singleextra={%
		\fill[black](P-|O)circle[radius=0.4em];
		\node at(P-|O){\color{white}\scriptsize\bf i};
		\draw[very thick](P-|O)++(0,-0.8em)--(O);
	}
}

\newtheorem{theorem}{Theorem}[section]

\newtheorem{proposition}{Proposition}[theorem]
\newtheorem{lemma}[theorem]{Lemma}
\theoremstyle{remark}
\newtheorem*{remark}{Remark}
\theoremstyle{definition}
\newtheorem{definition}{Definition}[section]

\title{Approximation for the invariant measure with applications for jump processes (convergence in total variation distance)}
\author[*]{Vlad Bally}
\author[*]{Yifeng Qin}

\affil[*]{Universit\'e Gustave Eiffel, LAMA (UMR CNRS, UPEMLV, UPEC), MathRisk INRIA, F-77454 Marne-la-Vall\'ee, France.
Email address: bally@univ-mlv.fr}
\date{2023}

\begin{document}

\maketitle

\textbf{Abstract} In this paper, we establish an abstract framework for the
approximation of the invariant probability measure for a Markov semigroup.
Following Pag$\grave{e}$s and Panloup $\cite{ref7}$ we use an Euler scheme
with decreasing step (unadjusted Langevin algorithm). Under some contraction
property with exponential rate and some regularization properties, we give
an estimate of the error in total variation distance. This abstract
framework covers the main results in $\cite{ref7}$ and $\cite{ref14}.$ As a
specific application we study the convergence in total variation distance to
the invariant measure for jump type equations. The main technical difficulty
consists in proving the regularzation properties - this is done under an
ellipticity condition, using Malliavin calculus for jump processes.


\tableofcontents

\section{Introduction}

The aim of this paper is to study the convergence to the invariant measure
of a Markov process. We refer to $\cite{ref23},\cite{ref34},\cite{ref35}$
for the existence of an invariant probability measure for a general Markov
process and to $\cite{ref38},$ $\cite{ref39}$ for some basic computation of
the invariant probability measure for a L\'{e}vy process. Following the
ideas from Pag$\grave{e}$s and Panloup $\cite{ref7}$ (see also Lamberton and
Pag$\grave{e}$s $\cite{ref53}$ $\cite{ref54}$) we use an Euler scheme with
decreasing step (known in the literature as the unadjusted Langevin
algorithm) in order to construct our algorithm (this has been studied in
depth in $\cite{ref41})$.

Our paper has two parts. In the first part we construct an abstract
framework which is appropriate in order to state and discuss our
approximation problem. We focus on the estimate of the error in total
variation distance. And the main achievement is to give some sufficient
regularization properties for the semigroup and for the Euler scheme, which
allow to treat bounded and measurable test functions. Furthermore, in order
to check such regularization properties, one has to use integration by parts
techniques inspired from Malliavin calculus. We give a regularization lemma
based on such arguments, which is the crucial step in our approach (it has
its own interest, beyond the application in this particular framework). Let
us mention that the abstract framework settled in our paper encompass the
following recent results: in $\cite{ref7}$, the authors use unadjusted
Langevin algorithm to approximate the invariant probability measure of a
diffusion process and study the Wasserstein and total variation distance
between them. In $\cite{ref14}$, the authors approximate the invariant
probability measure of a L\'{e}vy process but only study the Wasserstein
distance.

In the second part of the paper we illustrate our results in the case of
jump type $SDE^{\prime }s.$ In order to do it we recall the Malliavin
calculus for jump processes and prove estimates of the Sobolev norms and of
the Malliavin covariance matrix for the solution of such equations. These
estimates are rather long and technical, but at a certain extend they come
back on results already obtained in $\cite{ref20}$. Once these estimates are
proved, we apply the abstract results from the first part and obtain the
estimate of the error in total variation distance.

Let us present in more detail our results. We give in Section 2 the abstract
framework of the approximation for the invariant probability measure. We
denote $C_{b}^{l}(\mathbb{R}^{d})$ the space of $l-$times differential and
bounded functions on $\mathbb{R}^{d}$ with bounded derivatives up to order $%
l $. We consider a semigroup $P_{t},t\geq 0$ on the space $\mathcal{M}_{b}(%
\mathbb{R}^{d})$\ of the bounded measurable functions on $\mathbb{R}^{d}$
and assume that there exists at least one invariant probability measure $\nu 
$ for the semigroup $P_{t},t\geq 0$. We assume moreover the "exponential
Lipschitz property": there exists two constants $C_{0}\geq 1$ and $\rho >0$
such that for every $t>0$ and every $\varphi \in C_{b}^{1}(\mathbb{R}^{d})$%
\[
(L_{0})\quad \left\Vert \nabla P_{t}\varphi \right\Vert _{\infty }\leq
C_{0}\left\Vert \nabla \varphi \right\Vert _{\infty }e^{-\rho t}.
\]%
This immediately implies that $\nu $ is unique.

In order to approximate the invariant measure $\nu$, we introduce an Euler
scheme with decreasing time steps (unadjusted Langevin algorithm). For every 
$\gamma >0$ we give an operator $\overline{P}_{\gamma }:C_{b}^{\infty
}\rightarrow C_{b}^{\infty }$ such that $\Vert\overline{P}%
_\gamma\varphi\Vert_\infty\leq\Vert\varphi\Vert_\infty$ and which
approximates our semigroup in the following sense: for every $\gamma >0$%
\[
A(k_{0},\alpha )\quad \left\Vert (P_{\gamma }-\overline{P}_{\gamma })\varphi
\right\Vert _{\infty }\leq C_{k_{0}}\left\Vert \nabla \varphi \right\Vert
_{k_{0},\infty }\gamma ^{1+\alpha }.
\]%
Here $\alpha >0$ is a given number, $k_{0}\in \mathbb{N}$ and $\left\Vert
\psi \right\Vert _{k_{0},\infty }=\sum\limits_{\left\vert \alpha \right\vert
\leq k_{0}}\left\Vert \partial ^{\alpha }\psi \right\Vert _{\infty }. $ We
consider a decreasing sequence of time steps $\gamma _{n}\downarrow 0$ and
define the time grid $\Gamma _{n}=\sum\limits_{i=1}^{n}\gamma _{i}.$ We
assume that%
\[
(\Gamma )\quad \sum_{i=1}^{\infty }\gamma _{i}=\lim_{n\rightarrow \infty
}\Gamma _{n}=\infty .
\]%
We also introduce 
\[
\overline{\omega }=\overline{\omega }((\gamma _{n})_{n\in N})=\overline{%
\lim_{n\rightarrow \infty }}\frac{\gamma _{n}-\gamma _{n+1}}{\gamma
_{n+1}^{2}}<\infty.
\]%
The typical example is $\gamma _{n}=\frac{1}{n}$ and then $\overline{\omega }%
=1.$ In the following we denote $\{\Gamma \}=\{\Gamma _{n},n\in \mathbb{N}%
\}. $ And, for $\Gamma _{i}\leq t<\Gamma _{i+1}$ we denote $N(t)=i$ and $%
\tau (t)=\Gamma _{i}. $ Then, for $s\in \{\Gamma \}$ and $t\in \{\Gamma \}$
we define the Euler scheme%
\begin{equation}
\overline{P}_{s,t}=\prod_{i=N(s)}^{N(t)-1}\overline{P}_{\gamma _{i}},
\label{Pbar}
\end{equation}%
the product being understood in the sense of composition. This means that we
travel from $\tau (s)$ to $\tau (t)$ by using the Euler scheme associated to
the one step Euler scheme $\overline{P}_{\gamma }.$

So now we use the Euler scheme with decreasing time steps $\overline{P}%
_{0,\Gamma _{n}}$ (given in (\ref{Pbar})) to approximate the invariant
probability measure $\nu $. Our aim is to estimate the total variation
distance between them. To do so, we need some regularization properties.
First we give the regularization hypothesis concerning the semigroup $P_{t}$%
: 
\begin{eqnarray*}
R_{P}(k)\quad \sup_{1\leq t\leq 2}\left\Vert \nabla P_{t}\varphi \right\Vert
_{k-1,\infty } &\leq &C_{k}\left\Vert \varphi \right\Vert _{\infty },\quad and
\\
R_{P}^{\prime }(k)\quad \sup_{1\leq t\leq 2}\left\Vert \nabla P_{t}\varphi
\right\Vert _{k-1,\infty } &\leq &C_{k}^{\prime }\left\Vert \nabla \varphi
\right\Vert _{\infty }.
\end{eqnarray*}%
We also introduce the following variant of the Lipschitz property:
\begin{eqnarray}
\overline{L}_{k}\quad \left\Vert \nabla P_{t}\varphi \right\Vert _{k,\infty } &\leq
&C_{k}\left\Vert \nabla \varphi \right\Vert _{k,\infty },\quad 1\geq t>0. 
\nonumber
\end{eqnarray}

We give now the regularization properties for the Euler scheme $\overline{P}%
_{s,t}$. To begin, we introduce some notations. We fix a super kernel $\phi $
(see (18) for the precise definition), and, for $\delta \in (0,1]$ we denote 
$\phi _{\delta }(y)=\frac{1}{\delta ^{d}}\phi (\frac{y}{\delta })$. Moreover,
for a function $\varphi $ we denote $\varphi _{\delta }$ the regularization
by convolution with the super kernel: $\varphi _{\delta }=\varphi \ast \phi
_{\delta },$ with $\ast $ denoting convolution. For $\delta >0,\eta >0,$ and 
$q,\kappa ,p\in \mathbb{N}$ we denote%
\[
A_{q,\kappa ,p}^{\delta ,\eta }(h)=\frac{\delta ^{q}}{\eta ^{2q}}+\eta
^{-p}h^{p}+\eta ^{\kappa },\quad h>0.
\]%
Let $\beta >0$ and $p\geq 1$ be fixed and we assume the following
regularization property for the Euler scheme $\overline{P}_{s,t}$: we assume
that for every $q,\kappa \in \mathbb{N}$ there exists a constant $%
C=C_{q,\kappa ,p}$\ such that for every $\delta >0,\eta >0,$ every $%
1<t<r<t+2 $ and every bounded measurable function $\varphi $%
\begin{eqnarray*}
R_{\overline{\mathcal{P}}}(p,\beta )\quad &&\left\Vert \overline{P}%
_{t-1,t}P_{t,r}\varphi -\overline{P}_{t-1,t}P_{t,r}\varphi _{\delta
}\right\Vert _{\infty }+\left\Vert \overline{P}_{t-1,t}\overline{P}%
_{t,r}\varphi -\overline{P}_{t-1,t}\overline{P}_{t,r}\varphi _{\delta
}\right\Vert _{\infty } \\
&\leq &C_{q,\kappa ,p}\times A_{q,\kappa ,p}^{\delta ,\eta }(\gamma
_{N(t)}^{\beta })\left\Vert \varphi \right\Vert _{\infty }.
\end{eqnarray*}

Now we can give our main result (see \textbf{Proposition 2.1.1}). We assume
that an invariant probability measure $\nu $ exists for the semigroup $%
P_{t},t\geq 0$. We construct an Euler scheme with decreasing time steps $%
\overline{P}_{s,t}$ by (\ref{Pbar}). Suppose that $(L_{0})$ holds for some $\rho$, $%
A(k_{0},\alpha )$ holds for some $k_{0},\alpha $ with $\rho>\alpha\overline{\omega}$, $R_{P}(k)$, $%
R_{P}^{\prime }(k)$ and $\overline{L}_{k}$ hold for every $k$, and $R_{\overline{\mathcal{P}}%
}(p,\beta )$ holds true for some $p,\beta $. Then the invariant probability
measure $\nu $ is unique and for any $\varepsilon >0$,  for every $x\in 
\mathbb{R}^{d}$ and $n$ large enough,
\[
d_{TV}(\overline{P}_{0,\Gamma _{n}}(x,.),\nu )\leq C_{\varepsilon }(\gamma
_{n}^{((p\beta )\wedge \alpha )-\varepsilon }+\int_{\mathbb{R}%
^{d}}\left\vert x-y\right\vert d\nu (y)e^{-\rho \Gamma _{n}}).
\]%
We remark that we get the same speed of convergence as in $\cite{ref7}$ and $%
\cite{ref14}$, but in a more general framework.

\bigskip

We notice that we need some regularization properties (see $R_{P}(k)$, $%
R_{P}^{\prime }(k)$ and $R_{\overline{\mathcal{P}}}(p,\beta )$). In order to
obtain these properties, we introduce in Section 3 an abstract framework
built on a particular case of the Dirichlet form theory (see $\cite{ref25}$
and $\cite{ref2}$) in which such a property may be obtained by using some
integration by parts techniques. Those techniques are very similar to the
standard Malliavin calculus but are presented in a more general framework
which goes beyond the sole case of the Wiener space. In particular, we aim
at providing a minimalist setting leading to our regularization lemma. Our
unified framework includes the standard Malliavin calculus and different
known versions: the calculus based on the splitting method developed and
used in $\cite{ref26},\cite{ref21},\cite{ref27}$ as well as the $\Gamma -$%
calculus in $\cite{ref25}$. We also mention that our approach applies in the
case of the Malliavin calculus for jump type processes as settled by $\cite%
{ref5}$ and in the "lent particle" approach for Poisson point measures
developed by $\cite{ref28}$. 

\bigskip

In Section 4, we apply the results in Section 2 for jump processes. So we
consider the $d-$dimensional stochastic differential equation with jumps as
follows: 
\begin{eqnarray}
X_{t}&=&x+\int_{0}^{t}b(X_{r})dr +\int_{0}^{t}\int_{\mathbb{R}%
^d}c(z,X_{r-})N(dz,dr),  \label{Intro1}
\end{eqnarray}
where $N(dz,dr)$ is a Poisson point measure on the state space $\mathbb{R}^d$
with intensity measure $\widehat{N}(dz,dr)=\mu(dz)dr$, $x$ is the initial
value, $\mu$ is a positive $\sigma$-finite measure on $\mathbb{R}^d$, and $b:%
\mathbb{R}^d\rightarrow\mathbb{R}^d$, $c:\mathbb{R}^d\times\mathbb{R}%
^d\rightarrow\mathbb{R}^d$. Some basic background of jump processes can be
found in $\cite{ref6},\cite{ref8},\cite{ref36},\cite{ref37}$ and $\cite%
{ref51}$.

We need to give sufficient conditions to ensure the existence of an
invariant probability measure for the jump equation (\ref{Intro1}). We
recall by $\cite{ref23}$ the classical results of the existence of an
invariant probability measure for a general Markov process. Recently, $\cite%
{ref33}$ gives some specific criterias for the existence of an invariant
probability measure of a jump process and also discuss some ergodicity
properties. Here we suppose that (\textbf{Hypothesis 2.5}) 
\begin{eqnarray*}
i)\quad \left\langle x-y,b(x)-b(y)\right\rangle &\leq &-\overline{b}%
\left\vert x-y\right\vert ^{2} \\
ii)\quad \left\vert c(z,x)-c(z,y)\right\vert &\leq &\bar{c}(z)\left\vert
x-y\right\vert
\end{eqnarray*}%
and%
\[
iii)\quad 2\overline{b}-\int_{\mathbb{R}^{d}}(2\bar{c}(z)+\bar{c}^{2}(z))\mu
(dz):=\theta >0.
\]%
Our conditions are based on $\cite{ref23}$ and are essentially the same as
the conditions in $\cite{ref33}$. Indeed, the conditions above implies that
for some $\bar{\beta},\bar{\alpha}>0$ and a Lyapunov function $%
V(x)=\left\vert x\right\vert ^{2}$, we have ${L}V\leq \bar{\beta}-\bar{\alpha%
}V,$with ${L}$ denoting the infinitesimal operator of (\ref{Intro1}). This
guarantees the existence of an invariant probability measure $\nu $.

Moreover, in order to apply the Malliavin framework in Section 3 and obtain
regularization properties, we assume (see \textbf{Hypothesis 2.4} $b)$) that
the measure $\mu$ is absolutely continuous with respect to the Lebesgue
measure: $\mu(dz)=h(z)dz$, where $h$ is infinitely differentiable and $\ln h$
has bounded derivatives of any order. We also need some regularity and
ellipticity conditions on the coefficients (see \textbf{Hypothesis 2.1$\sim$%
2.3} for details). We mention that for every multi-indices $\beta_1,\beta_2$%
, we assume that there exists a non-negative function $\bar{c}:\mathbb{R}%
^d\rightarrow\mathbb{R}_{+}$ such that 
\[
\vert {c}(z,x)\vert+\vert \partial_z^{\beta_2}\partial_x^{\beta_1} {c}%
(z,x)\vert\leq\bar{c}(z),
\]%
with $\int_{\mathbb{R}^d}\vert \bar{c}(z)\vert^p\mu(dz)<\infty,\ \forall
p\geq1.$ We also assume that there exists a non-negative function $%
\underline{c}:\mathbb{R}^d\rightarrow\mathbb{R}_{+}$ such that for every $%
\zeta\in\mathbb{R}^d$, 
\[
\sum_{j=1}^d\langle\partial_{z_j}{c}(z,x),\zeta\rangle^{2}\geq \underline{c}%
(z)\vert\zeta\vert^2.
\]

Now we construct the Euler scheme. We take a partition with decreasing time
steps $\mathcal{P}=\{0=\Gamma_0<\Gamma_1<\cdots<\Gamma_{n-1}<\Gamma_n<\cdots%
\}$ with the time steps $\gamma_n=\Gamma_n-\Gamma_{n-1},\ n\in\mathbb{N}$
verifying some suitable conditions (see Section 4.3 for details). For $%
\Gamma _{n}\leq t<\Gamma _{n+1}$ we denote $\tau (t)=\Gamma _{n}.$ We
consider the Euler scheme: 
\begin{eqnarray*}
X^{\mathcal{P}}_{t}&=&x+\int_{0}^{t}b(X^{\mathcal{P}}_{\tau(r)})dr
+\int_{0}^{t}\int_{\mathbb{R}^d}{c}(z,X^{\mathcal{P}}_{\tau(r)-})N(dz,dr).
\end{eqnarray*}%
Some results concerning the convergence of the Euler scheme of a jump
equation can be found for example in $\cite{ref16},\cite{ref45},\cite{ref46},%
\cite{ref47},\cite{ref48},\cite{ref49}$ and $\cite{ref50}$.

Since $\mu (\mathbb{R}^{d})=\infty $ (which is a consequence of \textbf{%
Hypothesis 2.4 $a)$}), we have infinitely many jumps. So we construct the
truncated Euler scheme in order to have finite numbers of jumps for the sake
of simulation and Malliavin calculus. For $m\in \mathbb{N}$, we denote $%
B_{m}=\{z\in \mathbb{R}^{d}:|z|\leq m\}$ and denote 
\[
\varepsilon _{m}:=\int_{\{|z|>m\}}|\bar{c}(z)|^{2}\mu (dz)+|\int_{\{|z|>m\}}%
\bar{c}(z)\mu (dz)|^{2}.
\]%
For every $\gamma >0$, we define the truncation function $M(\gamma )\in \mathbb{N}$
to be the smallest integer such that 
\[
\varepsilon _{M(\gamma )}\leq \gamma ^{2}.
\]%
For $\Gamma _{n}<t\leq \Gamma _{n+1}$, we denote $M_{{\mathcal{P}}%
}(t)=M(\gamma _{n+1})$. Now we cancel the "big jumps" (the jumps of size $%
|z|>M_{{\mathcal{P}}}(t)$): 
\begin{equation}
X_{t}^{\mathcal{P},{M_{\mathcal{P}}}}=x +\int_{0}^{t}b(X_{\tau (r)}^{%
\mathcal{P},{M_{\mathcal{P}}}})dr+\int_{0}^{t}\int_{B_{M_{{\mathcal{P}}}(r)}}%
{c}(z,X_{\tau (r)-}^{\mathcal{P},{M_{\mathcal{P}}}})N(dz,dr).  \label{Intro2}
\end{equation}%
We remark that the solution of the equation (\ref{Intro2}) can be
constructed in an explicit way.

Then we apply the abstract framework in Section 2 for $X_{\Gamma _{n}}^{%
\mathcal{P},{M_{\mathcal{P}}}}$ and obtain the following main result (see 
\textbf{Theorem 4.1}): An invariant probability measure $\nu $ of the jump
equation (\ref{Intro1}) exists and is unique, and for any $\varepsilon >0$,
there exists a constant $C_{\varepsilon }$ such that for every $x\in 
\mathbb{R}^{d}$ and $n$ large enough, we have
\[
d_{TV}(\mathcal{L}(X_{\Gamma _{n}}^{\mathcal{P},{M_{\mathcal{P}}}}),\nu
)\leq C_{\varepsilon }(\gamma _{n}^{1-\varepsilon }+\int_{\mathbb{R}%
^{d}}\left\vert x-y\right\vert d\nu (y)e^{-\frac{\theta }{2}\Gamma _{n}}),
\]%
with $\mathcal{L}(X)$ denoting the law of a random variable $X$. We notice
that we obtain the same speed of convergence as in $\cite{ref7}$ but $\cite%
{ref7}$ concern the diffusion process driven by a Brownian motion while here
we consider the jump process. Comparing with the results in $\cite{ref14}$,
we also obtain the same speed of convergence but $\cite{ref14}$ only deals
with the Wasserstein distance while in our paper, we deal with the total
variation distance.

\bigskip\bigskip

\section{Approximation of the invariant measure: Abstract framework}

\subsection{The semigroup and the invariant measure}

We consider a semigroup $P_{t},t\geq 0$ on the space $\mathcal{M}_{b}(%
\mathbb{R}^{d})$\ of the bounded measurable functions on $\mathbb{R}^{d}$.
We denote $C_b^{l}(\mathbb{R}^d)$ the space of $l-$times differential and
bounded functions on $\mathbb{R}^d$ with bounded derivatives up to order $l$%
. We will use the following two hypotheses:

$(I)\qquad $We assume that there exists at least one invariant distribution
for the semigroup $P_{t},t\geq 0.$

Moreover we assume the following "exponential Lipschitz property": we assume
that there exists two constants $C_{0}\geq 1$ and $\rho >0$ such that for
every $t>0$ and every $\varphi \in C_{b}^{1}(\mathbb{R}^{d})$%
\begin{equation}
(L_{0})\quad \left\Vert \nabla P_{t}\varphi \right\Vert _{\infty }\leq
C_{0}\left\Vert \nabla \varphi \right\Vert _{\infty }e^{-\rho t}.
\label{app1}
\end{equation}

We also denote by $\mathcal{P}_{1}$\ the space of the probability measures
on $\mathbb{R}^{d}$ which have finite moment of order one $\int_{\mathbb{R}^d} \left\vert
x\right\vert \nu (dx)<\infty .$ This is a Banach space under the Wasserstein
distance $W_{1}$:%
\[
W_{1}(\nu ,\mu )=\sup \{\left\vert \int_{\mathbb{R}^d} \varphi d(\nu -\mu
)\right\vert :\left\Vert \nabla \varphi \right\Vert _{\infty }\leq 1\}.
\]

\begin{proposition}
Suppose that the semigroup $P_{t},t\geq 0$ has at least an invariant
probability measure $\nu $ and that (\ref{app1}) holds true. Then the
invariant probability measure is unique and moreover, for every $x\in 
\mathbb{R}^{d}$ 
\begin{equation}
W_{1}(\nu ,P_{t}(x,\cdot ))\leq C\int_{\mathbb{R}^{d}}\left\vert
x-y\right\vert \nu (dy)\times e^{-\rho t}.  \label{app2}
\end{equation}
\end{proposition}

\textbf{Proof. Step 1} We will prove that for sufficiently large $t,$ the
application $\nu \mapsto \nu P_{t}$ is a strict contraction on the
Wassertein space: using (\ref{app1}),%
\begin{eqnarray*}
\left\vert \int_{\mathbb{R}^d} \varphi (y)d(\nu P_{t}-\mu
P_{t})(dy)\right\vert &=&\left\vert \int_{\mathbb{R}^d} P_{t}\varphi
(x)d(\nu (x)-\mu (x))\right\vert \\
&\leq &\left\Vert \nabla P_{t}\varphi \right\Vert _{\infty }W_{1}(\nu ,\mu )
\\
&\leq &C_{0}\left\Vert \nabla \varphi \right\Vert _{\infty }e^{-\rho
t}W_{1}(\nu ,\mu ).
\end{eqnarray*}%
This means that, for large $t$ 
\[
W_{1}(\nu P_{t},\mu P_{t})\leq C_{0}e^{-\rho t}W_{1}(\nu ,\mu )\leq \frac{1}{%
2}W_{1}(\nu ,\mu )
\]%
and this guarantees the uniqueness of the invariant measure.

\textbf{Step 2 }Since $\nu $ is an invariant measure 
\[
\int_{\mathbb{R}^d} \varphi (z)\nu (dz)=\int_{\mathbb{R}^d} \int_{\mathbb{R}%
^d} P_{t}(z,dy)\varphi (y)\nu (dz)
\]%
which gives, for every fixed $x\in \mathbb{R}^{d}$ ($\nu $ is a probability) 
\begin{eqnarray}
\int_{\mathbb{R}^d} \varphi (z)\nu (dz)-\int_{\mathbb{R}^d}
P_{t}(x,dy)\varphi (y) &=&\int_{\mathbb{R}^d} \int_{\mathbb{R}^d}
(P_{t}(z,dy)-P_{t}(x,dy))\varphi (y)\nu (dz)  \label{app2'} \\
&=&\int_{\mathbb{R}^d} (P_{t}\varphi (z)-P_{t}\varphi (x))\nu (dz)  \nonumber
\end{eqnarray}%
so that 
\begin{eqnarray*}
\left\vert \int_{\mathbb{R}^d} \varphi (z)\nu (dz)-\int_{\mathbb{R}^d}
P_{t}(x,dy)\varphi (y)\right\vert &\leq &\left\Vert \nabla P_{t}\varphi
\right\Vert _{\infty }\int_{\mathbb{R}^{d}}\left\vert x-z\right\vert \nu (dz)
\\
&\leq &C_{0}e^{-\rho t}\left\Vert \nabla \varphi \right\Vert _{\infty }\int_{%
\mathbb{R}^{d}}\left\vert x-z\right\vert \nu (dz)
\end{eqnarray*}%
which yields (\ref{app2}). $\square $

\bigskip

\subsection{The Euler scheme}

We introduce now an Euler scheme with decreasing steps. First, for every $%
\gamma >0$ we give an operator $\overline{P}_{\gamma }:C_{b}^{\infty
}\rightarrow C_{b}^{\infty }$ such that $\Vert\overline{P}%
_\gamma\varphi\Vert_\infty\leq\Vert\varphi\Vert_\infty$ and which
approximates our semigroup in the following sense: for every $\gamma >0$%
\begin{equation}
A(k_{0},\alpha )\quad \left\Vert (P_{\gamma }-\overline{P}_{\gamma })\varphi
\right\Vert _{\infty }\leq C_{k_{0}}\left\Vert \nabla \varphi \right\Vert
_{k_{0},\infty }\gamma ^{1+\alpha }.  \label{app3}
\end{equation}%
Here $\alpha >0$ is a given number, $k_{0}\in \mathbb{N}$ and 
\[
\left\Vert \psi \right\Vert _{k_{0},\infty }=\sum_{\left\vert \alpha
\right\vert \leq k_{0}}\left\Vert \partial ^{\alpha }\psi \right\Vert
_{\infty }.
\]

Moreover, we consider a decreasing sequence of time steps $\gamma
_{n}\downarrow 0$ and define the time grid $\Gamma _{n}=\sum_{i=1}^{n}\gamma
_{i}.$ We assume that%
\begin{equation}
(\Gamma )\quad \sum_{i=1}^{\infty }\gamma _{i}=\lim_{n\rightarrow \infty
}\Gamma _{n}=\infty .  \label{app4}
\end{equation}%
We also introduce 
\[
\overline{\omega }=\overline{\omega }((\gamma _{n})_{n\in N})=\overline{%
\lim_{n\rightarrow \infty }}\frac{\gamma _{n}-\gamma _{n+1}}{\gamma
_{n+1}^{2}}<\infty.
\]%
The typical example is $\gamma _{n}=\frac{1}{n}$ and then $\overline{\omega }%
=1.$ In the following we denote $\{\Gamma \}=\{\Gamma _{n},n\in \mathbb{N}%
\}. $ And, for $\Gamma _{i}\leq t<\Gamma _{i+1}$ we denote%
\[
N(t)=i\quad and\quad \tau (t)=\Gamma _{i}.
\]%
In particular, for $t=\Gamma _{i}\in \{\Gamma \}$ we have $N(t)=i$ such that 
$t=\Gamma _{N(t)}.$ Then, for $s\in \{\Gamma \}$ and $t\in \{\Gamma \}$ we
define the Euler scheme%
\begin{equation}
\overline{P}_{s,t}=\prod_{i=N(s)}^{N(t)-1}\overline{P}_{\gamma _{i}}
\label{app6}
\end{equation}%
the product being understood in sense of composition. This means that we
travel from $\tau (s)$ to $\tau (t)$ by using the Euler scheme associated to
the one step Euler scheme $\overline{P}_{\gamma }.$ In the appendix 7.1 we
will prove the following lemma (which is a slight generalisation of the
lemma given by Pages and Panloup $\cite{ref7}$): for every $\rho >\alpha 
\overline{\omega },$ there exists $n_{\rho }$ and $C_{\rho }$ such that for $%
n\geq n_{\rho }$%
\begin{equation}
\sum_{i=1}^{n}\gamma _{i}^{1+\alpha }e^{-\rho (\Gamma _{n}-\Gamma _{i})}\leq
C_{\rho }\gamma _{n}^{\alpha }.  \label{app5}
\end{equation}%
Moreover, there exists $n_{\ast }$ such that, for $n_{\ast }\leq i\leq n$ 
\begin{equation}
\gamma _{i}\leq e^{2\overline{\omega }(\Gamma _{n}-\Gamma _{i})}\gamma _{n}.
\label{app}
\end{equation}

Notice that $P_{t},t\geq 0$ is a homogeneous semigroup, and we may define $%
P_{s,t}=P_{t-s}=P_{0,t-s}.$ In contrast, $\overline{P}_{s,t},s<t,$ is not
homogeneous: we do not have $\overline{P}_{s,t}=\overline{P}_{0,t-s}.$ This
is due to the fact that the greed $\Gamma _{i},i\in \mathbb{N}$ is not
uniform.

Finally we assume the following stronger variant of the Lipschitz property $%
L_{0}$:%
\begin{equation}
(L_{k_{0}})\quad \left\Vert \nabla P_{t}\varphi \right\Vert _{k_{0},\infty
}\leq C_{k_{0}}\left\Vert \nabla \varphi \right\Vert _{k_{0},\infty
}e^{-\rho t}  \label{app7}
\end{equation}%
where $k_{0}$ is the one from $A(k_{0},\alpha ).$

\begin{proposition}
Suppose that (\ref{app3}) and (\ref{app7}) hold true with $\rho >\alpha 
\overline{\omega }$. Then for $N(t)>n_\rho+1$, we have
\begin{equation}
\left\Vert (P_{s,t}-\overline{P}_{s,t})\varphi \right\Vert _{\infty }\leq
C_{k_{0}}\left\Vert \nabla \varphi \right\Vert _{k_{0},\infty }\gamma
_{N(t)}^{\alpha }.  \label{app8}
\end{equation}
\end{proposition}

\textbf{Proof }We use (\ref{app3}) first and (\ref{app7}) then%
\begin{eqnarray*}
\left\Vert (P_{s,t}-\overline{P}_{s,t})\varphi \right\Vert _{\infty } &\leq
&\sum_{i=N(s)}^{N(t)-1}\left\Vert \overline{P}_{s,\Gamma _{i-1}}(\overline{P}%
_{\gamma _{i}}-P_{\gamma _{i}})P_{\Gamma _{i},t}\varphi \right\Vert _{\infty
} \\
&\leq &\sum_{i=N(s)}^{N(t)-1}\left\Vert (\overline{P}_{\gamma
_{i}}-P_{\gamma _{i}})P_{\Gamma _{i},t}\varphi \right\Vert _{\infty } \\
&\leq &C_{k_{0}}\sum_{i=N(s)}^{N(t)-1}\left\Vert \nabla P_{\Gamma
_{i},t}\varphi \right\Vert _{k_{0},\infty }\gamma _{i}^{1+\alpha } \\
&\leq &C_{k_{0}}^{\prime }\sum_{i=N(s)}^{N(t)-1}\left\Vert \nabla \varphi
\right\Vert _{k_{0},\infty }\gamma _{i}^{1+\alpha }e^{-\rho (\Gamma
_{N(t)}-\Gamma _{i})} \\
&\leq &C_{k_{0}}^{\prime \prime }\left\Vert \nabla \varphi \right\Vert
_{k_{0},\infty }\gamma _{N(t)}^{\alpha }.
\end{eqnarray*}%
For the last inequality we have used (\ref{app5}). $\square $

\begin{remark}
Suppose that (\ref{app3}) and (\ref{app7}) hold with $k_0=0$. We also
suppose that an invariant probability measure $\nu $ of the semigroup $%
P_{t},t\geq 0$ exists and that (\ref{app1}) holds true. Then \textbf{%
Proposition 2.0.1} and \textbf{Proposition 2.0.2} give that for every $x\in 
\mathbb{R}^{d}$, we have 
\[
W_{1}(\nu ,\overline{P}_{0,t}(x,\cdot ))\leq C(\gamma _{N(t)}^{\alpha
}+\int_{\mathbb{R}^{d}}\left\vert x-y\right\vert \nu (dy)\times e^{-\rho t}).
\]%
For this result, we do not need any regularization properties. In order to
obtain the result for the total variation distance, we give some
regularization properties in the next subsection.
\end{remark}

\bigskip

\subsection{Regularization properties}

In this section we will assume that the semigroup and the Euler scheme have
some regularization properties which allow to obtain convergence in total
variation distance.

First we give the regularization hypothesis concerning the semigroup: 
\begin{eqnarray}
R_{P}(k)\quad \sup_{1\leq t\leq 2}\left\Vert \nabla P_{t}\varphi \right\Vert
_{k-1,\infty } &\leq &C_{k}\left\Vert \varphi \right\Vert _{\infty },\quad
and  \label{app9} \\
R_{P}^{\prime }(k)\quad \sup_{1\leq t\leq 2}\left\Vert \nabla P_{t}\varphi
\right\Vert _{k-1,\infty } &\leq &C_{k}^{\prime }\left\Vert \nabla \varphi
\right\Vert _{\infty },  \label{app9'}
\end{eqnarray}

Such a regularization property is proved using the integration by parts
formula in Malliavin calculus. Moreover, we suppose that we have the
following variant of the Lipschitz property: 
\begin{eqnarray}
\overline{L}_{k}\quad i)\quad \left\Vert \nabla P_{t}\varphi \right\Vert
_{\infty } &\leq &C_{k}\left\Vert \nabla \varphi \right\Vert _{k,\infty
}e^{-\rho t},\quad t\geq 1,  \label{app10} \\
ii)\quad \left\Vert \nabla P_{t}\varphi \right\Vert _{k,\infty } &\leq
&C_{k}\left\Vert \nabla \varphi \right\Vert _{k,\infty },\quad 1\geq t>0. 
\nonumber
\end{eqnarray}%
Notice that $\overline{L}_{k},i)$ is weaker then $L_{0}$ (see (\ref{app1}))
because we have $\left\Vert \nabla \varphi \right\Vert _{k,\infty }$ instead
of $\left\Vert \nabla \varphi \right\Vert _{\infty }.$ However, if the
regularization property $R_{P}^{\prime }(k)$ holds then $\overline{L}_{k},i)$
implies $L_{0}$ (for $t\geq 1).$ Indeed, $\overline{L}_{k}$ gives%
\begin{eqnarray*}
\left\Vert \nabla P_{t}\varphi \right\Vert _{\infty } &=&\left\Vert \nabla
(P_{t-1}P_{1}\varphi )\right\Vert _{\infty }\leq C\left\Vert \nabla
P_{1}\varphi \right\Vert _{k,\infty }e^{-\rho (t-1)} \\
&\leq &C\left\Vert \nabla \varphi \right\Vert _{\infty }e^{-\rho (t-1)},
\end{eqnarray*}%
the last inequality being the consequence of $R_{P}^{\prime }(k).$ In
particular, if an invariant probability measure $\nu $ exists, then it is
unique and we have (\ref{app2}). 

\begin{remark}
We also notice that $R_{P}^{\prime }(k+1)$ and $\overline{L}_{k}$ imply $%
L_{k}.$ Indeed, for $t\leq 1,\overline{L}_{k}\ ii)$ gives 
\[
\left\Vert \nabla P_{t}\varphi \right\Vert _{k,\infty }\leq C_{k}\left\Vert
\nabla \varphi \right\Vert _{k,\infty }\leq e^{\rho }C_{k}\left\Vert \nabla
\varphi \right\Vert _{k,\infty }e^{-\rho t}
\]%
and for $t\geq 1$ 
\begin{eqnarray*}
\left\Vert \nabla P_{t}\varphi \right\Vert _{k,\infty } &=&\left\Vert \nabla
(P_{1}P_{t-1}\varphi )\right\Vert _{k,\infty }\leq C\left\Vert \nabla
P_{t-1}\varphi \right\Vert _{\infty } \\
&\leq &C\left\Vert \nabla \varphi \right\Vert _{k,\infty }e^{-\rho (t-1)}.
\end{eqnarray*}
\end{remark}

\bigskip 

Moreover, for $t\geq 1$, $\overline{L}_{k}$ and $R_{P}(k+1)$ give 
\begin{equation}
d_{TV}(P_{t}(x,.),\nu )\leq C(\int_{\mathbb{R}^{d}}\left\vert x-y\right\vert
d\nu (y))e^{-\rho t},  \label{App5}
\end{equation}%
where $d_{TV}$ denotes the total variation distance: 
\[
d_{TV}(\mu ,\nu )=\sup\limits_{\Vert f\Vert _{\infty }\leq 1}\big\vert\int_{%
\mathbb{R}^{d}}f(x)\mu (dx)-\int_{\mathbb{R}^{d}}f(x)\nu (dx)\big\vert.
\]%
Indeed, 
\begin{eqnarray*}
\left\vert P_{t}\varphi (x)-P_{t}\varphi (y)\right\vert  &=&\left\vert
P_{t-1}P_{1}\varphi (x)-P_{t-1}P_{1}\varphi (y)\right\vert  \\
&\leq &C_{k}\left\Vert \nabla P_{1}\varphi \right\Vert _{k,\infty }e^{-\rho
(t-1)}\left\vert x-y\right\vert  \\
&\leq &C_{k}C_{k+1}e^{\rho }\left\Vert \varphi \right\Vert _{\infty
}e^{-\rho t}\left\vert x-y\right\vert .
\end{eqnarray*}%
Then we come back to (\ref{app2'}) and we obtain%
\[
\left\vert \int_{\mathbb{R}^{d}}\varphi (z)\nu (dz)-\int_{\mathbb{R}%
^{d}}P_{t}(x,dy)\varphi (y)\right\vert \leq C\left\Vert \varphi \right\Vert
_{\infty }\int_{\mathbb{R}^{d}}e^{-\rho t}\left\vert x-y\right\vert \nu (dy)
\]%
so (\ref{App5}) is proved. $\square $

\bigskip

We give now the regularization properties for the Euler scheme; this is a
more delicate subject, because we have some difficulties in order to use
directly the Malliavin calculus for the Euler scheme (the reason is that the
decomposition using the inverse of the tangent flow does not work, and so
the proof of the non degeneracy property is more difficult) .

We introduce some notations. We recall that a super kernel $\phi :\mathbb{R}%
^{d}\rightarrow \mathbb{R}$ is a function which belongs to the Schwartz
space and such that for every multi-indexes $\beta _{1}$ and $\beta _{2}$,
one has 
\begin{equation}
\int_{\mathbb{R}^d} \phi (x)dx=1,\quad \int_{\mathbb{R}^d} y^{\beta
_{1}}\phi (y)dy=0\quad \text{for}\quad |\beta _{1}|\geq 1,\quad \int_{%
\mathbb{R}^d} |y|^{m}|\partial _{\beta _{2}}\phi (y)|dy<\infty \quad \text{%
for}\quad m\in \mathbb{N}.  \label{SuperKernel}
\end{equation}%
We fix a super kernel $\phi $. For $\delta \in (0,1]$, we denote $\phi
_{\delta }(y)=\frac{1}{\delta ^{d}}\phi (\frac{y}{\delta }) $ and $\varphi
_{\delta }$ the regularization by convolution with a super kernel: 
\begin{equation}
\varphi _{\delta }=\varphi \ast \phi _{\delta },  \label{superkernel}
\end{equation}%
with $\ast $ denoting convolution.

As usual, for a multi-index $\beta _{1}=(\beta _{1}^{1},\cdots ,\beta
_{1}^{m})\in \{1,\cdots ,d\}^{m}$, one denotes $|\beta _{1}|=m$ and $%
y^{\beta _{1}}=\prod_{i=1}^{m}y_{\beta _{1}^{i}}$.

For $\delta >0,\eta >0,$ and $q,\kappa ,p\in \mathbb{N}$ we denote%
\[
A_{q,\kappa ,p}^{\delta ,\eta }(h)=\frac{\delta ^{q}}{\eta ^{2q}}+\eta
^{-p}h^{p}+\eta ^{\kappa },\quad h>0.
\]%
Then we assume the following:

Let $\beta >0$ and $p\geq 1$ be fixed. We assume that for every $q,\kappa
\in \mathbb{N}$ there exists a constant $C=C_{q,\kappa ,p}$\ such that for
every $\delta >0,\eta >0,$ every $1<t<r<t+2 $ and every bounded measurable
function $\varphi $%
\begin{eqnarray}
R_{\overline{\mathcal{P}}}(p,\beta )\quad &&\left\Vert \overline{P}%
_{t-1,t}P_{t,r}\varphi -\overline{P}_{t-1,t}P_{t,r}\varphi _{\delta
}\right\Vert _{\infty }+\left\Vert \overline{P}_{t-1,t}\overline{P}%
_{t,r}\varphi -\overline{P}_{t-1,t}\overline{P}_{t,r}\varphi _{\delta
}\right\Vert _{\infty }  \label{app11} \\
&&\leq C_{q,\kappa,p}\times A_{q,\kappa ,p}^{\delta ,\eta }(\gamma
_{N(t-1)}^{\beta })\left\Vert \varphi \right\Vert _{\infty }.  \nonumber
\end{eqnarray}

This represents the "regularization property for $\overline{P}_{t-1,t}".$ In
order to prove it, one employs \textbf{Lemma 3.5} (see (\ref{1*})) in Section 3.1.

As a consequence of these  properties, we obtain the following lemma. We recall $n_\rho$ and $n_\ast$ in (\ref{app5}) and (\ref{app}).

\begin{lemma}
We fix $\beta >0$ and $p\geq 1.$Suppose that (\ref{app3}) (\ref{app7}) hold with $\rho >\alpha 
\overline{\omega }$, and $%
R_{\overline{\mathcal{P}}}(p,\beta )$ (see (\ref{app11})) holds. Then, for
every $\varepsilon >0$ there exists a constant $C_{\varepsilon }\geq 1$ such
that for every $s<t-1<t<r<t+2$ with $N(r)>n_\rho+1$ and $N(t-1)>n_\ast$, and for every bounded measurable function $%
\varphi $ 
\begin{equation}
\left\Vert \overline{P}_{s,t}(\overline{P}_{t,r}-P_{t,r})\varphi \right\Vert
_{\infty }\leq C_{\varepsilon }\left\Vert \varphi \right\Vert _{\infty
}\gamma _{N(t)}^{((p\beta )\wedge \alpha )-\varepsilon }.  \label{app12}
\end{equation}
\end{lemma}

\bigskip

\textbf{Proof }We use (\ref{app11}) and (\ref{app}) in order to get 
\begin{eqnarray*}
\left\Vert \overline{P}_{s,t}(\overline{P}_{t,r}-P_{t,r})\varphi \right\Vert
_{\infty } &\leq &\left\Vert \overline{P}_{t-1,t}(\overline{P}%
_{t,r}-P_{t,r})\varphi \right\Vert _{\infty } \\
&\leq &C_{q,\kappa,p }\left\Vert \varphi \right\Vert _{\infty }\times
A_{q,\kappa ,p}^{\delta ,\eta }(\gamma _{N(t-1)}^{\beta })+b_{\delta }\\
&\leq &C_{q,\kappa,p }\left\Vert \varphi \right\Vert _{\infty }\times
A_{q,\kappa ,p}^{\delta ,\eta }(\gamma _{N(t)}^{\beta })+b_{\delta }
\end{eqnarray*}%
with%
\begin{eqnarray*}
b_{\delta } &=&\left\Vert \overline{P}_{t-1,t}(\overline{P}%
_{t,r}-P_{t,r})\varphi _{\delta }\right\Vert _{\infty }\leq \left\Vert (%
\overline{P}_{t,r}-P_{t,r})\varphi _{\delta }\right\Vert _{\infty }\leq \\
&\leq &C\left\Vert \nabla \varphi _{\delta }\right\Vert _{k_{0,\infty
}}\gamma _{N(r)}^{\alpha }\leq \frac{C}{\delta ^{1+k_{0}}}\left\Vert \varphi
\right\Vert _{\infty }\gamma _{N(t)}^{\alpha }.
\end{eqnarray*}%
Here we used (\ref{app8}) and $\gamma _{N(r)}\leq \gamma _{N(t)}$. We
conclude that 
\[
\left\Vert \overline{P}_{s,t}(\overline{P}_{t,r}-P_{t,r})\varphi \right\Vert
_{\infty }\leq C_{q,\kappa,p }\left\Vert \varphi \right\Vert _{\infty
}\times (A_{q,\kappa ,p}^{\delta ,\eta }(\gamma _{N(t)}^{\beta })+\frac{1}{%
\delta ^{1+k_{0}}}\gamma _{N(t)}^{\alpha }).
\]

\textbf{Optimization} For some fixed $\alpha ,\beta ,p,k_{0},\varepsilon $,
we optimize over $\delta ,\eta ,\kappa ,q$. Let $\Delta =\gamma
_{N(t)}^{\beta }.$ First we choose $\eta =\Delta ^{\frac{p}{p+\kappa }}$ so
that $\eta ^{-p}\Delta ^{p}=\eta ^{\kappa }.$ Then%
\[
A_{q,\kappa ,p}^{\delta ,\eta }(\gamma _{N(t)}^{\beta })=\frac{\delta ^{q}}{%
\Delta ^{\frac{2pq}{p+\kappa }}}+2\Delta ^{\frac{p\kappa }{p+\kappa }}.
\]%
Take now $\delta =\Delta ^{\frac{3p}{p+\kappa }}$ so that 
\[
A_{q,\kappa ,p}^{\delta ,\eta }(\gamma _{N(t)}^{\beta })=\Delta ^{\frac{pq}{%
p+\kappa }}+2\Delta ^{\frac{p\kappa }{p+\kappa }}.
\]%
With this choice 
\begin{eqnarray*}
A_{q,\kappa ,p}^{\delta ,\eta }(\gamma _{N(t)}^{\beta })+\frac{\gamma
_{N(t)}^{\alpha }}{\delta ^{1+k_{0}}} &=&\Delta ^{\frac{pq}{p+\kappa }%
}+2\Delta ^{\frac{p\kappa }{p+\kappa }}+\Delta ^{-\frac{3p(1+k_{0})}{%
p+\kappa }}\gamma _{N(t)}^{\alpha } \\
&=&\gamma _{N(t)}^{\frac{pq\beta }{p+\kappa }}+2\gamma _{N(t)}^{\frac{%
p\kappa \beta }{p+\kappa }}+\gamma _{N(t)}^{-\frac{3p(1+k_{0})\beta }{%
p+\kappa }}\times \gamma _{N(t)}^{\alpha }
\end{eqnarray*}%
We need%
\begin{eqnarray*}
i)\quad \frac{3p(1+k_{0})\beta }{p+\kappa } &<&\varepsilon , \\
ii)\quad \frac{\kappa }{p+\kappa } &\geq &1-\varepsilon \\
iii)\quad \frac{q}{p+\kappa } &\geq &1-\varepsilon .
\end{eqnarray*}%
We first choose $\kappa (\varepsilon )$ such that $i)$ and $ii)$ hold true.
Then we choose $q(\varepsilon )$ such that $\frac{q(\varepsilon )}{p+\kappa
(\varepsilon )}\geq 1-\varepsilon .$ With this choice we have 
\begin{eqnarray*}
\left\Vert \overline{P}_{s,t}(\overline{P}_{t,r}-P_{t,r})\varphi \right\Vert
_{\infty } &\leq &C_{q,\kappa,p }\left\Vert \varphi \right\Vert _{\infty
}\times (A_{q,\kappa ,p}^{\delta ,\eta }(\gamma _{N(t)}^{\beta })+\gamma
_{N(t)}^{-\frac{3p(1+k_{0})\beta }{p+\kappa }}\gamma _{N(t)}^{\alpha }) \\
&\leq &C_{q(\varepsilon ),\kappa (\varepsilon ),p}^{\prime }\left\Vert
\varphi \right\Vert _{\infty }\times (\gamma _{N(t)}^{p\beta (1-\varepsilon
)}+\gamma _{N(t)}^{\alpha -\varepsilon }) \\
&\leq &C_{q(\varepsilon ),\kappa(\varepsilon ),p}^{\prime }\left\Vert
\varphi \right\Vert _{\infty }\times \gamma _{N(t)}^{((p\beta )\wedge \alpha
)-\bar{\varepsilon}},
\end{eqnarray*}%
with $\bar{\varepsilon}=p\beta \varepsilon \vee \varepsilon $. $\square $

We give now \textbf{the main result}. We recall $n_\rho$ and $n_\ast$ in (\ref{app5}) and (\ref{app}).

\begin{proposition}
Let $\beta >0$ and $p\geq 1$ be fixed. Suppose that (\ref{app3}) holds for
some $\alpha,k_0$, (\ref{app9}),(\ref{app9'}),(\ref{app10}) hold for every $k
$ and some $\rho$ with $\rho >\alpha 
\overline{\omega }$, and $R_{\overline{\mathcal{P}}}(p,\beta )$ (see (\ref{app11})) holds. For
every $\varepsilon >0$ and every measurable and bounded function $\varphi $, for $n$ large enough such that $N(\Gamma_n-3)>n_\ast$ and $N(\Gamma_n-2)>n_\rho+1$, we have
\begin{equation}
\left\Vert (\overline{P}_{0,\Gamma _{n}}-P_{0,\Gamma _{n}})\varphi
\right\Vert _{\infty }\leq C_{\varepsilon }\left\Vert \varphi \right\Vert
_{\infty }\gamma _{n}^{((p\beta )\wedge \alpha )-\varepsilon }.
\label{app15}
\end{equation}%
Moreover, if an invariant probability measure $\nu $ exists, then the
invariant probability measure $\nu $ is unique and for every $x\in \mathbb{R}%
^{d}$, we have
\begin{equation}
d_{TV}(\overline{P}_{0,\Gamma _{n}}(x,.),\nu )\leq C_{\varepsilon }(\gamma
_{n}^{((p\beta )\wedge \alpha )-\varepsilon }+\int_{\mathbb{R}^d} \left\vert
x-y\right\vert d\nu (y)e^{-\rho \Gamma _{n}}).  \label{app16}
\end{equation}
\end{proposition}

\bigskip

\textbf{Proof }We fix $i<n$ such that $1< \Gamma _{i}$ and $\Gamma
_{i}+1\leq \Gamma _{n}\leq \Gamma _{i}+2$ and we write%
\begin{eqnarray*}
&&\left\Vert (\overline{P}_{0,\Gamma _{n}}-P_{0,\Gamma _{n}})\varphi
\right\Vert _{\infty } \\
&\leq &\left\Vert (\overline{P}_{0,\Gamma _{i}}\overline{P}_{\Gamma
_{i},\Gamma _{n}}-\overline{P}_{0,\Gamma _{i}}P_{\Gamma _{i},\Gamma
_{n}})\varphi \right\Vert _{\infty }+\left\Vert (\overline{P}_{0,\Gamma
_{i}}P_{\Gamma _{i},\Gamma _{n}}-P_{0,\Gamma _{i}}P_{\Gamma _{i},\Gamma
_{n}})\varphi \right\Vert _{\infty } \\
&=&:A+B.
\end{eqnarray*}%
First, since $\Gamma _{i}> 1,$ using (\ref{app12}) with $s=0,t=\Gamma
_{i} $ and $r=\Gamma _{n}$ we obtain 
\[
A\leq C_{\varepsilon }\left\Vert \varphi \right\Vert _{\infty }\times \gamma
_{i}^{((p\beta )\wedge \alpha )-\varepsilon }\leq C_{\varepsilon }\left\Vert
\varphi \right\Vert _{\infty }\times \gamma _{n}^{((p\beta )\wedge \alpha
)-\varepsilon },
\]%
where in the last inequality, we have used (\ref{app}).

Moreover, we recall that (\ref{app9'}) and (\ref{app10}) imply (\ref{app7}). So using (\ref{app8}) and the regularization property (\ref{app9})
(notice that $\Gamma _{n}-\Gamma _{i}\geq 1)$ we obtain%
\[
B\leq C\left\Vert \nabla P_{\Gamma _{i},\Gamma _{n}}\varphi \right\Vert
_{k_{0},\infty }\gamma _{i}^{\alpha }\leq C\left\Vert \varphi \right\Vert
_{\infty }\gamma _{i}^{\alpha }\leq C\left\Vert \varphi \right\Vert _{\infty
}\gamma _{n}^{\alpha },
\]%
the last inequality being obtained by (\ref{app}) (because $\Gamma
_{n}-\Gamma _{i}\leq 2)$.

Finally, in order to obtain (\ref{app16}) we use (\ref{App5}). The
uniqueness of the invariant probability measure $\nu $ comes directly from 
\textbf{Proposition 2.0.1}. $\square $

\bigskip\bigskip

\section{Abstract integration by parts framework}

Here we recall the abstract integration by parts framework in $\cite{ref2}$.

We denote $C_p^\infty(\mathbb{R}^d)$ to be the space of smooth functions
which, together with all the derivatives, have polynomial growth. We also
denote $C_p^q(\mathbb{R}^d)$ to be the space of $q-$times differentiable
functions which, together with all the derivatives, have polynomial growth.

We consider a probability space ($\Omega$,$\mathcal{F}$,$\mathbb{P}$), and a
linear subset $\mathcal{S}\subset\mathop{\bigcap}\limits_{p=1}^\infty
L^p(\Omega;\mathbb{R})$ such that for every $\phi\in C_p^\infty(\mathbb{R}%
^d) $ and every $F\in\mathcal{S}^d$, we have $\phi(F)\in\mathcal{S}$. A
typical example of $\mathcal{S}$ is the space of simple functionals, as in
the standard Malliavin calculus. Another example is the space of "Malliavin
smooth functionals", usually denoted by $\mathcal{D}_\infty$ (see $\cite%
{ref15}$).

Given a separable Hilbert space $\mathcal{H}$, we assume that we have a
derivative operator $D: \mathcal{S}\rightarrow\mathop{\bigcap}%
\limits_{p=1}^\infty L^p(\Omega;\mathcal{H})$ which is a linear application
which satisfies

$a)$ 
\begin{eqnarray}
D_hF:=\langle DF,h\rangle_\mathcal{H}\in\mathcal{S},\ for\ any\ h\in\mathcal{%
H},  \label{0.a}
\end{eqnarray}

$b)$ $\underline {Chain\ Rule}$: For every $\phi\in C_p^1(\mathbb{R}^d)$ and 
$F= (F_1,\cdots,F_d)\in \mathcal{S}^d$, we have

\begin{align}
D\phi(F)=\sum_{i=1}^d\partial_i\phi(F)DF_i ,  \label{0.00}
\end{align}

Since $D_hF\in\mathcal{S}$, we may define by iteration the derivative
operator of higher order $D^q:\mathcal{S}\rightarrow\mathop{\bigcap}%
\limits_{p=1}^\infty L^p(\Omega;\mathcal{H}^{\otimes q})$ which verifies $%
\langle D^qF,\otimes_{i=1}^q h_i\rangle_{\mathcal{H}^{\otimes
q}}=D_{h_q}D_{h_{q-1}}\cdots D_{h_1}F$. We also denote $D_{h_1,%
\cdots,h_q}^qF:=\langle D^qF,\otimes_{i=1}^q h_i\rangle_{\mathcal{H}%
^{\otimes q}}$,\quad for any $h_1,\cdots,h_q\in\mathcal{H}$. Then, $%
D_{h_1,\cdots,h_q}^qF=$ $D_{h_q}D_{h_1,\cdots,h_{q-1}}^{q-1}F$ ($q\geq2$).

We notice that since $\mathcal{H}$ is separable, there exists a countable
orthonormal base $(e_i)_{i\in\mathbb{N}}$. We denote 
\[
D_iF=D_{e_i}F=\langle DF,e_i\rangle_{\mathcal{H}}.
\]
Then 
\[
DF=\sum_{i=1}^{\infty}D_iF\times e_i\quad \text{and}\quad
D^qF=\sum_{i_1,\cdots,i_q}D_{i_1,\cdots,i_q}F\times\otimes_{j=1}^qe_j.
\]

For $F=(F_1,\cdots,F_d)\in\mathcal{S}^d$, we associate the Malliavin
covariance matrix 
\begin{eqnarray}
\sigma_F=(\sigma_F^{i,j})_{i,j=1,\cdots,d},\quad \text{with}\quad
\sigma_F^{i,j}=\langle DF_i,DF_j\rangle_\mathcal{H}.  \label{Mcov}
\end{eqnarray}
And we denote 
\begin{eqnarray}
\Sigma_p(F)=\mathbb{E}(1/ \det\sigma_F)^p.  \label{sigma}
\end{eqnarray}
We say that the covariance matrix of $F$ is non-degenerated if $%
\Sigma_p(F)<\infty$, $\forall p\geq1$.

We also assume that we have an Ornstein-Uhlenbeck operator $L:\mathcal{S}%
\rightarrow\mathcal{S}$ which is a linear operator satisfying

$a)$ $\underline {Duality}$: For every $F,G\in\mathcal{S}$,

\begin{align}
\mathbb{E}\langle DF,DG\rangle_\mathcal{H}=\mathbb{E}(FLG)=\mathbb{E}(GLF),
\label{0.01}
\end{align}

$b)$ $\underline {Chain\ Rule}$: For every $\phi\in C_p^2(\mathbb{R}^d)$ and 
$F= (F_1,\cdots,F_d)\in \mathcal{S}^d$, we have

\begin{align*}
L\phi(F)=\sum_{i=1}^d\partial_i\phi(F)LF_i-\sum_{i=1}^d\sum_{j=1}^d%
\partial_i\partial_j\phi(F)\langle DF_i,DF_j\rangle_\mathcal{H}.
\end{align*}
As an immediate consequence of the duality formula (\ref{0.01}), we know
that $L: \mathcal{S}\subset L^2(\Omega)\rightarrow L^2(\Omega)$ is closable.
But it is not clear that $D$ is also closable. We have to assume this and to
check it for each particular example.

\begin{definition}
If $D^q: \mathcal{S}\subset L^2(\Omega)\rightarrow L^2(\Omega;\mathcal{H}%
^{\otimes q})$, $\forall q\geq1$, are closable, then the triplet $(\mathcal{S%
},D,L)$ is called an IbP (Integration by Parts) framework.
\end{definition}

\begin{remark}
The bilinear forms $\Gamma(F,G)=\langle DF,DG\rangle_{\mathcal{H}}$ is
called "carr\'e du champ" operator in the theory of Dirichlet form. And $%
\mathcal{E}(F,G)=\mathbb{E}(\Gamma(F,G))$ is the Dirichlet form associated
to $\Gamma.$ So our Integration by Parts framework appears as a particular
case of the $\Gamma-$calculus, presented in $\cite{ref25}$ and $\cite{ref2}$.
\end{remark}

Now, we introduce the Sobolev norms. For any $l\geq1$, $F\in\mathcal{S}$, 
\begin{eqnarray}
\left\vert F\right\vert_{1,l} &=&\sum_{q=1}^{l}\left\vert D^{q}F\right\vert_{%
\mathcal{H}^{\otimes q}},\quad \left\vert F\right\vert_{l}=\left\vert
F\right\vert+\left\vert F\right\vert_{1,l},  \label{norm}
\end{eqnarray}

We put $\vert F\vert_{0}=\vert F\vert$, $\vert F\vert_{l}=0$ for $l<0$, and $%
\vert F\vert_{1,l}=0$ for $l\leq0$. For $F=(F_1,\cdots,F_d)\in\mathcal{S}^d$%
, we set 
\begin{eqnarray}
\left\vert F\right\vert_{1,l} &=&\sum_{i=1}^{d}\left\vert
F_i\right\vert_{1,l},\quad \left\vert
F\right\vert_{l}=\sum_{i=1}^{d}\left\vert F_i\right\vert_l,  \nonumber
\end{eqnarray}

Moreover, we associate the following norms. For any $l\geq0, p\geq1$, 
\begin{eqnarray}
\left\Vert F \right\Vert _{l,p}&=&(\mathbb{E}\left\vert F
\right\vert_{l}^{p})^{1/p}, \quad \left\Vert F \right\Vert _{p}=(\mathbb{E}%
\left\vert F \right\vert^{p})^{1/p} ,  \nonumber \\
\left\Vert F\right\Vert _{L,l,p}&=&\left\Vert F\right\Vert _{l,p}+\left\Vert
LF\right\Vert _{l-2,p}.  \label{sobnorm}
\end{eqnarray}

With these notations, we have the following lemma from $\cite{ref3}$ (lemma\
8 and lemma\ 10), which is a consequence of the chain rule.

\begin{lemma}
Let $F\in\mathcal{S}^d$. For every $l\in\mathbb{N},$ if $\phi: \mathbb{R}%
^d\rightarrow\mathbb{R}$ is a $C^{l}(\mathbb{R}^d)$ function ($l-$times
differentiable function), then there is a constant $C_l$ dependent on $l$
such that 
\[
a)\quad \vert\phi(F)\vert_{1,l}\leq\vert\nabla\phi(F)\vert\vert
F\vert_{1,l}+C_l\sup_{2\leq\vert\beta\vert\leq
l}\vert\partial^\beta\phi(F)\vert\vert F\vert_{1,l-1}^{l}.
\]
If $\phi\in C^{l+2}(\mathbb{R}^d)$, then 
\[
b)\quad \vert L\phi(F)\vert_{l}\leq\vert\nabla\phi(F)\vert\vert
LF\vert_{l}+C_l\sup_{2\leq\vert\beta\vert\leq
l+2}\vert\partial^\beta\phi(F)\vert(1+\vert F\vert_{l+1}^{l+2})(1+\vert
LF\vert_{l-1}).
\]
For $l=0$, we have 
\[
c)\quad \vert L\phi(F)\vert\leq\vert\nabla\phi(F)\vert\vert
LF\vert+\sup_{\vert\beta\vert=2}\vert\partial^\beta\phi(F)\vert\vert
F\vert_{1,1}^{2}.
\]
\end{lemma}

We denote by $\mathcal{D}_{l,p}$ the closure of $\mathcal{S}$ with respect
to the norm $\left\Vert \circ \right\Vert _{L,l,p}:$ 
\begin{equation}
\mathcal{D}_{l,p}=\overline{\mathcal{S}}^{\left\Vert \circ \right\Vert
_{L,l,p}},  \label{3.9}
\end{equation}%
and 
\begin{equation}
\mathcal{D}_{\infty}=\mathop{\bigcap}\limits_{l=1}^{\infty}\mathop{\bigcap}%
\limits_{p=1}^{\infty}\mathcal{D}_{l,p},\quad \mathcal{H}_{l}=\mathcal{D}%
_{l,2}.  \label{Dinf}
\end{equation}

For an IbP framework $(\mathcal{S},D,L)$, we now extend the operators from $%
\mathcal{S}$ to $\mathcal{D}_{\infty}$. For $F\in \mathcal{D}_{\infty}$, $%
p\geq2$, there exists a sequence $F_{n}\in \mathcal{S}$ such that $%
\left\Vert F-F_{n}\right\Vert _{p}\rightarrow 0$, $\left\Vert
F_{m}-F_{n}\right\Vert _{q,p}\rightarrow 0$ and $\left\Vert
LF_{m}-LF_{n}\right\Vert _{q-2,p}\rightarrow 0$. Since $D^{q} $ and $L$ are
closable, we can define 
\begin{equation}
D^{q}F=\lim_{n\rightarrow\infty}D^{q}F_{n}\quad in\quad L^p(\Omega;\mathcal{H%
}^{\otimes q}),\quad LF=\lim_{n\rightarrow\infty}LF_{n}\quad in\quad
L^p(\Omega).  \label{3.10}
\end{equation}%
We still associate the same norms and covariance matrix introduced above for 
$F\in\mathcal{D}_\infty$.

\begin{lemma}
The triplet $(\mathcal{D}_\infty,D,L)$ is an IbP framework.
\end{lemma}

\begin{proof}
The proof is standard and we refer to the lemma 3.1 in $\cite{ref4}$ for
details.
\end{proof}

The following lemma is useful in order to control the Sobolev norms and
covariance matrices when passing to the limit.

\begin{lemma}
\textbf{(A)}\quad We fix $p\geq 2,l\geq2.$ Let $F\in L^1(\Omega;\mathbb{R}%
^d) $ and let $F_{n}\in \mathcal{S}^d,n\in\mathbb{N}$ such that 
\begin{eqnarray*}
i)\quad \mathbb{E}\left\vert F_{n}-F\right\vert &\rightarrow &0, \\
ii)\quad \sup_{n}\left\Vert F_{n}\right\Vert _{L,l,p} &\leq &K_{l,p}<\infty.
\end{eqnarray*}%
Then for every $1\leq \bar{p}<p,$ we have $F\in \mathcal{D}_{l,\bar{p}}^d$
and $\left\Vert F\right\Vert _{L,l,\bar{p}}\leq K_{l,\bar{p}}$ . Moreover,
there exists a convex combination 
\[
G_{n}=\sum_{i=n}^{m_{n}}\gamma _{i}^{n}\times F_{i}\in \mathcal{S}^d,
\]%
with $\gamma _{i}^{n}\geq 0,i=n,....,m_{n}$ and $\sum\limits_{i=n}^{m_{n}}%
\gamma _{i}^{n}=1$, such that 
\[
\left\Vert G_{n}-F\right\Vert _{L,l,2}\rightarrow 0.
\]
\textbf{(B)}\quad For $F\in\mathcal{D}_\infty^d$, we denote 
\[
\lambda(F)=\inf_{\vert\zeta\vert=1}\langle \sigma_F\zeta,\zeta\rangle
\]
the lowest eigenvalue of the covariance matrix $\sigma_F$. We consider some $%
F$ and $F_n$ which verify $i), ii)$ in \textbf{(A)}. We also suppose that 
\[
iii)\quad(DF_n)_{n\in\mathbb{N}}\text{ is a Cauchy sequence in } L^2(\Omega;%
\mathcal{H}),
\]
and for every $p\geq1$, 
\begin{eqnarray}
iv)\quad \sup_{n}\mathbb{E}(\lambda^{-p}(F_n))\leq Q_p<\infty.  \label{iv}
\end{eqnarray}
Then we have 
\[
\mathbb{E}(\lambda^{-p}(F))\leq Q_p<\infty,\quad\forall p\geq1.
\]
\textbf{(C)}\quad We suppose that we have $(F, \bar{F})$ and $(F_n, \bar{F}%
_n)$ which verify the hypotheses of \textbf{(A)}. If we also have 
\begin{eqnarray}
v)\quad \sup_n\Vert DF_n-D\bar{F}_n\Vert_{L^2(\Omega;\mathcal{H})}\leq\bar{%
\varepsilon},  \label{v}
\end{eqnarray}%
then%
\[
\Vert DF-D\bar{F}\Vert_{L^2(\Omega;\mathcal{H})}\leq\bar{\varepsilon}.
\]
\end{lemma}

\begin{proof}
\textbf{Proof of (A)} For the sake of the simplicity of notations, we only
prove for the one dimensional case. We recall the notations in Section 3.
The Hilbert space $\mathcal{H}_{l}=\mathcal{D}_{l,2}$ equipped with the
scalar product 
\begin{eqnarray*}
\left\langle U,V\right\rangle _{L,l,2}&:=&\sum_{q=1}^{l}\mathbb{E} \langle
D^qU, D^qV\rangle_{\mathcal{H}^{\otimes q}}+\mathbb{E}( UV) \\
&+&\sum_{q=1}^{l-2}\mathbb{E} \langle D^qLU, D^qLV\rangle_{\mathcal{H}
^{\otimes q}}+\mathbb{E}( LU\times LV)
\end{eqnarray*}
is the space of the functionals which are $l-$times differentiable in $L^{2}$
sense. By $ii)$, for $p\geq2$, $\left\Vert F_{n}\right\Vert _{L,l,2}\leq
\left\Vert F_{n}\right\Vert _{L,l,p}\leq K_{l,p}$. Then, applying Banach
Alaoglu theorem, there exists $G\in\mathcal{H}_l$ and a subsequence (we
still denote it by $n$), such that $F_n\rightarrow G$ weakly in the Hilbert
space $\mathcal{H}_l$. This means that for every $Q\in\mathcal{H}_l$, $%
\langle F_n,Q\rangle_{L,l,2}\rightarrow\langle G,Q\rangle_{L,l,2}$.
Therefore, by Mazur theorem, we can construct some convex combination 
\[
G_{n}=\sum_{i=n}^{m_{n}}\gamma _{i}^{n}\times F_{i}\in \mathcal{S}
\]
with $\gamma _{i}^{n}\geq 0,i=n,....,m_{n}$ and $\sum\limits_{i=n}^{m_{n}}
\gamma _{i}^{n}=1$, such that 
\[
\left\Vert G_{n}-G\right\Vert _{L,l,2}\rightarrow 0.
\]
In particular we have 
\[
\mathbb{E}\left\vert G_{n}-G\right\vert \leq \left\Vert G_{n}-G\right\Vert
_{L,l,2}\rightarrow 0.
\]
Also, we notice that by i), 
\[
\mathbb{E}\left\vert G_{n}-F\right\vert \leq \sum_{i=n}^{m_{n}}\gamma
_{i}^{n}\times \mathbb{E}\left\vert F_{i}-F\right\vert \rightarrow 0.
\]
So we conclude that $F=G\in \mathcal{H}_{l}.$ We also have 
\[
\Vert G_n\Vert_{L,l,p}\leq\sum_{i=n}^{m_n}\gamma_i^n\Vert
F_i\Vert_{L,l,p}\leq K_{l,p}.
\]
Then a standard argument gives, for every $\bar{p}\in[1,p)$, 
\[
\left\Vert F\right\Vert _{L,l,\bar{p}}\leq K_{l, \bar{p}}.
\]

\bigskip\bigskip

\textbf{Proof of (B)} We consider for a moment some general $F,G\in\mathcal{%
D }_{\infty}^d$. Notice that $\langle
\sigma(F)\zeta,\zeta\rangle=\vert\langle DF,\zeta\rangle\vert_{ \mathcal{H}%
}^2 , $ so $\lambda(F)=\inf_{\vert\zeta\vert=1}\vert\langle
DF,\zeta\rangle\vert_{ \mathcal{H}}^2. $ It is easy to check that 
\begin{eqnarray}
\vert\sqrt{\lambda(F)}-\sqrt{\lambda(G)}\vert\leq\vert D(F-G)\vert_{\mathcal{%
\ H}}.  \label{beautiful}
\end{eqnarray}
We now come back to our framework. Recalling that $G_{n}=\sum%
\limits_{i=n}^{m_{n}} \gamma _{i}^{n}\times F_{i}$, we observe that 
\[
\Vert DG_n-DF_n\Vert_{L^2(\Omega;\mathcal{H})}\leq \sum_{i=n}^{m_{n}}\gamma
_{i}^{n}\Vert DF_{i}-DF_n\Vert_{L^2(\Omega;\mathcal{H})}\rightarrow0.
\]
Here we use the fact that $(DF_n)_{n\in\mathbb{N}}$ is a Cauchy sequence in $%
L^2(\Omega;\mathcal{H})$. Meanwhile, we know from \textbf{(A)} that $\Vert
DG_n-DF\Vert_{L^2(\Omega;\mathcal{H})}\rightarrow0. $ So we conclude that $%
\Vert DF-DF_n\Vert_{L^2(\Omega;\mathcal{H} )}\rightarrow0 $. Thus, by (\ref%
{beautiful}), $\mathbb{E}\vert\sqrt{ \lambda(F)}-\sqrt{\lambda(F_n)}%
\vert\rightarrow0.$ This gives that there exists a subsequence (also denote
by $n$) such that $\sqrt{\lambda(F_n)}$ converges to $\sqrt{\lambda(F)}$
almost surely, and consequently $\vert\lambda(F_n)\vert^{-p}$ converges to $%
\vert\lambda(F)\vert^{-p}$ almost surely. Since we have (\ref{iv}), $%
(\vert\lambda(F_n)\vert^{-p})_{n\in \mathbb{N}}$ is uniformly integrable. It
follows that 
\begin{eqnarray*}
\mathbb{E}(\vert \lambda(F)\vert^{-p})=\lim_{n\rightarrow\infty}\mathbb{E}
(\vert \lambda(F_n)\vert^{-p})\leq Q_p.
\end{eqnarray*}

\bigskip\bigskip

\textbf{Proof of (C)} Since the couples $(F, \bar{F})$ and $(F_n, \bar{F}_n)$
verify the hypotheses of \textbf{(A)}, we know by \textbf{(A)} that we may
find a convex combination such that 
\[
\overline{\lim}_{n\rightarrow\infty}\Vert\sum\limits_{i=n}^{m_{n}}\gamma
_{i}^{n}(D F_{i},D \bar{F}_{i})-(DF,D\bar{F})\Vert_{L^2(\Omega;\mathcal{H}
)}=0.
\]
Then it follows by (\ref{v}) that 
\begin{eqnarray*}
\Vert DF-D\bar{F}\Vert_{L^2(\Omega;\mathcal{H})}&\leq& \overline{\lim}
_{n\rightarrow\infty}\Vert\sum\limits_{i=n}^{m_{n}}\gamma _{i}^{n}(D F_{i}-D 
\bar{F}_i)\Vert_{L^2(\Omega;\mathcal{H})} \\
&\leq& \overline{\lim}_{n\rightarrow\infty}\sum\limits_{i=n}^{m_{n}}\gamma
_{i}^{n}\Vert D F_{i}-D\bar{F}_i\Vert_{L^2(\Omega;\mathcal{H})} \\
&\leq&\bar{\varepsilon}.
\end{eqnarray*}
\end{proof}

\bigskip

\subsection{Main consequences}

We will use the abstract framework presented above for the IbP framework $(%
\mathcal{D}_\infty,D,L)$, with $D$ and $L$ defined in (\ref{3.10}). We
recall the notations $\Vert F\Vert_{L,l,p}$ in (\ref{sobnorm}), $\Sigma_p(F)$
in (\ref{sigma}) and $\sigma_F$ in (\ref{Mcov}). For any $\eta>0$, we take $%
\Upsilon_\eta(x):(0,\infty)\rightarrow\mathbb{R}$ to be a function of class $%
C^\infty_b$ such that 
\[
\mathbbm{1}_{[\frac{\eta}{2},\infty)}\leq\Upsilon_\eta\leq\mathbbm{1}%
_{[\eta,\infty)}.
\]%
We remark that $\sigma_F$ is invertible on the set $\{\Upsilon_\eta(\det
\sigma_F)>0\}.$ We give the following lemma, which is stated in lemma 2.4 of 
$\cite{ref2}$ and is proved in the Appendix of $\cite{ref21}$, based on some
integration by parts formula.

\begin{lemma}
Let $F=(F_{1},\cdots,F_{d})\in \mathcal{D}_\infty^{d}$ and $G\in\mathcal{D}%
_\infty$. We fix $q\in\mathbb{N}$.

\textbf{(A)} Suppose that there exists a constant $C_q$ (dependent on $q,d$)
such that $\Vert F\Vert_{L,q+2,8dq}+\Sigma_{4q}(F)+\Vert G\Vert_{q,4}\leq
C_q $. Then for any multi-index $\beta$ with $\vert\beta\vert= q$ and any
function $f\in C_b^q(\mathbb{R}^d)$, 
\begin{eqnarray}
(\mathbf{B}_q)\quad\vert\mathbb{E}(\partial^\beta f(F)G)\vert\leq C_q\Vert
f\Vert_\infty,\quad\forall \vert\beta\vert= q.  \label{IBP}
\end{eqnarray}

\textbf{(B)} Suppose that there exists a constant $C_q^\prime$ (dependent on 
$q,d$) such that $\Vert F\Vert_{L,q+2,(4d+1)q}+\Vert G\Vert_{q,4}\leq
C_q^\prime$. Then for any $\eta>0$, any multi-index $\beta$ with $%
\vert\beta\vert= q$ and any function $f\in C_b^q(\mathbb{R}^d)$, 
\[
(\mathbf{B}_q^\prime)\quad\vert\mathbb{E}(\partial^\beta
f(F)\Upsilon_\eta(\det \sigma_F)G)\vert\leq C_q^\prime\Vert
f\Vert_\infty\times\frac{1}{\eta^{2q}},\quad\forall \vert\beta\vert= q.
\]
\end{lemma}

\begin{remark}
In \textbf{(A)}, we assume the non-degeneracy condition for $F$, so we can
give the estimate based on the standard integration by parts formula. In 
\textbf{(B)}, we no longer suppose non-degeneracy condition for $F$, so we
can only obtain an estimate based on a localized form of integration by
parts formula.
\end{remark}

\begin{remark}
If the property ($\mathbf{B}_q$) (respectively ($\mathbf{B}_q^\prime$))
holds for a random variable $F$, then it also holds for $F+x$ for every $x$
in $\mathbb{R}^d$, with the same constant $C_q$ (respectively $C_q^\prime$).
In order to see this, given a test function $f$, one defines $f_x(y)=f(x+y)$
so that $f(F+x)=f_x(F)$. And one notice that the infinite norm of $f_x$ is
the same as the infinite norm of $f$.
\end{remark}

\bigskip

We give now a regularization lemma which plays a crucial role in our paper.
We consider the $d-$dimensional super kernel $\phi_\delta$ in (\ref%
{SuperKernel}) and (\ref{superkernel}) and we denote%
\[
f_{\delta }(x)=f\ast \phi _{\delta }(x)=\int_{\mathbb{R}^{d}}f(y)\phi
_{\delta }(x-y)dy.
\]%
Then we have the following regularization lemma.

\begin{lemma}
We fix some $q,d\in\mathbb{N}$ and $\kappa,p\geq1$. We suppose that $F\in%
\mathcal{D}_\infty^d$ such that $\Vert F\Vert_{L,q+2,(4d+1)q}<\infty$. We
also consider an auxiliary random variable $Q\in\mathcal{D}_\infty^d$ such
that $\Sigma_{\kappa}(Q)<\infty$. Then there exists a constant $C$ depending
on $p,q,\kappa$ and $d$ (but not on $Q$) such that for any $\eta>0$ and $%
\delta>0$, for any function $f\in C_b^{q}(\mathbb{R}^d)$, we have 
\begin{equation}
\left\vert \mathbb{E}(f(F))-\mathbb{E}(f_{\delta }(F))\right\vert \leq
C\left\Vert f\right\Vert _{\infty }\times (\frac{\delta ^{q}}{\eta^{2q}}%
+\eta^{-p}{\mathbb{E}(\vert\det\sigma_F-\det\sigma_Q\vert}^p)+\eta^\kappa%
\mathbb{E}(\vert\det\sigma_Q\vert^{-\kappa})).  \label{1*}
\end{equation}
\end{lemma}

\begin{remark}
We remark that we do not assume the non-degeneracy condition for $F$, but we
need to assume that we have another random variable $Q$ which is
non-degenerated such that $\det\sigma_Q$ is close to $\det \sigma_F$. Then
we obtain the regularization lemma (\ref{1*}). The regularization lemma here
is originally from the paper $\cite{ref2}$.
\end{remark}

\begin{remark}
If the property (\ref{1*}) holds for a random variable $F$, then it also
holds for $F+x$ for every $x$ in $\mathbb{R}^d$, with the same constant $C$.
\end{remark}

\begin{proof}
We denote%
\[
R_{q}(\delta ,x)=\frac{1}{q!}\sum_{\left\vert \alpha \right\vert
=q}\int_{0}^{1}d\lambda (1-\lambda )^{q}\int_{\mathbb{R}^{d}}dy\phi _{\delta
}(y)y^{\alpha }\partial ^{\alpha }f(x+\lambda y)
\]%
with $y^{\alpha }=\prod_{i=1}^{q}y_{\alpha _{i}}$ for $\alpha =(\alpha
_{1},...,\alpha _{q}).$ Notice that if $F$ satisfies $(\mathbf{B}%
^\prime_{q}) $ with $G=1$, then 
\begin{equation}
\left\vert \mathbb{E}(R_{q}(\delta ,F)\Upsilon_\eta(\det
\sigma_F))\right\vert \leq C^\prime_{q }\frac{\left\Vert f\right\Vert
_{\infty }}{\eta^{2q}}\int_{\mathbb{R}^{d}}dy\phi _{\delta }(y)\left\vert
y\right\vert ^{q}=C^\prime_{q }\int_{\mathbb{R}^{d}}\phi (y)\left\vert
y\right\vert ^{q}dy\left\Vert f\right\Vert _{\infty }\frac{\delta ^{q}}{%
\eta^{2q}}.  \label{5*}
\end{equation}

We use a development in Taylor series of order $q$ in order to get%
\begin{eqnarray*}
\mathbb{E}(f(F)\Upsilon_\eta(\det \sigma_F))-\mathbb{E}(f_{\delta
}(F)\Upsilon_\eta(\det \sigma_F)) &=&\mathbb{E}(\int_{\mathbb{R}^{d}}dy\phi
_{\delta }(y)(f(F+y)-f(y))\Upsilon_\eta(\det \sigma_F)) \\
&=&\mathbb{E}(R_{q}(\delta ,F)\Upsilon_\eta(\det \sigma_F)).
\end{eqnarray*}%
Here we have used the property of a super kernel: $\int_{\mathbb{R}%
^{d}}y^\beta\phi(y)dy=0,\ \forall\vert\beta\vert\leq q$. Using (\ref{5*}),
we have 
\begin{eqnarray}
\vert\mathbb{E}(f(F)\Upsilon_\eta(\det \sigma_F))-\mathbb{E}(f_{\delta
}(F)\Upsilon_\eta(\det \sigma_F))\vert \leq C\left\Vert f\right\Vert
_{\infty }\frac{\delta ^{q}}{\eta^{2q}}.  \label{5**}
\end{eqnarray}
Following the idea from $\cite{ref52}$ p14, we denote 
\[
R=\frac{\det \sigma _F-\det \sigma_Q}{\det \sigma_Q}.
\]%
For an arbitrary $\eta$, we write 
\begin{eqnarray}
\mathbb{P}(\det \sigma_F<\eta)\leq \mathbb{P}(\det \sigma_F<\eta, \vert
R\vert<\frac{1}{4})+\mathbb{P}(\vert R\vert\geq\frac{1}{4}).  \label{BT}
\end{eqnarray}
When $\vert R\vert<\frac{1}{4}$, $\vert\det\sigma_F-\det\sigma_Q\vert<\frac{1%
}{4}\det\sigma_Q$. This implies that $\det\sigma_F>\frac{1}{2}\det\sigma_Q$.
Recalling that $Q$ is non-degenerated and using Markov inequality, for every 
$\kappa\in\mathbb{N}$, it follows that 
\begin{eqnarray}
\mathbb{P}(\det \sigma_F<\eta, \vert R\vert<\frac{1}{4})\leq\mathbb{P}%
(\det\sigma_Q<2\eta)\leq 2^\kappa\eta^\kappa\mathbb{E}(\vert\det\sigma_Q%
\vert^{-\kappa}).  \label{BT1}
\end{eqnarray}%
For any $\eta>0$, $\kappa\in\mathbb{N}$, we write 
\begin{eqnarray}
\mathbb{P}(\vert R\vert\geq\frac{1}{4})&=&\mathbb{P}(\vert\det\sigma_F-\det%
\sigma_Q\vert\geq\frac{1}{4}\det\sigma_Q)  \nonumber \\
&\leq&\mathbb{P}(\det\sigma_Q\leq\eta)+\mathbb{P}(\vert\det\sigma_F-\det%
\sigma_Q\vert>\frac{1}{4}\eta)  \nonumber \\
&\leq&C(\eta^\kappa\mathbb{E}(\vert\det\sigma_Q\vert^{-\kappa})+\eta^{-p}%
\mathbb{E}(\vert \det\sigma_F-\det\sigma_Q\vert^p)).  \label{BT2}
\end{eqnarray}%
So we conclude that 
\begin{eqnarray}
\mathbb{P}(\det \sigma_F<\eta)&\leq&C(\eta^\kappa\mathbb{E}%
(\vert\det\sigma_Q\vert^{-\kappa})+\eta^{-p}\mathbb{E}(\vert
\det\sigma_F-\det\sigma_Q\vert^p)).  \label{BT0}
\end{eqnarray}
Then we have 
\begin{eqnarray}
\vert\mathbb{E}((1-\Upsilon_\eta(\det \sigma_F))f(F))\vert\leq\Vert
f\Vert_\infty\mathbb{P}(\det \sigma_F<\eta)\leq C\Vert
f\Vert_\infty(\eta^\kappa\mathbb{E}(\vert\det\sigma_Q\vert^{-\kappa})+%
\eta^{-p}\mathbb{E}(\vert \det\sigma_F-\det\sigma_Q\vert^p)).  \label{5***}
\end{eqnarray}%
Similarly, we also have 
\begin{eqnarray}
\vert\mathbb{E}((1-\Upsilon_\eta(\det \sigma_F))f_\delta(F))\vert\leq C\Vert
f\Vert_\infty(\eta^\kappa\mathbb{E}(\vert\det\sigma_Q\vert^{-\kappa})+%
\eta^{-p}\mathbb{E}(\vert \det\sigma_F-\det\sigma_Q\vert^p)).  \label{5****}
\end{eqnarray}%
We conclude by combining (\ref{5**}), (\ref{5***}) and (\ref{5****}).
\end{proof}

\bigskip\bigskip

\section{Application for jump equations}

\subsection{Basic notations and the main equation}

To begin, we introduce some notations which will be used in the following
sections. For a multi-index $\beta$, we denote $\vert\beta\vert$ to be the
length of $\beta$. We denote $C_b^{l}(\mathbb{R}^d)$ the space of $l-$times
differential and bounded functions on $\mathbb{R}^d$ with bounded
derivatives up to order $l$, and $\left\Vert f \right\Vert _{l,\infty
}:=\sum\limits_{ \vert\beta\vert\leq l}\left\Vert \partial^\beta
f\right\Vert _{\infty }$ for a function $f\in C_b^{l}(\mathbb{R}^d)$. We
also denote $\mathcal{P}_l(\mathbb{R}^d)$ the space of all probability
measures on $\mathbb{R}^d$ with finite $l-$moment. For $\rho_1,\rho_2\in%
\mathcal{P}_1(\mathbb{R}^d)$, we define the Wasserstein distance $W_1$ by 
\begin{eqnarray}
W_1(\rho_1,\rho_2)=\sup\limits_{Lip(f)\leq1}\big\vert\int_{\mathbb{R}%
^d}f(x)\rho_1(dx)-\int_{\mathbb{R}^d}f(x)\rho_2(dx)\big\vert,  \label{W1}
\end{eqnarray}%
with $Lip(f):=\sup\limits_{x\neq y}\frac{\vert f(x)-f(y)\vert}{\vert x-y\vert%
}$ the Lipschitz constant of $f$, and we define the total variation distance 
$d_{TV}$ by 
\begin{eqnarray}
d_{TV}(\rho_1,\rho_2)=\sup\limits_{\Vert f\Vert_{\infty}\leq1}\big\vert\int_{%
\mathbb{R}^d}f(x)\rho_1(dx)-\int_{\mathbb{R}^d}f(x)\rho_2(dx)\big\vert.
\label{dTV}
\end{eqnarray}
For $F,G\in L^1(\Omega)$, we also denote $W_1(F,G)=W_1(\mathcal{L}(F),%
\mathcal{L}(G))$ and $d_{TV}(F,G)=d_{TV}(\mathcal{L}(F),\mathcal{L}(G))$,
with $\mathcal{L}(F)$(respectively $\mathcal{L}(G)$) the law of the random
variable $F$(respectively $G$). We refer to $\cite{ref22}$ and $\cite{ref19}$
the basic properties of these distances. In addition, along the paper, $C$
will be a constant which may change from a line to another. It may depend on
some parameters and sometimes the dependence is precised in the notation
(ex. $C_l $ is a constant depending on $l$).

In this paper, we consider the $d-$dimensional stochastic differential
equation with jumps 
\begin{eqnarray}
X_{t}&=&x+\int_{0}^{t}b(X_{r})dr +\int_{0}^{t}\int_{\mathbb{R}%
^d}c(z,X_{r-})N(dz,dr),  \label{1.1}
\end{eqnarray}
where $N(dz,dr)$ is a Poisson point measure on the state space $\mathbb{R}^d$
with intensity measure $\widehat{N}(dz,dr)=\mu(dz)dr$, $x$ is the initial
value, $\mu$ is a positive $\sigma$-finite measure on $\mathbb{R}^d$, and $b:%
\mathbb{R}^d\rightarrow\mathbb{R}^d$, $c:\mathbb{R}^d\times\mathbb{R}%
^d\rightarrow\mathbb{R}^d$.

\subsection{Hypotheses}

Here we give our hypotheses.

\textbf{Hypothesis 2.1} (\textbf{Regularity}) We assume that the function $%
x\mapsto b(x)$ is infinitely differentiable with bounded derivatives of any
orders. We also assume that the function $(z,x)\mapsto c(z,x)$ is infinitely
differentiable and for every multi-indices $\beta_1,\beta_2$, there exists a function $\bar{c}:\mathbb{R}^d\rightarrow\mathbb{R}_{+}$
depending on $\beta_1,\beta_2$ such that we have 
\begin{eqnarray}
\sup_{x\in\mathbb{R}^d}(\vert {c}(z,x)\vert+\vert
\partial_z^{\beta_2}\partial_x^{\beta_1} {c}(z,x)\vert)\leq\bar{c}(z),
\quad\forall z\in\mathbb{R}^d,  \label{CBar}
\end{eqnarray}
\begin{eqnarray}
with\quad \int_{\mathbb{R}^d}\vert \bar{c}(z)\vert^p\mu(dz):=\bar{c}%
_p<\infty,\quad\forall p\geq1.  \label{cbar}
\end{eqnarray}

\begin{remark}
We will use several times the following consequence of (\ref{cbar}) and of
Burkholder inequality (see for example the Theorem 2.11 in $\cite{ref12}$,
see also  $\cite{ref13}$): Let ${\Phi}(s,z,\omega):[0,T]\times%
\mathbb{R}^d\times\Omega\rightarrow\mathbb{R}_{+}$ and ${\varphi}(s,\omega):%
[0,T]\times\Omega\rightarrow\mathbb{R}_{+}$ be two non-negative
functions. The Burkholder inequality states that for any $p\geq2$, we have 
\begin{eqnarray}
&&\mathbb{E}\vert\int_{0}^t\int_{\mathbb{R}^d}{\Phi}(s,z,\omega)N(dz,ds)%
\vert^p  \nonumber \\
&&\leq C[\mathbb{E}(\int_{0}^t\int_{\mathbb{R}^d}\vert{\Phi}%
(s,z,\omega)\vert^2 \mu(dz)ds)^{\frac{p}{2}} +\mathbb{E}\int_{0}^t\int_{%
\mathbb{R}^d}\vert{\Phi}(s,z,\omega)\vert^p \mu(dz)ds  \nonumber \\
&&+\mathbb{E}\vert\int_{0}^t\int_{\mathbb{R}^d}\vert{\Phi}(s,z,\omega)\vert
\mu(dz)ds\vert^{p}].  \label{Burk}
\end{eqnarray}
If we have 
\[
\vert{\Phi}(s,z,\omega)\vert\leq\vert\bar{c}(z)\vert\vert{\varphi}%
(s,\omega)\vert,
\]
then for any $p\geq2$, 
\begin{eqnarray}
\mathbb{E}\Big\vert\int_{0}^t\int_{\mathbb{R}^d}{\Phi}(s,z,\omega)N(dz,ds)%
\Big\vert^p \leq C \mathbb{E}\int_{0}^t\vert{\varphi}(s,\omega)\vert^p ds,
\label{BurkN}
\end{eqnarray}
where $C$ is a constant depending on $p$, $\bar{c}_1$, $\bar{c}_2$, $\bar{c}%
_p$ and $T$.
\end{remark}

\begin{proof}
By compensating $N$ and using Burkholder inequality and (\ref{cbar}), we
have 
\begin{eqnarray}
&&\mathbb{E}\vert\int_{0}^t\int_{\mathbb{R}^d}{\Phi}(s,z,\omega)N(dz,ds)%
\vert^p  \nonumber \\
&&\leq C[\mathbb{E}(\int_{0}^t\int_{\mathbb{R}^d}\vert{\Phi}%
(s,z,\omega)\vert^2 \mu(dz)ds)^{\frac{p}{2}} +\mathbb{E}\int_{0}^t\int_{%
\mathbb{R}^d}\vert{\Phi}(s,z,\omega)\vert^p \mu(dz)ds  \nonumber \\
&&+\mathbb{E}\vert\int_{0}^t\int_{\mathbb{R}^d}\vert{\Phi}(s,z,\omega)\vert
\mu(dz)ds\vert^{p}]  \nonumber \\
&&\leq C\mathbb{E}\int_{0}^t\vert{\varphi}(s,\omega)\vert^pds.  \nonumber
\end{eqnarray}
\end{proof}

For the sake of simplicity of notations, in the following, for a constant $C$%
, we do not precise the dependence on the regularity constants of the
function $b$ and $c$ (such as $\Vert\nabla_xb\Vert_{\infty}$, $L_b$ and $%
\bar{c}_p$).

\textbf{Hypothesis 2.2} We assume that there exists a non-negative function $%
\breve{c}:\mathbb{R}^d\rightarrow\mathbb{R}_{+}$ such that $\int_{\mathbb{R}%
^d}\vert \breve{c}(z)\vert^p\mu(dz):=\breve{c}_p<\infty,\ \forall p\geq1$,
and 
\[
\left\Vert\nabla_{x}{c}(z,x)(I_d+\nabla_{x}{c}(z,x))^{-1}\right\Vert\leq%
\breve{c}(z),\quad\forall x\in\mathbb{R}^d, z\in\mathbb{R}^d,
\]%
with $I_d$ the $d-$dimensional identity matrix. To avoid overburdening
notation, since both hypotheses 2.1 and 2.2 apply, we take $\breve{c}(z)=%
\bar{c}(z)$ and $\breve{c}_p=\bar{c}_p$.

\begin{remark}
We need this hypothesis to prove the regularity of the inverse tangent flow
(see Section 5.2).
\end{remark}

\textbf{Hypothesis 2.3} (\textbf{Ellipticity}) There exists a non-negative
function $\underline{c}:\mathbb{R}^d\rightarrow\mathbb{R}_{+}$ such that for
every $x\in\mathbb{R}^d, z\in\mathbb{R}^d, \zeta\in\mathbb{R}^d$, we have 
\[
\sum_{j=1}^d\langle\partial_{z_j}{c}(z,x),\zeta\rangle^{2}\geq \underline{c}%
(z)\vert\zeta\vert^2.
\]

\begin{remark}
We notice that together with \textbf{Hypothesis 2.1}, we have $\underline{c}%
(z)\leq\vert\bar{c}(z)\vert^2,\ \forall z\in\mathbb{R}^d$.
\end{remark}

\textbf{Hypothesis 2.4}

We give some supplementary hypotheses concerning the function $\underline{c}$
and the measure $\mu$.

$a)$ We assume that 
\begin{eqnarray}
\underline{\lim}_{u\rightarrow+\infty}\frac{1}{\ln u}\overline{\mu}\{%
\underline{c}\geq \frac{1}{u}\}=\infty,  \label{undergamma}
\end{eqnarray}%
with 
\[
\overline{\mu}(dz)=\sum_{k=1}^{\infty}\mathbbm{1}_{[k-\frac{3}{4},k-\frac{1}{%
4}]}(\vert z\vert)\mu(dz).
\]
This means that $\underline{c}$ could not be too small so that we could have
enough noises to deduce the non-degeneracy of the Malliavin covariance
matrix (see Section 5.2).

\begin{remark}
If $\mu(\mathbb{R}^d)<\infty$, then $\underline{\lim}_{u\rightarrow+\infty}%
\frac{1}{\ln u}\overline{\mu}\{\underline{c}\geq \frac{1}{u}\}=0$. So (\ref%
{undergamma}) implies that $\mu(\mathbb{R}^d)=\infty$.
\end{remark}

\bigskip

$b)$ We assume that $\mu$ is absolutely continuous with respect to the
Lebesgue measure: $\mu(dz)=h(z)dz$, where $h$ is infinitely differentiable
and $\ln h$ has bounded derivatives of any order.

\begin{remark}
We need this hypothesis to construct the integration by parts framework for
the jump equations.
\end{remark}

\textbf{Hypothesis 2.5}

We give some conditions which ensure the existence and uniqueness of the
invariant measure and the "exponential Lipschitz property" (\ref{app1}).

Suppose that 
\begin{eqnarray}
i)\quad \left\langle x-y,b(x)-b(y)\right\rangle &\leq &-\overline{b}%
\left\vert x-y\right\vert ^{2}  \nonumber \\
ii)\quad \left\vert c(z,x)-c(z,y)\right\vert &\leq &\bar{c}(z)\left\vert
x-y\right\vert  \label{2.5ab}
\end{eqnarray}%
and%
\begin{equation}
iii)\quad 2\overline{b}-\int_{\mathbb{R}^d}(2\bar{c}(z)+\bar{c}%
^{2}(z))\mu(dz):=\theta >0.  \label{2.5c}
\end{equation}

\textbf{Hypothesis 2.6}

We assume that $\mathcal{P}$ is a partition with decreasing time steps: $%
\mathcal{P}=\{0=\Gamma_0<\Gamma_1<\cdots<\Gamma_{n-1}<\Gamma_n<\cdots\}$. We
denote $\gamma_n=\Gamma_n-\Gamma_{n-1},\ n\in\mathbb{N}$ and assume that $%
\gamma_n\downarrow0$. We also introduce 
\[
\overline{\omega }=\overline{\lim_{n\rightarrow \infty }}\frac{\gamma
_{n}-\gamma _{n+1}}{\gamma _{n+1}^{2}},
\]%
and assume that $\overline{\omega }<\frac{\theta}{2}$, with $\theta$ given
in (\ref{2.5c}).

\begin{remark}
A typical example is $\gamma _{n}=\frac{1}{n}$ and so $\overline{\omega }=1.$
\end{remark}

\bigskip\bigskip

\subsection{The truncated Euler scheme}

Now we construct the Euler scheme. For some technical reasons, we take a
general partition $\mathcal{P}=\{0=\Gamma_0<\Gamma_1<\cdots<\Gamma_{n-1}<%
\Gamma_n<\cdots\}$ (without assuming \textbf{Hypothesis 2.6} at this
moment). We denote $\gamma_n=\Gamma_n-\Gamma_{n-1},\ n\in\mathbb{N}$ and
denote $\vert\mathcal{P}\vert:=\max\limits_{n\in\mathbb{N}%
}(\Gamma_{n+1}-\Gamma_n)$. We assume that $\vert\mathcal{P}\vert\leq1$, and 
\[
\sum_{i=1}^{\infty }\gamma _{i}=\lim_{n\rightarrow \infty }\Gamma
_{n}=\infty .
\]%
For $\Gamma _{n}\leq t<\Gamma _{n+1}$ we denote $N(t)=n$ and $\tau
(t)=\Gamma _{n}. $ We consider the Euler scheme: 
\begin{eqnarray}
X^{\mathcal{P}}_{t}&=&x+\int_{0}^{t}b(X^{\mathcal{P}}_{\tau(r)})dr
+\int_{0}^{t}\int_{\mathbb{R}^d}{c}(z,X^{\mathcal{P}}_{\tau(r)-})N(dz,dr).
\label{1.3}
\end{eqnarray}

Since we have $\mu(\mathbb{R}^d)=\infty$ (which is a consequence of (\ref%
{undergamma})), we have infinitely many jumps. We use a truncation argument
in order to have finite numbers of jumps and obtain a representation by
means of a compound Poisson process. This is necessary in order to obtain a
scheme which may be simulated. We construct the truncated Euler scheme as
below. To begin, we give some notations.

We denote 
\begin{eqnarray}
\varepsilon_m:=\int_{\{\vert z\vert>m\}}\vert\bar{c}(z)\vert^2\mu(dz)+\vert%
\int_{\{\vert z\vert>m\}}\bar{c}(z)\mu(dz)\vert^2,\quad\forall m\in\mathbb{N}%
.  \label{epsM}
\end{eqnarray}
For every $\gamma>0$, we define the truncation function $M(\gamma)\in\mathbb{N}$ to
be the smallest integer such that 
\begin{eqnarray}
\varepsilon_{M(\gamma)}\leq \gamma^2.  \label{Mgamma}
\end{eqnarray}

For $m\in\mathbb{N}$, we denote $B_m=\{z\in\mathbb{R}^d:\vert z\vert\leq m\}$%
. For $\Gamma_{n-1}< t\leq\Gamma_{n}$, we denote $M_{{\mathcal{P}}%
}(t)=M(\gamma_{n})$. We remark that we have $\lim\limits_{\gamma%
\rightarrow0}M(\gamma)=\infty $ and for $\Gamma_{n-1}<t\leq \Gamma_n$, we
have $M_{{\mathcal{P}}}(t)=M(\gamma_n)\geq M(\vert{\mathcal{P}}%
\vert)\rightarrow\infty,\ \text{as }\vert{\mathcal{P}}\vert\rightarrow0. $
Now we discard the "big jumps" (the jumps of size $\vert z\vert>M_{{\mathcal{%
P}}}(t)$): 
\begin{eqnarray}
X^{\mathcal{P},{M_\mathcal{P}}}_t&=&x+\int_{0}^{t}b(X^{\mathcal{P},{M_%
\mathcal{P}}}_{\tau(r)})dr +\int_{0}^{t}\int_{B_{M_{{\mathcal{P}}}(r)}}{c}%
(z,X^{\mathcal{P},{M_\mathcal{P}}}_{\tau(r)-})N(dz,dr).  \label{truncation}
\end{eqnarray}

The advantage of considering $X^{\mathcal{P},{M_\mathcal{P}}}_{t}$ is that
we may represent it by means of compound Poisson processes. For $k\in\mathbb{%
N}$, we denote $I_1=B_1$, $I_k=B_{k}\backslash B_{k-1}$ for $k\geq2$ and
take $({J}_t^{k})_{t\geq0}$ a Poisson process of intensity $\mu (I_k)$. We
denote by $({T}_i^{k})_{i\in\mathbb{N}}$ the jump times of $({J}%
_t^{k})_{t\geq0}$ and we consider a sequences of independent random
variables ${Z}_i^{k}\sim\mathbbm{1}_{I_k}(z)\frac{\mu(dz)}{\mu(I_k)},k,i\in%
\mathbb{N}$. Moreover, $({J}_t^{k}) _{\substack{ t\geq0  \\ k\in\mathbb{N}}}$
and $({Z}_i^{k})_{k,i\in\mathbb{N}}$ are taken to be independent. Then we
represent the jump's part of the equation (\ref{truncation}) by compound
Poisson processes. We write 
\begin{eqnarray*}
X^{\mathcal{P},{M_\mathcal{P}}}_{t}&=&x+\int_{0}^{t}b(X^{\mathcal{P},{M_%
\mathcal{P}}}_{\tau(r)})dr +\sum_{k=1}^{\infty}\sum_{i=1}^{J^k_t}\mathbbm{1}%
_{B_{M_{\mathcal{P}}(T^k_i)}}(Z^k_i){c}({Z}_i^{k},X^{\mathcal{P},{M_\mathcal{%
P}}}_{\tau(T^k_i)-}) \\
&=&x+\int_{0}^{t}b(X^{\mathcal{P},{M_\mathcal{P}}}_{\tau(r)})dr
+\sum_{k=1}^{\infty}\sum_{i=1}^{J^k_t}\sum_{n=0}^{N(t)}\mathbbm{1}%
_{B_{M(\gamma_{n+1})}}(Z^k_i){c}({Z}_i^{k},X^{\mathcal{P},{M_\mathcal{P}}%
}_{\tau(T^k_i)-}).
\end{eqnarray*}%
Since $Z^k_i\in B_k\backslash B_{k-1}$, it follows that
$Z^k_i\in B_{M(\gamma_{n+1})}$ is equivalent to $k\leq
M(\gamma_{n+1})$. Then 
\begin{eqnarray}
X^{\mathcal{P},{M_\mathcal{P}}}_{t}&=&x+\int_{0}^{t}b(X^{\mathcal{P},{M_%
\mathcal{P}}}_{\tau(r)})dr
+\sum_{n=0}^{N(t)}\sum_{k=1}^{M(\gamma_{n+1})}\sum_{i=1}^{J^k_t}\mathbbm{1}_{\{\Gamma_n<T^k_i\leq%
\Gamma_{n+1}\wedge t\}}{c}({Z}_i^{k},X^{\mathcal{P},{M_\mathcal{P}}%
}_{\tau(T^k_i)-}).  \label{fictiveshock}
\end{eqnarray}%
We remark that the solution of the equation (\ref{fictiveshock}) can be
constructed in an explicit way.

We recall the notation $\theta$ in \textbf{Hypothesis 2.5}. We also recall $n_\rho=n_{\frac{\theta}{2}}$ in (\ref{app5}) (with $\rho=\frac{\theta}{2}$ in our case) and $n_\ast$ in (\ref{app}). We obtain the
following error estimate for $X^{\mathcal{P},{M_\mathcal{P}}}_{t}$, which represents the main result in our paper.

\begin{theorem}
Assume that \textbf{Hypothesis 2.1$\sim$2.5} hold and the partition $%
\mathcal{P}$ satisfies \textbf{Hypothesis 2.6}. Then an invariant
probability measure $\nu $ exists and is unique, and for $n>\max\{n_{\frac{\theta}{2}}+3,n_\ast+3\}$, for any $\varepsilon>0$
there exists a constant $C_\varepsilon$ such that 
\begin{equation}
d_{TV}(\mathcal{L}(X^{\mathcal{P},{M_\mathcal{P}}}_{\Gamma_n}),\nu )\leq
C_{\varepsilon }(\gamma _{n}^{1-\varepsilon }+\int_{\mathbb{R}^d} \left\vert
x-y\right\vert d\nu (y)e^{-\frac{\theta}{2} \Gamma _{n}}).
\end{equation}
\end{theorem}

The proof of this theorem will be given in Section 6 by using some Malliavin
integration by parts techniques introduced in Section 5.

In order to apply the Malliavin framework which will be presented in Section
5, we introduce additionally an auxiliary equation as follows (see (\ref{Xm}%
) below).

For $\Gamma_{n}< t\leq \Gamma_{n+1}$, we define 
\begin{eqnarray}
a^{{\mathcal{P}}}_t=(\sum\limits_{\substack{ 1\leq i\leq n }}%
\gamma_i\int_{\{\vert z\vert\geq M(\gamma_i)\}}\underline{c}%
(z)\mu(dz)+(t-\Gamma_n)\int_{\{\vert z\vert\geq M(\gamma_{n+1})\}}\underline{%
c}(z)\mu(dz))^\frac{1}{2},  \label{abPt}
\end{eqnarray}%
where $\underline c$ is given in \textbf{Hypothesis 2.3}. We notice that $%
\vert a^{\mathcal{P}}_t\vert\leq \sqrt{t\times\varepsilon_{M(\vert\mathcal{P}%
\vert)}}\leq \sqrt{t}\times\vert\mathcal{P}\vert.$

Now we cancel the big jumps in equation (\ref{1.1}) and replace them by a ($d-$dimensional) Gaussian
random variable $\Delta$ which is independent of the Poisson point measure $%
N(dz,ds)$: 
\begin{eqnarray}
{X}^{M_{{\mathcal{P}}}}_{t}&=&x+a^{{\mathcal{P}}}_t\Delta+\int_{0}^{t}b({X}%
^{M_{{\mathcal{P}}}}_{s})ds +\int_{0}^{t}\int_{B_{M_{{\mathcal{P}}}(s)}}{c}%
(z,{X}^{M_{{\mathcal{P}}}}_{s-})N(dz,ds).  \label{Xm}
\end{eqnarray}%
We remark that $\Delta$ is necessary in order to obtain the non degeneracy
of the covariance matrix (see Section 5.2 for details).

Following the same idea as above, we represent the jump's parts of the
equation (\ref{Xm}) by compound Poisson processes: 
\begin{eqnarray}  
{X}^{M_{{\mathcal{P}}}}_{t}&=&x+a^{{\mathcal{P}} }_{t}\Delta+\int_{0}^{t}b({X%
}^{M_{{\mathcal{P}}}}_{s})ds
+\sum_{n=0}^{N(t)}\sum_{k=1}^{M(\gamma_{n+1})}\sum_{i=1}^{J^k_t}\mathbbm1_{\{\Gamma_n<T^k_i\leq
\Gamma_{n+1}\wedge t\}} {c}({Z}_i^{k},{X}^{M_{{\mathcal{P}}}}_{T^k_i-}). \label{MXm}
\end{eqnarray}

We sometimes write ${X}^{{\mathcal{P}},M_{{\mathcal{P}}}}_{t}(x)$ (resp. ${X}%
^{M_{{\mathcal{P}}}}_{t}(x)$, ${X}_{t}(x)$) instead of ${X}^{{\mathcal{P}}%
,M_{{\mathcal{P}}}}_{t}$ (resp. ${X}^{M_{{\mathcal{P}}}}_{t}$, ${X}_{t}$) to
stress the dependence on the initial value $x$.

\bigskip\bigskip

\subsection{Some examples}

We give some typical examples to illustrate our main results.

\textbf{Example 1} We take $h=1$ so the measure $\mu$ is the Lebesgue
measure. We consider two types of behaviour for $c$.

\textbf{i) Exponential decay} We assume that $\vert\bar{c}%
(z)\vert^2=e^{-a_1\vert z\vert^p}$ and $\underline{c}(z)=e^{-a_2\vert
z\vert^p}$ with some constants $0<a_1\leq a_2$, $p>0$. We only check \textbf{%
Hypothesis 2.4} here. We have 
\[
\overline{\mu}\{\underline{c}>\frac{1}{u}\}= \overline{\mu}\{\vert z\vert<(%
\frac{\ln u}{a_2})^{\frac{1}{p}}\}\geq\frac{r_d}{2}(\frac{\ln (u-1)}{a_2})^{%
\frac{d}{p}},
\]%
with $r_d$ the volume of the unit ball in $\mathbb{R}^d$, so that 
\[
\frac{1}{\ln u}\overline{\mu}\{\underline{c}>\frac{1}{u}\}\geq \frac{r_d}{%
2(a_2)^{\frac{d}{p}}}\frac{(\ln (u-1))^{\frac{d}{p}}}{\ln u}.
\]
We notice that $\underline{\lim}_{u\rightarrow+\infty}\frac{1}{\ln u}%
\overline{\mu}\{\underline{c}\geq \frac{1}{u}\}=\infty$ when $0<p<d$.
Therefore, when $p\geq d$, we can say nothing; when $0<p<d$, the results in 
\textbf{Theorem 4.1} are true.

\textbf{ii) Polynomial decay} We assume that $\vert\bar{c}(z)\vert^2=\frac{%
a_1}{1+\vert z\vert^p}$ and $\underline{c}(z)=\frac{a_2}{1+\vert z\vert^p}$
for some constants $0<a_2\leq a_1$ and $p>d$. Then 
\[
\overline{\mu}\{\underline{c}>\frac{1}{u}\}= \overline{\mu}\{\vert
z\vert<(a_2u-1)^{\frac{1}{p}}\}\geq\frac{r_d}{2}(a_2(u-1)-1)^{\frac{d}{p}},
\]%
so that 
\[
\frac{1}{\ln u}\overline{\mu}\{\underline{c}>\frac{1}{u}\}\geq \frac{r_d}{2}%
\frac{(a_2(u-1)-1)^{\frac{d}{p}}}{\ln u}.
\]%
We notice that in this case, $\underline{\lim}_{u\rightarrow+\infty}\frac{1}{%
\ln u}\overline{\mu}\{\underline{c}\geq \frac{1}{u}\}=\infty$. Thus, the
results in \textbf{Theorem 4.1} hold true.

\bigskip

\textbf{Example 2} We consider the ($1-$dimensional) truncated $\alpha-$%
stable process: $X_t=X_0+\int_0^t\sigma(X_{r-})dU_r.$ Here $(U_t)_{t\geq0}$
is a (pure jump) L\'{e}vy process with intensity measure 
\[
\mathbbm{1}_{\{\vert z\vert\leq1\}}\frac{1}{\vert z\vert^{1+\alpha}}dz,\quad
0\leq \alpha<1.
\]
We assume that $\sigma\in C_b^\infty(\mathbb{R})$, $0<\underline{\sigma}%
\leq\sigma(x)\leq\bar{\sigma}$ and $-1<\underline{a}\leq\sigma^{\prime
}(x)\leq\bar{\sigma},\ \forall x\in\mathbb{R}$, for some universal constants 
$\bar{\sigma},\underline{\sigma},\underline{a}$, where $\sigma^{\prime }$ is
the differential of $\sigma$ in $x$. Then by a change of variable $z\mapsto%
\frac{1}{z}$, we come back to the setting of this paper with $%
c(r,v,z,x,\rho)=\sigma(x)\times\frac{1}{z}$ and ${\mu}(dz)=\mathbbm{1}%
_{\{\vert z\vert\geq1\}}\frac{1}{\vert z\vert^{1-\alpha}}dz$. We only check 
\textbf{Hypothesis 2.4} here. In this case, $\underline{c}(z)=\underline{%
\sigma}\times\frac{1}{\vert z\vert^4}$, then 
\[
\frac{1}{\ln u}\overline{\mu}\{\underline{c}>\frac{1}{u}\}\geq\frac{1}{\ln u}%
\int_1^{(\underline{\sigma}(u-1))^{\frac{1}{4}}}\frac{1}{\vert
z\vert^{1-\alpha}}dz=\frac{(\underline{\sigma}(u-1))^{\frac{\alpha}{4}}-1}{%
\alpha\ln u},
\]%
so that $\underline{\lim}_{u\rightarrow+\infty}\frac{1}{\ln u}\overline{\mu}%
\{\underline{c}\geq \frac{1}{u}\}=\infty$. Thus we can apply \textbf{Theorem
4.1}.

\bigskip\bigskip

\section{Malliavin framework for jump
equations}

We take time $t\in[0,3]$ throughout this section and we use the notations
from Section 4. We recall $(X_t)_{t\in[0,3]}$ in (\ref{1.1}), $(X^{\mathcal{P%
},M_{\mathcal{P}}}_t)_{t\in[0,3]}$ in (\ref{truncation}) and $(X^{M_{%
\mathcal{P}}}_t)_{t\in[0,3]}$ in (\ref{Xm}), where $\mathcal{P}%
=\{0=\Gamma_0<\Gamma_1<\cdots<\Gamma_{N(3)}\leq3 \}$ is a general partition (which is not supposed to verify \textbf{Hypothesis 2.6}).

\begin{lemma}
Suppose that \textbf{Hypothesis 2.1} holds true. Then we have the followings.

$i)$ For every $t\in[0,3]$, we have 
\[
\mathbb{E}\vert{X}^{{{{\mathcal{P}}}},M_{{\mathcal{P}} }}_t-X_t\vert%
\rightarrow0, \ \text{as } \vert{\mathcal{P}}\vert\rightarrow0;
\]

$ii)$ For every fixed $t\in[0,3]$ and every $p\geq2$, we have 
\[
\mathbb{E}\vert {X}^{M_{{\mathcal{P}}}}_t-X_t\vert^p\rightarrow0,\ \text{as }
\vert{\mathcal{P}}\vert\rightarrow0;
\]
$iii)$ For every fixed $t\in[0,3]$ and every multi-index $\beta$, we have 
\[
\mathbb{E}\vert \partial_x^\beta {X}^{M_{{\mathcal{P}} }}_t-\partial_x^\beta
X_t\vert\rightarrow0,\ \text{as }\vert{\mathcal{P}} \vert\rightarrow0.
\]
\end{lemma}

\begin{proof}
The proof of this lemma is standard and straightforward by Gronwall lemma
and Buckholder inequality. So we leave it out.
\end{proof}

Now we use Malliavin calculus for ${X}^{ {\ \mathcal{P}},M_{{\mathcal{P}}%
}}_{t}$, ${X}^{M_{{\mathcal{P}} }}_{t}$ and ${X}_t$. There are several
approaches given in $\cite{ref5}$, $\cite{ref9}$, $\cite{ref10}$, $\cite%
{ref11}$, $\cite{ref15}$, $\cite{ref17}$ and $\cite{ref18}$ for example.
Here we give a framework analogous to $\cite{ref3}$.

To begin we define a regularization function. 
\begin{eqnarray}
a(y) &=&1-\frac{1}{1-(4y-1)^{2}}\quad for\quad y\in \lbrack \tfrac{1}{4}, 
\tfrac{1}{2}),  \label{3.14} \\
\psi (y) &=&\mathbbm{1}_{\{\left\vert y\right\vert \leq \frac{1}{4}\}}+ %
\mathbbm{1}_{\{\frac{1}{4}< \left\vert y\right\vert \leq \frac{1}{2}
\}}e^{a(\left\vert y\right\vert )}.  \label{3.14*}
\end{eqnarray}
We notice that $\psi \in C_{b}^{\infty }(\mathbb{R})$ and that its support
is included in $[-\frac{1}{2},\frac{1}{2}]$. We denote 
\begin{equation}
\Psi_k(y)=\psi(\vert y\vert-(k-\tfrac{1}{2})),\ \forall k\in\mathbb{N}.
\label{3.19*}
\end{equation}
Then for any $l\in\mathbb{N}$, there exists a constant $C_l$ such that 
\begin{eqnarray}
\sup_{k\in\mathbb{N}}\Vert\Psi_k\Vert_{l,\infty}\leq C_l<\infty.
\label{psiborn}
\end{eqnarray}

We focus on ${X}^{{\mathcal{P}},M_{{\mathcal{P}}}}_{t}(x)$ and ${X}^{M_{{%
\mathcal{P}}}}_{t}(x)$ (solutions of (\ref{fictiveshock}) and (\ref{MXm}))
which are functions of random variables $T^k_i,$ $Z^k_i$ and $\Delta$ .

Now we introduce the space of simple functionals $\mathcal{S}.$ We take $%
\mathcal{G}=\sigma(T^k_i: k,i\in \mathbb{N})$\ to be the $\sigma-$algebra
associated to the noises which will not be involved in our calculus. In the
following, we will do the calculus based on $Z^k_i=(Z^k_{i,1},
\cdots,Z^k_{i,d}),k,i\in\mathbb{N}$ and $\Delta=(\Delta_1,\cdots,\Delta_d)$.
We denote by $C_{\mathcal{G},p}$ the space of the functions $f:\Omega \times 
\mathbb{R}^{m\times m^{\prime}\times d+d}\rightarrow\mathbb{R}$ such that
for each $\omega ,$ the function $(z_{1,1}^1,...,z_{m,d}^{m^{\prime}},
\delta_1,\cdots,\delta_d)\mapsto f(\omega
,z_{1,1}^1,...,z_{m,d}^{m^{\prime}},\delta_1,\cdots,\delta_d)$ belongs to $%
C_{p}^{\infty }(\mathbb{R}^{m\times m^{\prime}\times d+d})$ (the space of
smooth functions which, together with all the derivatives, have polynomial
growth), and for each $(z_{1,1}^1,...,z_{m,d}^{m^{\prime}},\delta_1,\cdots,
\delta_d)$, the function $\omega \mapsto f(\omega
,z_{1,1}^1,...,z_{m,d}^{m^{\prime}},\delta_1,\cdots,\delta_d)$ is $\mathcal{%
G }$-measurable. And we consider the weights 
\begin{eqnarray}
\xi_{i}^k=\Psi_k({Z}^k_i).\label{xi}
\end{eqnarray}
Then we define the space of simple functionals 
\[
\mathcal{S}=\{F=f (\omega ,({Z}_{i}^{k})_{\substack{ 1\leq k\leq m^{\prime} 
\\ 1\leq i\leq m}},\Delta) :f\in C_{\mathcal{G},p}, m,m^{\prime}\in \mathbb{%
N }\}.
\]

\begin{remark}
Take $m^{\prime}=\max\limits_{t\leq3}M_{\mathcal{P}}(t)$ and $m=\max\limits_{k\leq
m^{\prime}}J^k_t$. Then ${X}^{M_{ {\mathcal{P} }}}_{t}(x)$ (solution of (\ref%
{MXm})) is a function of $T^k_i,$ $Z^k_i$ and of $\Delta$, with $k\leq m^{\prime}$
and $i\leq m$. So it is a simple functional (the same for ${X}^{ {\ \mathcal{%
P}},M_{{\mathcal{P}}}}_{t}(x)$ (solution of (\ref{fictiveshock}))).
\end{remark}

On the space $\mathcal{S}$, for $t\geq1$, we define the derivative operator $%
DF=(D^ZF,D^{\Delta}F)$, where 
\begin{eqnarray}
D^Z_{(\bar{k},\bar{i},\bar{j})}F&=&\xi _{\bar{i}}^{\bar{k}}\frac{\partial f}{
\partial z^{\bar{k}}_{\bar{i},\bar{j}}}(\omega ,({Z}_{i}^{k})_{\substack{ %
1\leq k\leq m^{\prime}  \\ 1\leq i\leq m}},\Delta),\quad \bar{k},\bar{i}\in 
\mathbb{N},\bar{j}\in\{1,\cdots,d\},  \label{D} \\
D^{\Delta}_{\Tilde{j}}F&=&\frac{\partial f}{\partial \delta_{\Tilde{j}}}
(\omega ,({Z}_{i}^{k})_{\substack{ 1\leq k\leq m^{\prime}  \\ 1\leq i\leq m}}
,\Delta),\quad \Tilde{j}\in\{1,\cdots,d\}.  \nonumber
\end{eqnarray}
We regard $D^ZF$ as an element of the Hilbert space $l_{2}$\ (the space of
the sequences $u=(u_{k,i,j})_{k,i\in \mathbb{N},{j}\in\{1,\cdots,d\}}$ with $%
\left\vert u\right\vert _{l_{2}}^{2}:=\sum_{k=1}^{\infty }\sum_{i=1}^{\infty
}\sum_{j=1}^{d}\vert u_{k,i,j}\vert^{2}<\infty$) and $DF$ as an element of $%
l_{2}\times\mathbb{R}^d$, so we have 
\begin{eqnarray}
\left\langle DF,DG\right\rangle_{l_2\times\mathbb{R}^d}
=\sum_{j=1}^dD_j^{\Delta}F\times D_j^{\Delta}G+\sum_{k=1}^{\infty
}\sum_{i=1}^{\infty }\sum_{j=1}^{d }D^Z_{(k,i,j)}F\times D^Z_{(k,i,j)}G.
\label{Cov}
\end{eqnarray}
We also denote $D^1F=DF$, and we define the derivatives of order $q\in 
\mathbb{N}$ recursively: $D^{q}F:=DD^{q-1}F.$ And we denote $D^{Z,q}$
(respectively $D^{\Delta,q}$) as the derivative $D^Z$ (respectively $%
D^\Delta $) of order $q$.

We recall that $\mu(dz)=h(z)dz$ with $h\in C^\infty(\mathbb{R}^d)$ (see 
\textbf{Hypothesis 2.4 $b)$}). We  define the Ornstein-Uhlenbeck
operator $LF=L^ZF+L^{\Delta}F$ with 
\begin{eqnarray}
L^ZF&=&-\sum_{k=1}^{m^{\prime}}\sum_{i=1}^{m }\sum_{j=1}^{d
}(\partial_{z^k_{i,j}}(\xi^k_{i}D^Z_{(k,i,j)}F)+ D^Z_{(k,i,j)}F\times
D^Z_{(k,i,j)}\ln [{h(Z^k_i)}]),  \label{L} \\
L^{\Delta}F&=&\sum_{j=1}^dD^{\Delta}_jF\times\Delta_j-\sum_{j=1}^dD^{
\Delta}_jD^{\Delta}_jF.  \nonumber
\end{eqnarray}

One can check that the triplet $(\mathcal{S},D,L)$ is consistent with the
IbP framework given in Section 3.1. In particular the duality formula (\ref%
{0.01}) holds true. We refer to $\cite{ref4}$(Appendix 5.3). We say that $F$
is a "Malliavin smooth functional" if $F\in\mathcal{D}_\infty$ (with the
definition given in (\ref{Dinf})).

\bigskip

We recall ${X}^{{\mathcal{P}},M_{{\mathcal{P}}}}_{t}(x)$ in (\ref%
{fictiveshock}), ${X}^{M_{{\mathcal{P}}}}_{t}(x)$ in (\ref{MXm}) and ${X}%
_{t}(x)$ in (\ref{1.1}). We denote 
\begin{eqnarray}
{F}_t^{{\mathcal{P}},M_{{\mathcal{P}}}}(x)={X}^{ {\mathcal{P}},M_{{\mathcal{P%
}}}}_{t}(x)-x, {F}_t^{M_{ { \mathcal{P}}}}(x)={X}^{M_{{\mathcal{P}}%
}}_{t}(x)-x\ \text{and}\ {F} _t(x)={X}_{t}(x)-x  \label{Ftx}
\end{eqnarray}
In the following subsections, we will give some lemmas concerning the
Sobolev norms and the covariance matrices. We recall (see (\ref{Mcov})) that 
$\sigma_F$ denotes the covariance matrix of $F$, and recall the Sobolev
norms defined in (\ref{norm}) and (\ref{sobnorm}).

\bigskip

\subsection{Sobolev norms}

We recall the notations ${F}_t^{{\mathcal{P}},M_{{\mathcal{P}}}}(x)$, ${F}%
_t^{M_{{\mathcal{P}}}}(x)$ and ${F}_t(x)$ in (\ref{Ftx}).

\begin{lemma}
Assuming \textbf{Hypothesis 2.1} and \textbf{Hypothesis 2.4} $b)$, for all $%
p\geq1, l\geq 0$, there exists a constant $C_{l,p}$ depending on $l,p,d$,
such that for any $t\in[0,3]$, 
\[
i)\quad\sup\limits_{{\mathcal{P}}}\sup\limits_{x}(\Vert {F}^{ {\ \mathcal{P}}%
,M_{{\mathcal{P}}}}_{t}(x)\Vert_{L,l,p}+\Vert {F} _{t}^{M_{{\mathcal{P}}%
}}(x)\Vert_{L,l,p})\leq C_{l,p}.
\]
Moreover, ${F}_{t}(x)$ belongs to $\mathcal{D}_{\infty}^d$ and 
\[
ii)\quad\sup\limits_{x}\Vert {F}_{t}(x)\Vert_{L,l,p}\leq C_{l,p}.
\]
For all $p,q\geq1, l\geq 0$, there exists a constant $C_{l,p,q}$ depending
on $l,p,q,d$, such that for every multi-index $\beta$ with $%
\vert\beta\vert=q $, we also have 
\[
iii)\quad\sup\limits_{x}\Vert \partial_x^\beta({X}_{t}(x))\Vert_{l,p}\leq
C_{l,p,q}.
\]
\end{lemma}

\begin{remark}
Since $Dx=0,\ \forall x\in\mathbb{R}^d$, we also have 
\[
\quad\sup\limits_{{\mathcal{P}}}\sup\limits_{x}(\mathbb{E}\vert {X}^{ {%
\mathcal{P}},M_{{\mathcal{P}}}}_{t}(x)\vert_{1,l}^p+\mathbb{E}\vert {X}%
_{t}^{M_{{\mathcal{P}}}}(x)\vert_{1,l}^p+\mathbb{E}\vert {X}%
_{t}(x)\vert_{1,l}^p)\leq C_{l,p}.
\]
\end{remark}

\begin{proof}
We first notice that for any $l,p$, $\sup\limits_{{\mathcal{P}}%
}\sup\limits_{x}(\Vert {F}^{ {\ \mathcal{P}},M_{{\mathcal{P}}%
}}_{t}(x)\Vert_{L,l,p}+\Vert {F} _{t}^{M_{{\mathcal{P}}}}(x)\Vert_{L,l,p})%
\leq C_{l,p}$ This is a slight variant of the proof of Lemma 3.7 $i)$ in $%
\cite{ref20}$. The difference in that the truncation function $M$ is
constant in $\cite{ref20}$ while here it depends on the time. But this does
not change anything. In a similar way, for every multi-index $\beta$ with $%
\vert\beta\vert=q$, we have $\sup\limits_{{{\mathcal{P}}} }\sup\limits_x%
\Vert \partial_x^\beta({X}_{t}^{M_{{\mathcal{P}} }}(x))\Vert_{l,p}\leq
C_{l,p,q}$.

Afterwards, we consider an increasing sequence of partition ${\mathcal{P} }%
_n,\ n\in\mathbb{N}$, (${\mathcal{P}}_n\subset{\mathcal{P}}_{n+1}$ ), such
that $\vert{\mathcal{P}}_n\vert\downarrow0$. In particular, $\forall t,\ M_{ 
{\mathcal{P}_n}}(t)\uparrow\infty$. Noticing by \textbf{\ Lemma 5.1} $ii)$
that $\mathbb{E}\vert {F}^{M_{{\mathcal{P}} _n}}_t-F_t\vert\rightarrow0$ as $%
n\rightarrow0$, and applying \textbf{Lemma 3.3 (A)} with $F_n={F}_{t}^{M_{{%
\mathcal{P}}_n}}$ and $F= {F}_t$, we get that ${F}_{t}$ belongs to $\mathcal{%
D}_{\infty}^d$ and $\sup\limits_x\Vert {F}_{t}(x)\Vert_{L,l,p}\leq C_{l,p}$.

Furthermore, noticing by \textbf{Lemma 5.1} $iii)$ that $\mathbb{E}\vert
\partial_x^\beta {X}^{M_{{\mathcal{P}} _n}}_t-\partial_x^\beta
X_t\vert\rightarrow0$ as $n\rightarrow0$, and applying \textbf{Lemma 3.3 (A)}
with $F_n=\partial_x^\beta {X}_{t}^{M_{ {\mathcal{P}}_n}}$ and $%
F=\partial_x^\beta {X}_t$, we obtain that $\partial_x^\beta{X}_{t}$ belongs
to $\mathcal{D}_{\infty}^d$ and $\sup\limits_x\Vert \partial_x^\beta({X}%
_{t}(x))\Vert_{l,p}\leq C_{l,p,q}$.

\end{proof}

\bigskip

\subsection{Covariance matrix}

\begin{lemma}
Assume that \textbf{Hypothesis 2.1, 2.2, 2.3} and \textbf{2.4} hold true. We
denote the lowest eigenvalue of the Malliavin covariance matrix $\sigma_{ {X}%
^{M_{{\mathcal{P}}}}_t}$ by $\lambda_{t}^{M_{{\mathcal{P}}}}$ . Then for
every $p\geq1$, $1\leq t\leq3$, we have 
\begin{eqnarray*}
i)\quad \sup\limits_{{\mathcal{P}}}\sup\limits_{x}\mathbb{E}(1/ \det\sigma_{{%
X}_{t}^{M_{{\mathcal{P}}}}(x)})^p\leq\sup\limits_{ {\ \mathcal{P}}%
}\sup\limits_{x}\mathbb{E}(\vert \lambda_{t}^{M_{ {\mathcal{P} }%
}}\vert^{-dp})\leq C_p, \\
ii)\quad\sup\limits_{x}\mathbb{E}(1/ \det\sigma_{{X}_{t}(x)})^p\leq C_p,
\end{eqnarray*}
with $C_p$ a constant depending on $p,d$.
\end{lemma}

\begin{remark}
We recall the notations ${F}^{M_{{\mathcal{P}}}}_{t}(x)={X}^{M_{ {\mathcal{P}%
}}}_{t}(x)-x$ and ${F}_{t}(x)={X}_{t}(x)-x$.
Since $Dx=0,\ \forall x\in\mathbb{R}^d$, the above results are equivalent to 
\begin{eqnarray*}
i)\quad \sup\limits_{{\mathcal{P}}}\sup\limits_{x}\mathbb{E}(1/ \det\sigma_{ 
{F}_{t}^{M_{{\mathcal{P}}}}(x)})^p\leq\sup\limits_{ {\ \mathcal{P}}%
}\sup\limits_{x}\mathbb{E}(\vert \lambda_{t}^{M_{ {\mathcal{P} }%
}}\vert^{-dp})\leq C_p, \\
ii)\quad\sup\limits_{x}\mathbb{E}(1/ \det\sigma_{ {F}_{t}(x)})^p\leq C_p.
\end{eqnarray*}
\end{remark}

\bigskip

\textbf{Proof of $i)$} We proceed in 4 steps.

\textbf{Step 1} We notice by the definition (\ref{D}) that for any $k_0,{i_0}
\in\mathbb{N},{j}\in\{1,\cdots,d\}$, 
\begin{eqnarray}
&&D^Z_{(k_0,{i_0},{j})}{X}^{M_{{\mathcal{P}}}}_{t}=\int_{{T}_{{i_0}
}^{k_0}}^{t}\nabla_xb({X}^{M_{{\mathcal{P}}}}_{r})D^Z_{(k_0,{i_0},{j} )}{X}%
^{M_{{\mathcal{P}}}}_{r}dr  \nonumber \\
&&+\sum_{n=0}^{N(t)}\mathbbm{1}_{\{\Gamma_n<T_{{i_0}}^{k_0}\leq
\Gamma_{n+1}\wedge t\}}\mathbbm{1}_{\{1\leq k_0\leq M(\gamma_{n+1})\}}\xi_{{%
i_0}}^{k_0}\partial_{z^{k_0}_{i_0, {j}}}{c}({Z}_{{i_0}}^{k_0},{X}^{M_{{%
\mathcal{P}}}}_{{T}_{{i_0} }^{k_0}-})  \nonumber \\
&&+\sum_{n=0}^{N(t)}\sum_{k=1}^{M(\gamma_{n+1})}\sum_{i=1}^{J^k_t}\mathbbm{1}_{\{\Gamma_n\vee T^{k_0}_{i_0}<T^k_i\leq%
\Gamma_{n+1}\wedge t\}}\nabla_x{c}({Z}_{i}^{k},{X}%
^{M_{ {\mathcal{P}}}}_{{T}_{i}^{k}-})D^Z_{(k_0,{i_0},{j})}{X}^{M_{ {\ 
\mathcal{P}}}}_{{T}_{i}^{k}-}, \quad\quad  \label{cM1}
\end{eqnarray}
\begin{eqnarray}
&&D^\Delta_j{X}^{M_{{\mathcal{P}}}}_t=a^{{\mathcal{P}}}_{t} \bm{e_j}%
+\int_{0}^{t}\nabla_xb({X}^{M_{{\mathcal{P}} }}_{r})D^\Delta_{j}{X}^{M_{{%
\mathcal{P}}}}_{r}dr
+\sum_{n=0}^{N(t)}\sum_{k=1}^{M(\gamma_{n+1})}\sum_{i=1}^{J^k_t}\mathbbm{1}_{\{\Gamma_n<T^k_i\leq%
\Gamma_{n+1}\wedge t\}}\nabla_x{c}({Z}_{i}^{k},{X}^{M_{{\mathcal{P}}}}_{{T}
_{i}^{k}-})D^\Delta_j{X}^{M_{{\mathcal{P}}}}_{{T}_{i}^{k}-},\nonumber\\ 
\label{cM1*}
\end{eqnarray}
where $\bm{e_j}=(0,\cdots,0,1,0,\cdots,0)$ with value $1$ at the $j-$th
component.

Now we introduce $({Y}^{M_{{\mathcal{P}}}}_{t})_{t\geq0}$ (this is so-called
the tangent flow) which is the matrix solution of the linear equation 
\begin{eqnarray*}
{Y}^{M_{{\mathcal{P}}}}_{t}=I_d&+&\int_{0}^{t}\nabla_xb({X}^{M_{ { \mathcal{%
P}}}}_{r}){Y}^{M_{{\mathcal{P}}}}_{r}dr+\sum_{n=0}^{N(t)}
\sum_{k=1}^{M(\gamma_{n+1})}\sum_{i=1}^{J^k_t}\mathbbm{1}_{\{\Gamma_n<T^k_i\leq%
\Gamma_{n+1}\wedge t\}}\nabla_x{c}({Z}_{{\ i}}^{k},{X}^{M_{{\mathcal{P}}}}_{{T}_{{i}}^{k}-}){Y}%
^{M_{ { \mathcal{P}}}}_{{T}_{{i}}^{k}-}.
\end{eqnarray*}
And using It$\hat{o}$'s formula, the inverse matrix $\widetilde{{Y}}^{M_{ {%
\mathcal{P}}}}_{t}=({Y}^{M_{{\mathcal{P}}}}_{t})^{-1}$ verifies the equation 
\begin{eqnarray}
\widetilde{{Y}}^{M_{{\mathcal{P}}}}_{t}=I_d&-&\int_{o}^{t}\widetilde{{Y}}
^{M_{{\mathcal{P}}}}_{r}\nabla_xb({X}^{M_{{\mathcal{P}} }}_{r})dr-%
\sum_{n=0}^{N(t)}\sum_{k=1}^{M(\gamma_{n+1})}\sum_{i=1}^{J^k_t}\mathbbm{1}_{\{\Gamma_n<T^k_i\leq%
\Gamma_{n+1}\wedge t\}}\widetilde{{Y}}^{M_{{\mathcal{P}}}}_{{T}_{{i}
}^{k}-}\nabla _{x}{c}(I_d+\nabla _{x}{c})^{-1}({Z}_{{i}}^{k},{X}^{M_{ {%
\mathcal{P}}}}_{{T}_{{i}}^{k}-}).\nonumber\\
\label{inverse}
\end{eqnarray}

\begin{remark}
We notice that ${Y}^{M_{{\mathcal{P}}}}_{t}=\nabla_x({X}^{M_{{\mathcal{P}}%
}}_{t}(x))$. If instead we consider the gradient of the Euler scheme ${Y}^{{%
\mathcal{P}},M_{{\mathcal{P}}}}_{t}=\nabla_x({X}^{{\mathcal{P}},M_{{\mathcal{%
P}}}}_{t}(x))$, the matrix ${Y}^{{\mathcal{P}},M_{{\mathcal{P}}}}_{t}$ is
not invertible, and this is a specific difficulty when we deal with the
Euler scheme. This is why we have to work with ${X}^{M_{{\mathcal{P}}}}_{t}$
only.
\end{remark}

\bigskip

Applying \textbf{Hypothesis 2.1} and \textbf{Hypothesis 2.2}, one also has 
\begin{equation}
\mathbb{E}(\sup_{0<t\leq 2}(\left\Vert {Y}^{M_{{\mathcal{P}}
}}_{t}\right\Vert ^{p}+\left\Vert \widetilde{{Y}}^{M_{{\mathcal{P}}
}}_{t}\right\Vert ^{p}))\leq C_p<\infty.  \label{cM2}
\end{equation}%
The proof of (\ref{cM2}) is straightforward and we leave it out.

Then using the uniqueness of solution to the equation (\ref{cM1}) and (\ref%
{cM1*}), one obtains 
\begin{equation}
D^Z_{(k,{i},{j})}{X}^{M_{{\mathcal{P}}}}_{t}=\sum_{n=0}^{N(t)}\mathbbm{1}%
_{\{\Gamma_n<T_{{i}}^{k}\leq \Gamma_{n+1}\wedge t\}}\mathbbm{1}_{\{1\leq
k\leq M(\gamma_{n+1})\}}\xi _{{i}}^k{Y}^{M_{{\mathcal{P}}}}_{t}\widetilde{{%
Y}}^{M_{ {\mathcal{P}}}}_{T_{{i}}^{k}}\partial_{z^k_{i,{j}}}{c}({Z}_{{i}%
}^{k},{X} ^{M_{{\mathcal{P}}}}_{{T}_{{i}}^{k}-}),  \label{cM3}
\end{equation}
and $D^\Delta_j{X}^{M_{{\mathcal{P}}}}_{t}=a^{{\mathcal{P}}}_{t}{\ Y}^{M_{{%
\mathcal{P}}}}_{t}\bm{e_j}$.

We recall that we denote the lowest eigenvalue of the Malliavin covariance
matrix $\sigma_{{X}^{M_{{\mathcal{P}}}}_t}$ by $\lambda_{t}^{M_{{\mathcal{P}}%
}}$. Then we have (recalling the definitions (\ref{Mcov}) and (\ref{Cov})) 
\[
\lambda_{t}^{M_{{\mathcal{P}}}}=\inf\limits_{\vert\zeta\vert=1}\langle
\sigma_{{X}^{M_{{\mathcal{P}}}}_t}\zeta,\zeta\rangle\geq\inf
\limits_{\vert\zeta\vert=1}\sum_{n=0}^{N(t)}\sum_{k=1}^{M(\gamma_{n+1})}
\sum_{i=1}^{J^k_t}\mathbbm{1}_{\{\Gamma_n<T^k_i\leq%
\Gamma_{n+1}\wedge t\}}\sum_{j=1}^{d }\langle
D^Z_{(k,i,j)}{\ {X}^{M_{{\mathcal{P}}}}_t},\zeta\rangle^2+\inf\limits_{\vert%
\zeta \vert=1}\sum_{j=1}^{d }\langle D^\Delta_{j}{{X}^{M_{{\mathcal{P}}}}_t }%
,\zeta\rangle^2.
\]
By (\ref{cM3}), 
\begin{eqnarray*}
\lambda _{t}^{M_{{\mathcal{P}}}}
&\geq&\inf\limits_{\vert\zeta\vert=1}\sum_{n=0}^{N(t)}
\sum_{k=1}^{M(\gamma_{n+1})}\sum_{i=1}^{J^k_t}\mathbbm{1}_{\{\Gamma_n<T^k_i\leq%
\Gamma_{n+1}\wedge t\}}\sum_{j=1}^{d }\vert\xi_i^k\vert^2\langle \partial_{z^k_{i,{j}}}{c}({Z}_{{i%
}}^{k},{X}^{M_{ {\mathcal{P}}}}_{{T}_{{i}}^{k}-}),({Y}^{M_{{\mathcal{P}}%
}}_{t} \widetilde{{Y}}^{M_{{\mathcal{P}}}}_{T_{{i}}^{k}})^{\ast}\zeta%
\rangle^2\\
&+& \inf\limits_{\vert\zeta\vert=1}\sum_{j=1}^{d }\vert a^{{\mathcal{P%
}} }_{t}\vert^2\langle \bm{e_j},({Y}^{M_{{\mathcal{P}}}}_{t})^{\ast}\zeta
\rangle^2,
\end{eqnarray*}
where $Y^\ast$ denotes the transposition of a matrix $Y$.

We recall the ellipticity hypothesis (\textbf{Hypothesis 2.3}): there exists
a non-negative function $\underline{c}(z)$ such that 
\[
\sum_{j=1}^d\langle\partial_{z_j}{c}(z,x),\zeta\rangle^{2}\geq \underline{c}
(z)\vert\zeta\vert^2.
\]
So we deduce that 
\begin{eqnarray*}
\lambda_{t}^{M_{{\mathcal{P}}}}&\geq&
\inf\limits_{\vert\zeta\vert=1}(\sum_{n=0}^{N(t)}\sum_{k=1}^{M(%
\gamma_{n+1})} \sum_{i=1}^{J^k_t}\mathbbm{1}_{\{\Gamma_n<T^k_i\leq%
\Gamma_{n+1}\wedge t\}}\vert\xi_i^k\vert^2 \underline{c}({Z} _{i}^{k})\vert({Y}^{M_{{\mathcal{P}}%
}}_{t}\widetilde{{Y}}^{M_{ {\ \mathcal{P}}}}_{T_{{i}}^{k}})^{\ast}\zeta%
\vert^2)+\vert a^{{\mathcal{P}} }_{t}\vert^2\inf\limits_{\vert\zeta\vert=1}%
\vert({Y}^{M_{{\mathcal{P}} }}_{t})^{\ast}\zeta\vert^2.
\end{eqnarray*}
For every invertible matrix $A$ and every vector $y$, one has $\vert
Ay\vert\geq\frac{1}{\Vert A^{-1}\Vert}\vert y\vert$, so that 
\begin{eqnarray*}
\lambda_{t}^{M_{{\mathcal{P}}}} &\geq&
(\sum_{n=0}^{N(t)}\sum_{k=1}^{M(\gamma_{n+1})}\sum_{i=1}^{J^k_t}\mathbbm{1}_{\{\Gamma_n<T^k_i\leq%
\Gamma_{n+1}\wedge t\}} \vert\xi_i^k\vert^2\underline{c}({Z}^k_i)\Vert%
\widetilde{{Y} }^{M_{{\mathcal{P}}}}_{t}\Vert^{-2}\Vert{Y}^{M_{{\mathcal{P}}%
}}_{T_{{\ i}}^{k}}\Vert^{-2})+ \vert a^{{\mathcal{P}}}_{t}\vert^2\Vert%
\widetilde{{Y }}^{M_{{\mathcal{P}}}}_{t}\Vert^{-2} \\
&\geq& (\inf\limits_{0<t\leq 2}\Vert\widetilde{{Y}}^{M_{{\mathcal{P}}
}}_{t}\Vert^{-2}\Vert{Y}^{M_{{\mathcal{P}}}}_{t}\Vert^{-2})((
\sum_{n=0}^{N(t)}\sum_{k=1}^{M(\gamma_{n+1})}\sum_{i=1}^{J^k_t}\mathbbm{1}_{\{\Gamma_n<T^k_i\leq%
\Gamma_{n+1}\wedge t\}} \vert\xi_i^k\vert^2\underline{c}({Z}^k_i))+ \vert a^{{%
\mathcal{P}} }_{t}\vert^2).
\end{eqnarray*}

We denote 
\begin{eqnarray}
\chi _{t}^{M_{{\mathcal{P}}}}=\sum_{n=0}^{N(t)}\sum_{k=1}^{M(\gamma_{n+1})}
\sum_{i=1}^{J^k_t}\mathbbm{1}_{\{\Gamma_n<T^k_i\leq%
\Gamma_{n+1}\wedge t\}}\vert\xi_i^k\vert^2\underline{%
c}({Z} ^k_i).  \label{chi}
\end{eqnarray}
By (\ref{cM2}), $(\mathbb{E}\sup\limits_{0\leq t\leq 2}\Vert\widetilde{{Y}}
^{M_{{\mathcal{P}}}}_{t}\Vert^{4dp}\Vert{Y}^{M_{{\mathcal{P}}
}}_{t}\Vert^{4dp})^{1/2}\leq C_{d,p}<\infty,$ so that using Schwartz
inequality, we have 
\begin{eqnarray}
\mathbb{E}\vert\frac{1}{\det \sigma_{{X}^{M_{{\mathcal{P}}}}_{t}}}
\vert^p\leq\mathbb{E}(\vert \lambda_{t}^{M_{{\mathcal{P}} }}\vert^{-dp})\leq
C(\mathbb{E}(\vert \chi _{t}^{M_{{\mathcal{P}} }}+\vert a^{{\mathcal{P}}%
}_{t}\vert^2 \vert^{-2dp}))^{\frac{1}{2}}.  \label{Step1}
\end{eqnarray}

\textbf{Step 2 }Since it is not easy to compute $\mathbb{E}(\vert \chi
_{t}^{M_{{\mathcal{P}}}}+\vert a^{{{\mathcal{P}}}}_{t}\vert^2 \vert^{-2dp}))$
directly, we make the following argument where the idea comes originally
from $\cite{ref5}$. Let $\Gamma(p)=\int_0^\infty {s} ^{p-1}e^{-{s}}d{s}$ be
the Gamma function. By a change of variables, we have the numerical equality 
\[
\frac{1}{\vert \chi _{t}^{M_{{\mathcal{P}}}}+\vert a^{{{\mathcal{P}}}
}_{t}\vert^2\vert^{2dp}}=\frac{1}{\Gamma (2dp)}\int_{0}^{\infty }{s}
^{2dp-1}e^{-{s}(\chi _{t}^{M_{{\mathcal{P}}}}+\vert a^{{{\mathcal{P}} }%
}_{t}\vert^2)}d{s},
\]
which, by taking expectation, gives 
\begin{eqnarray}
\mathbb{E}(\frac{1}{\vert \chi _{t}^{M_{{\mathcal{P}}}}+\vert a^{{\ {\ 
\mathcal{P}}}}_{t}\vert^2\vert^{2dp}})=\frac{1}{\Gamma (2dp)}
\int_{0}^{\infty }{s}^{2dp-1}\mathbb{E}(e^{-{s}(\chi _{t}^{M_{ {\mathcal{P }}%
}}+\vert a^{{{\mathcal{P}}}}_{t}\vert^2)})d{s}.  \label{Gamma}
\end{eqnarray}

\textbf{Step 3} Now we compute $\mathbb{E}(e^{-{s}(\chi _{t}^{M_{ {\ 
\mathcal{P}}}}+\vert a^{{{\mathcal{P}}}}_{t}\vert^2)})$ for any ${s}>0$. We
recall that $I_1=B_1$, $I_k=B_{k}-B_{k-1}, k\geq2$ (given in Section 4.3),
and $\xi_i^k=\Psi_k(Z^k_i)$ (see (\ref{xi})). Then 
\begin{eqnarray*}
\chi _{t}^{M_{{\mathcal{P}}}}
&=&\sum_{n=0}^{N(t)}\sum_{k=1}^{M(\gamma_{n+1})}\int_{\Gamma_n}^{%
\Gamma_{n+1}\wedge t}\int_{I_k} \vert\Psi_k(z)\vert^2 {\underline{c}}%
(z)N(dz,dr)=\int_{0}^{t} \int_{B_{M_{{\mathcal{P}}}(r)}} \Psi(z) {\underline{%
c}}(z)N(dz,dr),
\end{eqnarray*}
with 
\[
\Psi(z)=\sum\limits_{k=1}^{\infty}\vert\Psi_k(z)\vert^2\mathbbm{1}
_{I_k}(z)\geq \sum\limits_{k=1}^{\infty}\mathbbm{1}_{[k-\frac{3}{4 },k-\frac{%
1}{4}]}(\vert z\vert)\mathbbm{1} _{I_k}(z).
\]
Using It\^{o} formula, 
\begin{eqnarray*}
\mathbb{E}(e^{-{s}\chi _{t}^{M_{{\mathcal{P}}}}}) &=&1+\mathbb{E}
\int_0^t\int_{B_{M_{{\mathcal{P}}}(r)}} (e^{-{s}({\chi _{r-}^{M_{ {\ 
\mathcal{P}}}}}+\Psi(z)\underline{c}(z))}-e^{-{s}{\chi _{r-}^{M_{ { 
\mathcal{P}}}}}})\widehat{N}(dz,dr) \\
&=&1-\sum_{n=0}^{N(t)}(\int_{\Gamma_n}^{\Gamma_{n+1}\wedge t}\mathbb{E}(e^{-{%
s}{\chi _{r}^{M_{{\mathcal{P}}}}}})dr\sum_{k=1}^{M(\gamma_{n+1})}\int_{I_k}
(1-e^{-{s} \vert\Psi_k(z)\vert^2\underline{c}(z)})\mu(dz)).
\end{eqnarray*}
Solving the above equation we obtain 
\begin{eqnarray*}
\mathbb{E}(e^{-{s}{\chi _{t}^{M_{{\mathcal{P}}}}}})&=&\exp(-
\sum_{n=0}^{N(t)}(((\Gamma_{n+1}\wedge
t)-\Gamma_n)\sum_{k=1}^{M(\gamma_{n+1})}\int_{I_k} (1-e^{-{s}%
\vert\Psi_k(z)\vert^2\underline{c}(z)})\mu(dz))) \\
&\leq&\exp (-\sum_{n=0}^{N(t)}(((\Gamma_{n+1}\wedge
t)-\Gamma_n)\sum_{k=1}^{M(\gamma_{n+1})}\int_{I_k} (1-e^{-{s}\mathbbm{1}_{[k-%
\frac{3}{4 },k-\frac{1}{4}]}(\vert z\vert)\underline{c}(z)})\mu(dz))) \\
&=&\exp (-\sum_{n=0}^{N(t)}(((\Gamma_{n+1}\wedge
t)-\Gamma_n)\sum_{k=1}^{M(\gamma_{n+1})}\int_{I_k} (1-e^{-{s}\underline{c}%
(z)}) \mathbbm{1}_{[k-\frac{3}{4},k-\frac{1}{4}]}(\vert z\vert)\mu(dz))) \\
&=&\exp (-\sum_{n=0}^{N(t)}(((\Gamma_{n+1}\wedge
t)-\Gamma_n)\int_{B_{M(\gamma_{n+1})}} (1-e^{-{s}\underline{c}(z)})\overline{%
\mu}(dz)),
\end{eqnarray*}
with 
\[
\overline{\mu}(dz)=\sum_{k=1}^{\infty}\mathbbm{1}_{[k-\frac{3}{4},k-\frac{1}{
4}]}(\vert z\vert)\mu(dz).
\]

On the other hand, we denote 
\begin{eqnarray*}
\bar{\chi} _{t}^{M_{{\mathcal{P}}}}=\int_{0}^{t}\int_{B_{M_{ {\mathcal{P}}%
}(r)}^c} \Psi(z) {\underline{c}}(z)N(dz,dr),
\end{eqnarray*}
where $B_{m}^c$ denote the complementary set of $B_{m}$. Then in the same
way, 
\begin{eqnarray*}
\mathbb{E}(e^{-{s}\bar{\chi} _{t}^{M_{{\mathcal{P}}}}}) \leq \exp
(-\sum_{n=0}^{N(t)}(((\Gamma_{n+1}\wedge t)-\Gamma_n)\int_{B^c_{M({%
\gamma_{n+1}})}} (1-e^{-{s} \underline{c}(z)})\overline{\mu}(dz))).
\end{eqnarray*}
We recall by (\ref{abPt}) that for $\Gamma_{n}< t\leq \Gamma_{n+1}$, 
\[
a^{{\mathcal{P}}}_t=(\sum\limits_{\substack{ 1\leq i\leq n }}
\gamma_i\int_{\{\vert z\vert\geq M(\gamma_i)\}}\underline{c}%
(z)\mu(dz)+(t-\Gamma_n)\int_{\{ \vert z\vert\geq M(\gamma_{n+1})\}}%
\underline{c}(z)\mu(dz))^\frac{1}{2}.
\]
Then 
\[
a^{{\mathcal{P}}}_t\geq \sqrt{\mathbb{E}\bar{\chi} _{t}^{M_{ { \mathcal{P}}%
}}}.
\]
Using Jensen inequality for the convex function $f(x)=e^{-{s}x}$, ${s},x>0$,
we have 
\[
e^{-{s}\vert a^{M_{{\mathcal{P}}}}_{t}\vert^2}\leq e^{-{s}\mathbb{E}\bar{
\chi} _{t}^{M_{{\mathcal{P}}}}}\leq \mathbb{E}(e^{-{s}\bar{\chi} _{t}^{M_{{%
\mathcal{P}}}}})\leq \exp (-\sum_{n=0}^{N(t)}(((\Gamma_{n+1}\wedge
t)-\Gamma_n)\int_{B^c_{M(\gamma_{n+1})}} (1-e^{-{s}\underline{c}(z)})%
\overline{\mu} (dz))).
\]
So we deduce that 
\begin{eqnarray}
&&\mathbb{E}(e^{-{s}(\chi _{t}^{M_{{\mathcal{P}}}}+\vert a^{ {\ \mathcal{P}}%
}_{t}\vert^2)})=\mathbb{E}(e^{-{s}\chi _{t}^{M_{{\mathcal{P}} }}})\times e^{-%
{s}\vert a^{{{\mathcal{P}}}}_{t}\vert^2}  \nonumber \\
&&\leq \exp (-\sum_{n=0}^{N(t)}(((\Gamma_{n+1}\wedge
t)-\Gamma_n)\int_{B_{M(\gamma_{n+1})}} (1-e^{-{s}\underline{c}(z)})\overline{%
\mu}(dz)))  \nonumber \\
&&\times \exp (-\sum_{n=0}^{N(t)}(((\Gamma_{n+1}\wedge
t)-\Gamma_n)\int_{B^c_{M(\gamma_{n+1})}} (1-e^{-{s}\underline{c}(z)})%
\overline{\mu} (dz)))  \nonumber \\
&&= \exp (-t\int_{\mathbb{R}^d} (1-e^{-{s}\underline{c}(z)})\overline{\mu}
(dz)),  \label{Step3}
\end{eqnarray}
and the last term does not depend on ${M_{{\mathcal{P}}}}(t).$

\bigskip

\textbf{Step 4} Now we use the Lemma 14 from $\cite{ref3}$, which states the
following.

\begin{lemma}
We consider an abstract measurable space $B$, a $\sigma$-finite measure $%
\mathcal{M}$ on this space and a non-negative measurable function $%
f:B\rightarrow\mathbb{R}_+$ such that $\int_Bfd\mathcal{M}<\infty.$ For $t>0$
and $p\geq1$, we note 
\[
\beta_f({s})=\int_B(1-e^{-{s}f(x)})\mathcal{M}(dx)\quad and\quad
I_{t}^p(f)=\int_0^\infty {s}^{p-1}e^{-t\beta_f({s})}d{s}.
\]
We suppose that for some $t>0$ and $p\geq1$, 
\begin{equation}
\underline{\lim }_{u\rightarrow \infty }\frac{1}{\ln u}\mathcal{M}(f\geq 
\frac{1}{u})>\frac{p}{t},  \label{cm}
\end{equation}
then $I_{t}^p(f)<\infty.$
\end{lemma}

\bigskip

We will use the above lemma for $\mathcal{M}(dz)=\overline{\mu}(dz)$, $f(z)= 
\underline{c}(z)$ and $B=\mathbb{R}^d$. Thanks to (\ref{undergamma}) in 
\textbf{Hypothesis 2.4}, 
\begin{equation}
\underline{\lim }_{u\rightarrow \infty }\frac{1}{\ln u}\overline{\mu}( 
\underline{c}\geq \frac{1}{u})=\infty.  \label{cm6}
\end{equation}
Then for every $p\geq 1, 1\leq t\leq 3$, we deduce from (\ref{Step1}),(\ref%
{Gamma}),(\ref{Step3}) and \textbf{Lemma 5.4} that 
\begin{eqnarray}
\mathbb{E}\vert\frac{1}{\det \sigma_{{X}^{M_{{\mathcal{P}}}}_{t}}}
\vert^p&\leq&\mathbb{E}(\vert \lambda_{t}^{M_{{\mathcal{P}}
}}\vert^{-dp})\leq C(\mathbb{E}(\vert \chi _{t}^{M_{{\mathcal{P}} }}+\vert
a^{{{\mathcal{P}}}}_{t}\vert^2 \vert^{-2dp}))^{\frac{1}{2}}  \nonumber \\
&\leq& C(\frac{1}{\Gamma (2dp)}\int_{0}^{\infty }{s}^{2dp-1}\mathbb{E}(e^{-{%
s }(\chi _{t}^{M_{{\mathcal{P}}}}+\vert a^{{{\mathcal{P}}} }_{t}\vert^2)})d{s%
})^\frac{1}{2}  \nonumber \\
&\leq &C(\frac{1}{\Gamma (2dp)}\int_{0}^{\infty }{s}^{2dp-1}\exp (-t\int_{ 
\mathbb{R}^d }(1-e^{-{s}\underline{c}(z)})\overline{\mu}(dz)d{s})^\frac{1}{2}
<\infty.\quad \quad\quad  \label{detXt}
\end{eqnarray}

\bigskip

\textbf{Proof of $ii)$} We consider an increasing sequence of partition ${%
\mathcal{P}}_n,\ n\in\mathbb{N}$, (${\mathcal{P}}_n\subset {\ \mathcal{P}}%
_{n+1}$), such that $\vert{\mathcal{P}}_n\vert\downarrow0$. In particular, $%
\forall t,\ M_{ {\mathcal{P}_n}}(t)\uparrow\infty$.

We recall the notations ${F}^{M_{{\mathcal{P}}}}_{t}(x)={X}^{M_{ {\mathcal{P}%
}}}_{t}(x)-x$ and ${F}_{t}(x)={X}_{t}(x)-x$. We notice by \textbf{Lemma 5.1 $%
ii)$} that $\mathbb{E}\vert {F}^{M_{ {\mathcal{P}}_n}}_{t}- {F}%
_t\vert\rightarrow0$ as $n\rightarrow\infty$ , and by \textbf{Lemma 5.2}
that $\sup\limits_n\sup\limits_x\Vert {F} ^{M_{{\mathcal{P}}%
_n}}_{t}(x)\Vert_{L,l,p}\leq C_{l,p}$.

Moreover, by \textbf{Lemma 5.5} $ii)$ (given immediately below), we know
that $(D {F}^{M_{ { \mathcal{P}}_n}}_{t})_{n\in\mathbb{N}}$ is a Cauchy
sequence in $L^2(\Omega;l^2\times\mathbb{R}^d)$. Then applying \textbf{Lemma
3.3 (B)} with $F_{n}= {F}^{M_{{\mathcal{P}}_n}}_{t}$ and $F= {F}_t$, \textbf{%
Lemma 5.3 $i)$} implies \textbf{Lemma 5.3 $ii)$}.

\qed

\bigskip

\subsection{Auxiliary results}

Besides the lemmas concerning the Sobolev norms and covariance matrices, we
establish an auxiliary result. We recall $\varepsilon_m$ given in (\ref{epsM}%
).

\begin{lemma}
We assume that \textbf{Hypothesis 2.1} and \textbf{Hypothesis 2.4} $b)$ hold
true.

$i)$ Then for any $\varepsilon_0>0$, there exists a constant $C$ dependent
on $d,\varepsilon_0$ such that for every $t\in[0,3]$ and every stating point 
$x\in\mathbb{R}^d$, we have 
\[
\mathbb{E}\vert \det\sigma_{{X}_{t}^{{\mathcal{P}},M_{ {\mathcal{P}}%
}}}-\det\sigma_{{X}_{t}}\vert^{\frac{ 2}{1+\varepsilon_0}}\leq C\Vert D{X}%
_{t}^{{{\mathcal{P}}},M_{ { \mathcal{P}}}}-D{X}_{t}\Vert_{L^2(\Omega;l_2
\times\mathbb{R}^d)}^{\frac{2}{1+\varepsilon_0}}\leq C(\vert{\mathcal{P}}
\vert+\varepsilon_{M(\vert\mathcal{P}\vert)})^{\frac{2}{(2+\varepsilon_0)(1+%
\varepsilon_0)}}.
\]

$ii)$ We consider an increasing sequence of partition ${\mathcal{P}}_n,\ n\in%
\mathbb{N}$, (${\mathcal{P}}_n\subset{\mathcal{P}}_{n+1}$), such that $\vert{%
\mathcal{P}}_n\vert\downarrow0$. In particular, $\forall t,\ M_{ {\mathcal{P}%
_n}}(t)\uparrow\infty$. We denote 
\[
F_n(x)={X}_{t}^{M_{{\mathcal{P}}_n}}(x).
\]
Then for each $x\in\mathbb{R}^d$, the sequence $DF_n(x),\ n\in\mathbb{N}$ is
Cauchy in $L^2(\Omega;l_2\times\mathbb{R}^d)$, uniformly with respect to $x$
: 
\[
\sup_x\Vert DF_n(x)-DF_m(x)\Vert_{L^2(\Omega;l_2\times\mathbb{R}
^d)}\rightarrow0,\ \text{as }n,m\rightarrow\infty.
\]
\end{lemma}

\begin{proof}
\textbf{Proof of $i)$}

By \textbf{Lemma 5.2}, we know that $\Vert D{X}^{{ { \mathcal{P}}},M_{{%
\mathcal{P}}}}_t\Vert_{L^2(\Omega;l_2\times\mathbb{R}^d)}$ and $\Vert D{X}%
_t\Vert_{L^2(\Omega;l_2\times\mathbb{R}^d)}$ are bounded, uniformly with
respect to $x$. Then using H$\ddot{o}$lder's inequality with conjugates $%
1+\varepsilon_0$ and $\frac{1+\varepsilon_0}{\varepsilon_0}$, we get 
\begin{eqnarray}
&&\mathbb{E}(\vert\det\sigma_{{X}^{{{\mathcal{P}}},M_{ {\mathcal{P} }%
}}_t}-\det\sigma_{{X}_t}\vert^\frac{2}{ 1+\varepsilon_0})\leq C\Vert D{X}%
_{t}^{{{\mathcal{P}}},M_{ {\ \mathcal{P}}}}-D{X}_{t}\Vert_{L^2(\Omega;l_2
\times\mathbb{R}^d)}^{\frac{2}{1+\varepsilon_0}}.
\end{eqnarray}
Now we only need to prove that 
\begin{eqnarray}
\Vert D{X}_{t}^{{{\mathcal{P}}},M_{{\mathcal{P}}}}-D{X} _{t}\Vert_{L^2(%
\Omega;l_2\times\mathbb{R}^d)}^{\frac{ 2}{1+\varepsilon_0}}\leq C(\vert{%
\mathcal{P}} \vert+\varepsilon_{M(\vert\mathcal{P}\vert)})^{\frac{2}{
(2+\varepsilon_0)(1+\varepsilon_0)}}.  \label{DDD}
\end{eqnarray}
The proof of (\ref{DDD}) is a slight variant of the proof of Lemma 3.9 $iii)$
in the paper $\cite{ref20}$. The difference in that the truncation function $%
M$ is constant in $\cite{ref20}$ while here it may vary on different time
intervals. We do not discuss in detail here. So we conclude that \textbf{\
Lemma 5.5 $i)$} holds.

\bigskip

\textbf{Proof of $ii)$}

We consider an increasing sequence of partition ${\mathcal{P}}_n,\ n\in 
\mathbb{N}$, (${\mathcal{P}}_n\subset{\mathcal{P}}_{n+1}$), such that $\vert{%
\mathcal{P}}_n\vert\downarrow0$. In particular, $\forall t,\ M_{ {\mathcal{P}%
_n}}(t)\uparrow\infty$. We need to prove that 
\begin{eqnarray}
\Vert D{X}^{M_{{\mathcal{P}}_n}}_{t}-D{X}_{t}^{M_{ {\mathcal{P }}%
_m}}\Vert_{L^2(\Omega;l_2\times\mathbb{R}^d)}\rightarrow0,\ \text{as }
n,m\rightarrow\infty.  \label{Cauchy}
\end{eqnarray}

The proof of (\ref{Cauchy}) is a slight variant of the proof of (148)
p.47-49 in $\cite{ref20}$, so we omit it.

\bigskip

\end{proof}

\bigskip\bigskip

\section{Proof of Theorem 4.1}

In this section, we give the proof of \textbf{Theorem 4.1}. We apply \textbf{%
Proposition 2.1.1} in Section 2. For a measurable function $f$, we denote $%
\overline{P}_tf(x)=\mathbb{E}f(X^{\mathcal{P},{M_\mathcal{P}}}_{t}(x))$ and $%
P_tf(x)=\mathbb{E}f(X_{t}(x))$. In the following subsections, we will check
the conditions of \textbf{Proposition 2.1.1}.

\subsection{Euler: condition (\protect\ref{app3})}

For every $\gamma>0$, we recall in (\ref{Mgamma}) that we define $%
M(\gamma)\in\mathbb{N}$ such that 
\begin{eqnarray*}
\varepsilon_{M(\gamma)}\leq \gamma^2.
\end{eqnarray*}

We recall the basic equation $X_t$ (see (\ref{1.1})). We denote by $\Tilde{X}%
^{M_\mathcal{P}}_t$ the one step truncated Euler scheme: 
\[
\Tilde{X}^{M_\mathcal{P}}_{t}(x)=x+\int_{0}^{t}\int_{B_{M({\gamma}%
)}}c(z,x)dN(z,s)+\int_{0}^{t}b(x)ds.
\]

Then, 
\begin{eqnarray*}
\mathbb{E}\vert \Tilde{X}^{M_\mathcal{P}}_{\gamma}-X_{\gamma}\vert&\leq& 
\mathbb{E}\int_0^\gamma\int_{\{\vert z\vert\geq M(\gamma)\}}\vert
c(z,X_s)\vert dN(z,s)+\mathbb{E}\int_0^{\gamma}\int_{B_{M({\gamma})}}\vert
c(z,x)-c(z,X_s)\vert dN(z,s) \\
&+&\mathbb{E}\int_0^{\gamma}\vert b(x) -b(X_s)\vert ds \\
&\leq& \gamma\int_{\{\vert z\vert\geq M(\gamma)\}}\bar{c}(z)\mu(dz)+C%
\int_0^{\gamma}\mathbb{E}\vert x-X_s\vert ds \\
&\leq& \gamma \sqrt{\varepsilon_{M(\gamma)}}+C\times\gamma^2\leq
C\times\gamma^2
\end{eqnarray*}

So 
\[
W_1(\Tilde{X}^{M_\mathcal{P}}_{\gamma},X_{\gamma})\leq \mathbb{E}\vert 
\Tilde{X}^{M}_{\gamma}-X_{\gamma}\vert\leq C\times{\gamma} ^{2 }.
\]
So we conclude that (\ref{app3}) holds for $\alpha=1$ and $k_0=0$.

\bigskip

\subsection{Lipschitz: condition (\protect\ref{app1}) and the existence of
an invariant measure}

We recall that $X$ is the solution of the equation (\ref{1.1}). 

\begin{lemma}
Suppose that \textbf{Hypothesis 2.5} (see (\ref{2.5ab}) and (\ref{2.5c}))
holds.

a)\qquad Then, for a Lipschitz continuous function $f$%
\[
\left\vert \mathbb{E}(f(X_{t}(x))-\mathbb{E}(f(X_{t}(y))\right\vert \leq
L_{f}e^{-\frac{\theta }{2}t}\left\vert x-y\right\vert ,
\]%
with $L_{f}$ the Lipschitz constant of $f$.

b)\qquad Moreover, there exists at least one invariant probability.
\end{lemma}

\textbf{Proof a)} We fix $x,y\in \mathbb{R}^{d}$ and we construct on the
same probability space, with the same Poisson point measure $N$ the solution 
$X_{t}^{M}(y)$ which starts from $y.$ Then we denote%
\begin{eqnarray*}
Y_{t} &=&X_{t}(x)-X_{t}(y),\quad \\
\Delta _{t}^{c}(z) &=&c(z,X_{s-}(x))-c(z,X_{s-}(y)) \\
\Delta _{t}^{b} &=&b(X_{s-}(x))-b(X_{s-}(y))
\end{eqnarray*}%
and we have%
\[
Y_{t}=x-y+\int_{0}^{t}\int_{\mathbb{R}^{d}}\Delta
_{s}^{c}(z)dN(z,s)+\int_{0}^{t}\Delta _{s}^{b}ds.
\]%
Using It\^{o}'s formula for $\Phi (t,u)=e^{\lambda t}\left\vert u\right\vert
^{2}$ we obtain 
\begin{eqnarray*}
\Phi (t,Y_{t}) &=&\left\vert x-y\right\vert ^{2}+\lambda \int_{0}^{t}\Phi
(s,Y_{s})ds+\int_{0}^{t}2e^{\lambda s}\left\langle Y_{s},\Delta
_{s}^{b}\right\rangle ds \\
&&+\int_{0}^{t}\int_{\mathbb{R}^{d}}(\Phi (s,Y_{s-}+\Delta _{s}^{c}(z))-\Phi
(s,Y_{s-}))dN(z,s) \\
&=&\left\vert x-y\right\vert ^{2}+\lambda \int_{0}^{t}\Phi
(s,Y_{s})ds+\int_{0}^{t}2e^{\lambda s}\left\langle Y_{s},\Delta
_{s}^{b}\right\rangle ds \\
&&+M_{t}+\int_{0}^{t}\int_{\mathbb{R}^{d}}(\Phi (s,Y_{s-}+\Delta
_{s}^{c}(z))-\Phi (s,Y_{s-}))d\mu (z)ds
\end{eqnarray*}%
with $M_{t}$ a martingale. Taking the expectation we get 
\[
e^{\lambda t}\mathbb{E}\left\vert Y_{t}\right\vert ^{2}\leq \left\vert
x-y\right\vert ^{2}+\int_{0}^{t}e^{\lambda s}\mathbb{E}(\Psi _{s})ds
\]%
with 
\begin{eqnarray*}
\Psi _{s} &=&\lambda \left\vert Y_{s}\right\vert ^{2}+2\left\langle
Y_{s},\Delta _{s}^{b}\right\rangle +\int_{\mathbb{R}^{d}}\left\vert
Y_{s}+\Delta _{s}^{c}(z)\right\vert ^{2}-\left\vert Y_{s}\right\vert ^{2}\mu
(dz) \\
&=&\lambda \left\vert Y_{s}\right\vert ^{2}+2\left\langle Y_{s},\Delta
_{s}^{b}\right\rangle +\int_{\mathbb{R}^{d}}\left\langle \Delta
_{s}^{c}(z),2Y_{s}+\Delta _{s}^{c}(z)\right\rangle \mu (dz).
\end{eqnarray*}%
We need to prove that $\mathbb{E}(\Psi _{s})\leq 0. $ We recall that we
assume \textbf{Hypothesis 2.5} $i)ii)$ (see (\ref{2.5ab})). We also have%
\[
\left\vert \left\langle \Delta _{s}^{c}(z),2Y_{s}+\Delta
_{s}^{c}(z)\right\rangle \right\vert \leq (2\bar{c}(z)+\bar{c}%
^{2}(z))\left\vert Y_{s}\right\vert ^{2},
\]%
so that%
\[
\Psi _{s}\leq \left\vert Y_{s}\right\vert ^{2}(\lambda +\int_{\mathbb{R}%
^{d}}(2\bar{c}(z)+\bar{c}^{2}(z))\mu (dz)-2\overline{b}).
\]%
Thanks to \textbf{Hypothesis 2.5} $iii)$ (see (\ref{2.5c})), taking $\lambda
\leq \theta $, we have 
\[
e^{\lambda t}\mathbb{E}\left\vert Y_{t}\right\vert ^{2}\leq \left\vert
x-y\right\vert ^{2}+\int_{0}^{t}e^{\lambda s}\mathbb{E}(\Psi _{s})ds\leq
\left\vert x-y\right\vert ^{2}
\]%
so that 
\[
\mathbb{E}\left\vert X_{t}(x)-X_{t}(y)\right\vert ^{2}\leq e^{-\theta
t}\left\vert x-y\right\vert ^{2}.
\]%
Then, for a Lipschitz continuous function $f$,%
\[
\left\vert \mathbb{E}(f(X_{t}(x))-\mathbb{E}(f(X_{t}(y))\right\vert \leq
L_{f}\mathbb{E}\left\vert X_{t}(x)-X_{t}(y)\right\vert \leq L_{f}e^{-\frac{%
\theta }{2}t}\left\vert x-y\right\vert .
\]

b)\qquad We denote ${L}$ to be the infinitesimal operator of (\ref{1.1}). We
take $V(x)=\left\vert x\right\vert ^{2}$ and we will prove that 
\[
{L}V\leq \bar{\beta} -\bar{\alpha} V
\]
for some $\bar{\beta} ,\bar{\alpha} >0$ (the Lyapunov mean reverting
condition). This implies ${L}V\leq C$ and $\lim\limits_{\left\vert
x\right\vert \rightarrow \infty }{L}V(x)=-\infty .$ Then we use Theorem 9.3
and Lemma 9.7 from $\cite{ref23}$ (with $\varphi =V$ and $\psi ={L}V$ )\
which guarantees existence of an invariant distribution. We have%
\begin{eqnarray*}
{L}V(x) &=&2\left\langle x,b(x)\right\rangle +\int_{\mathbb{R}%
^{d}}(V(x+c(z,x))-V(x))\mu (dz) \\
&\leq &-2\overline{b}\left\vert x\right\vert ^{2}+\int_{\mathbb{R}%
^{d}}(2\left\langle x,c(z,x)\right\rangle +\left\vert c(z,x)\right\vert
^{2})\mu (dz) \\
&\leq &-2\overline{b}\left\vert x\right\vert ^{2}+(\left\vert x\right\vert
^{2}+1)\int_{\mathbb{R}^{d}}\bar{c}(z)\mu (dz)+\int_{\mathbb{R}^{d}}\bar{c}%
^{2}(z)\mu (dz) \\
&=&\int_{\mathbb{R}^{d}}(\bar{c}(z)+\bar{c}^{2}(z))\mu (dz)-(2\overline{b}%
-\int_{\mathbb{R}^{d}}\bar{c}(z)\mu (dz))\left\vert x\right\vert ^{2}.
\end{eqnarray*}%
$\square $

\bigskip

\subsection{Regularity: conditions (\protect\ref{app9}), (\protect\ref{app9'}%
), (\protect\ref{app10}) and (\protect\ref{app11})}

Firstly, we deal with (\ref{app9}). Let $t\in[1,2]$. For any $k$ and any
multi-index $\beta_0$ with $\vert\beta_0\vert=k$, we write 
\[
\partial^{\beta_0}_xP_t\varphi(x)=\mathbb{E}[\partial^{\beta_0}_x(%
\varphi(X_t(x))]=\sum_{\vert\alpha_0\vert\leq\vert\beta_0\vert}\mathbb{E}%
[(\partial^{\alpha_0}\varphi)(X_t(x))\mathbf{P}_{\alpha_0}(x)],
\]%
with $\mathbf{P}_{\alpha_0}(x)$ a polynomial of $\partial^{\alpha_1}_xX_t(x)$%
, $\vert\alpha_1\vert\leq\vert\beta_0\vert$.

In the following, we use the results from Section 5. In \textbf{Lemma 5.2},
we prove that the Sobolev norms of each $\partial^{\alpha_1}_xX_t(x)$ are
bounded, uniformly with respect to $x$. It follows that this is also true
for $\mathbf{P}_{\alpha_0}(x)$.

We denote that $F_t(x)=X_t(x)-x$. In \textbf{Lemma 5.2}, we have proved that
the Sobolev norms of each $F_t(x)$ are bounded, uniformly with respect to $x$%
. Moreover, in \textbf{Lemma 5.3}, we prove that $F_t(x)$ is
non-degenerated, uniformly with respect to $x$, that is $\Sigma_p(F_t(x))<%
\infty$ for each $p$ (see (\ref{Mcov})).

Then we use \textbf{Lemma 3.4 (A)} which asserts that $(\bm{B}_k)$ is true
for $F=F_t(x)$ and $G=\mathbf{P}_{\alpha_0}(x)$. By the remark of \textbf{%
Lemma 3.4}, $(\bm{B}_k)$ is also true for $F=X_t(x)=F_t(x)+x$ and $G=\mathbf{%
P}_{\alpha_0}(x)$. This reads 
\[
\vert\mathbb{E}[(\partial^{\alpha_0}\varphi)(X_t(x))\mathbf{P}%
_{\alpha_0}(x)]\vert\leq C\Vert\varphi\Vert_{\infty},
\]
which gives (\ref{app9}).

In a similar way, we can obtain (\ref{app9'}).

For (\ref{app10}), $i)$ is a direct consequence of (\ref{app1}) which has been proved in Section 6.2. For (\ref{app10}) $ii)$, we take $t\in(0,1]$. For any $k$ and any
multi-index $\beta_0$ with $\vert\beta_0\vert=k$, we notice that \[
\vert\partial^{\beta_0}_x\nabla P_t\varphi(x)\vert=\vert\mathbb{E}[\partial^{\beta_0}_x(%
\nabla\varphi(X_t(x))]\vert=\vert\sum_{\vert\alpha_0\vert\leq\vert\beta_0\vert}\mathbb{E}%
[(\partial^{\alpha_0}\nabla\varphi)(X_t(x))\mathbf{P}_{\alpha_0}(x)]\vert\leq\Vert\nabla\varphi\Vert_{k,\infty}\sum_{\vert\alpha_0\vert\leq\vert\beta_0\vert}\mathbb{E}%
\vert\mathbf{P}_{\alpha_0}(x)\vert,
\]with $\mathbf{P}_{\alpha_0}(x)$ a polynomial of $\partial^{\alpha_1}_xX_t(x)$%
, $\vert\alpha_1\vert\leq\vert\beta_0\vert$. In $\cite{ref13}$, Kunita has shown in Theorem 3.4.1 and Theorem 3.4.2 the regularity of the flow associated with the jump-diffusion. So in our case, we have $\mathbb{E}%
\vert\mathbf{P}_{\alpha_0}(x)\vert<\infty$ and thus (\ref{app10}) $ii)$ holds true.

\bigskip

Now we prove (\ref{app11}). In order to prove (\ref{app11}), we need to
represent $\overline{P}_{s,t}\varphi(x)$ and $P_{s,t}\varphi(x)$. So we
consider the following equations.

We denote $X^{\mathcal{P},{M_\mathcal{P}}}_{s,t}$ and $X_{s,t}$ the
solutions of the following equations respectively: 
\begin{eqnarray}
X^{\mathcal{P},{M_\mathcal{P}}}_{s,t}&=&x+\int_{s}^{t}b(X^{\mathcal{P},{M_%
\mathcal{P}}}_{s,\tau(r)})dr +\int_{s}^{t}\int_{B_{M_{\mathcal{P}}(r)}}{c}%
(z,X^{\mathcal{P},{M_\mathcal{P}}}_{s,\tau(r)-})N(dz,dr);
\label{truncationEst}
\end{eqnarray}
\begin{eqnarray}
X_{s,t}&=&x+\int_{s}^{t}b(X_{s,r})dr +\int_{s}^{t}\int_{\mathbb{R}^d}{c}%
(z,X_{s,r-})N(dz,dr).  \label{Xst}
\end{eqnarray}
We sometimes write ${X}^{{\mathcal{P}},{M_\mathcal{P}}}_{s,t}(x)$ (and ${X}%
_{s,t}(x)$) instead of ${X}^{{\mathcal{P}},{M_\mathcal{P}}}_{s,t}$ (and ${X}%
_{s,t}$) to stress the dependence on the initial value $x$. And we denote $%
\overline{P}_{s,t}\varphi(x)=\mathbb{E}\varphi({X}^{{\mathcal{P}},{M_%
\mathcal{P}}}_{s,t}(x))$ and $P_{s,t}\varphi(x)=\mathbb{E}\varphi(X_{s,t}(x))
$.

Let $1<t<r<t+2$. We recall that $\mathcal{P}=\{0=\Gamma_0<\Gamma_1<\cdots<%
\Gamma_{l-1}<\Gamma_l<\cdots\}$, $\gamma_l=\Gamma_l-\Gamma_{l-1}$ and for $%
\Gamma _{l}\leq t<\Gamma _{l+1}$, $N(t)=l.$ We denote 
\begin{eqnarray}
F_{r-t+1}^{\mathcal{P},{M_\mathcal{P}}}(x)=X_{t,r}\circ X_{t-1,t}^{\mathcal{P%
},{M_\mathcal{P}}}(x)-x\ \text{and}\ F_{r-t+1}(x)=X_{t,r}\circ
X_{t-1,t}(x)-x=X_{t-1,r}(x)-x.  \label{Fregul}
\end{eqnarray}
We also denote $\vert\mathcal{P}^{t-1,t}\vert:=\max\limits_{\substack{ l\in%
\mathbb{N}\ s.t.  \\ \Gamma_{l+1}>t-1,\Gamma_{l}< t }}((\Gamma_{l+1}\wedge
t)-(\Gamma_l\vee (t-1)))$. Before we give the proof of (\ref{app11}), we
state the following lemma concerning $F_{r-t+1}^{\mathcal{P},{M_\mathcal{P}}%
}(x)$ and $F_{r-t+1}(x)$ given in (\ref{Fregul}).

\begin{lemma}
Under the \textbf{Hypothesis 2.1$\sim$2.4}, we have these results.

$i)$ For all $p\geq1, q\geq 0$, there exists a constant $C_{q,p}$ depending
on $q,p,d$, such that $F_{r-t+1}^{\mathcal{P},{M_\mathcal{P}}}(x)$ and $%
F_{r-t+1}(x)$ belong to $\mathcal{D}_{\infty}^d$ and 
\[
\sup\limits_{x}\sup\limits_{\mathcal{P}}\Vert F_{r-t+1}^{\mathcal{P},{M_%
\mathcal{P}}}(x)+F_{r-t+1}(x)\Vert_{L,q,p}\leq C_{q,p}.
\]

$ii)$ For every $p\geq1$, we have 
\begin{eqnarray*}
\sup\limits_{{\mathcal{P}}}\sup\limits_{x}\mathbb{E}(1/
\det\sigma_{F_{r-t+1}(x)})^p\leq C_p,
\end{eqnarray*}
with $C_p$ a constant depending on $p,d$.

$iii)$ For any $\varepsilon_0>0$, there exists a constant $C$ dependent on $%
d,\varepsilon_0$ such that 
\[
\sup_x\mathbb{E}\vert \det\sigma_{F_{r-t+1}^{\mathcal{P},{M_\mathcal{P}}%
}(x)}-\det\sigma_{F_{r-t+1}(x)}\vert^{\frac{2}{1+\varepsilon_0}}\leq C\vert{%
\mathcal{P}^{t-1,t}}\vert^{\frac{2}{(2+\varepsilon_0)(1+\varepsilon_0)}}.
\]
\end{lemma}

\begin{proof}
Firstly, we will construct an approximation scheme for $X_{t,r}\circ
X_{t-1,t}^{\mathcal{P},{M_\mathcal{P}}}(x)$. We take an integer $N_0$ such
that $\frac{1}{2^{N_0}}\leq \vert\mathcal{P}\vert$. For $n>N_0$, we take a
"mixed partition" 
\begin{eqnarray*}
&&\mathcal{P}_n=\{t-1<\Gamma_{N(t-1)+1}<\cdots<\Gamma_{N(t)}\leq t \\
&&<t+\frac{1}{2^n}(r-t)<t+\frac{2}{2^n}(r-t)<\cdots<t+\frac{l}{2^n}(r-t)<t+%
\frac{l+1}{2^n}(r-t)<\cdots<r\} \\
&&:=\{t-1=s_0<s_1<\cdots<s_{n_0}=r\}.
\end{eqnarray*}%
We remark that we take the partition $\{ \Gamma_l\}$ on $[t-1,t]$ and take the partition $\{\frac{l}{2^n}\}$ on $[t,r]$.
We denote $\vert\mathcal{P}_n\vert:=\max\limits_{k\in\{1,\cdots,n_0%
\}}s_k-s_{k-1}.$ We construct $M_{\mathcal{P}_n}(t)=M(s_{l+1}-s_{l})$ when $%
s_l<t\leq s_{l+1}$ with the truncation function $M(\bullet)$ given in (\ref%
{Mgamma}). And we denote $\tau^{{\mathcal{P}}_n} (t)=s _{l}$ when $s_l<t\leq
s_{l+1}$. Then we consider the truncated Euler scheme based on ${\mathcal{P}%
_n},\ M_{\mathcal{P}_n}$: 
\begin{eqnarray*}
X^{\mathcal{P}_n,{M_{\mathcal{P}_n}}}_{t-1,r}&=&x+\int_{t-1}^{r}b(X^{%
\mathcal{P}_n,{M_{\mathcal{P}_n}}}_{t-1,\tau^{{\mathcal{P}}_n} (s)})ds
+\int_{t-1}^{r}\int_{B_{M_{\mathcal{P}_n}(s)}}{c}(z,X^{\mathcal{P}_n,{M_{%
\mathcal{P}_n}}}_{t-1,\tau^{{\mathcal{P}}_n} (s)-})N(dz,ds).
\end{eqnarray*}

We denote 
\begin{eqnarray}
F^{\mathcal{P}_n,{M_{\mathcal{P}_n}}}_{r-t+1}(x)=X^{\mathcal{P}_n,{M_{%
\mathcal{P}_n}}}_{t-1,r}(x)-x.  \label{Fnregul}
\end{eqnarray}
We notice that we can apply the results from Section 5 for $F^{\mathcal{P}_n,%
{M_{\mathcal{P}_n}}}_{r-t+1}(x)$, $F^{\mathcal{P},{M_{\mathcal{P}}}%
}_{r-t+1}(x)$ and $F_{r-t+1}(x)$ defined in (\ref{Fregul}) and (\ref{Fnregul}%
).

Since $r-t+1<3$, by \textbf{Lemma 5.2 $i)$}, the Sobolev norms of $F^{%
\mathcal{P}_n,{M_{\mathcal{P}_n}}}_{r-t+1}(x)$ and $F_{r-t+1}(x)$ are
bounded, uniformly with respect to $x$. One can check that $F^{\mathcal{P}_n,%
{M_{\mathcal{P}_n}}}_{r-t+1}(x)\rightarrow F^{\mathcal{P},{M_{\mathcal{P}}}%
}_{r-t+1}(x)$  in $L^1(\Omega)$, as $n\rightarrow\infty$ (which is a variant
of \textbf{Lemma 5.1 $i)$}). So we can apply \textbf{Lemma 3.3 (A)} for $%
F_n=F^{\mathcal{P}_n,{M_{\mathcal{P}_n}}}_{r-t+1}(x)$ and $F=F^{\mathcal{P},{%
M_{\mathcal{P}}}}_{r-t+1}(x)$ in order to get that $F_{r-t+1}^{\mathcal{P},{%
M_\mathcal{P}}}(x)\in\mathcal{D}_{\infty}^d$ and $\sup\limits_{x}\sup%
\limits_{\mathcal{P}}\Vert F_{r-t+1}^{\mathcal{P},{M_\mathcal{P}}%
}(x)\Vert_{L,q,p}\leq C_{q,p}.$ Hence, \textbf{Lemma 6.2} $i)$ is proved.

Moreover, since $r-t+1>1$, using \textbf{Lemma 5.3} $ii)$ we have $%
\sup\limits_{{\mathcal{P}}}\sup\limits_{x}\mathbb{E}(1/
\det\sigma_{F_{r-t+1}^{{M_\mathcal{P}}}(x)})^p\leq C_p$. So \textbf{Lemma
6.2 $ii)$} is proved.

Finally, by \textbf{Lemma 5.5 $i)$} and recalling by (\ref{Mgamma}) that $%
\varepsilon_{M(\gamma)}\leq \gamma^2$, we have 
\[
\Vert DF^{\mathcal{P}_n,{M_{\mathcal{P}_n}}}_{r-t+1}(x)-DF_{r-t+1}(x)%
\Vert_{L^2(\Omega;l_2\times\mathbb{R}^d)}^{\frac{2}{ 1+\varepsilon_0}}\leq
C(\vert\mathcal{P}_n\vert+\varepsilon_{M(\vert\mathcal{P}_n\vert)})^{\frac{2%
}{(2+ \varepsilon_0)(1+\varepsilon_0)}}\leq C\vert\mathcal{P} ^{t-1,t}\vert^{%
\frac{2}{(2+\varepsilon_0)(1+\varepsilon_0)}},
\]%
where the last equality is true since $\frac{1}{2^n}\leq\vert\mathcal{P}%
^{t-1,t}\vert$ for every $n>N_0$. Then we can apply \textbf{Lemma 3.3 (C)}
for $(F_n,\Bar{F}_n)=(F^{\mathcal{P}_n,{M_{\mathcal{P}_n}}%
}_{r-t+1}(x),F_{r-t+1}(x))$ and $(F,\bar{F})=(F^{\mathcal{P},{M_{\mathcal{P}}%
}}_{r-t+1}(x),F_{r-t+1}(x))$. So $\sup\limits_x\mathbb{E}\vert
\det\sigma_{F_{r-t+1}^{\mathcal{P},{M_\mathcal{P}}}(x)}-\det%
\sigma_{F_{r-t+1}(x)}\vert^{\frac{2}{1+\varepsilon_0}}\leq C\Vert DF^{%
\mathcal{P},{M_{\mathcal{P}}}}_{r-t+1}(x)-DF_{r-t+1}(x)\Vert_{L^2(\Omega;l_2%
\times\mathbb{R}^d)}^{\frac{2}{1+\varepsilon_0}}\leq C\vert{\mathcal{P}%
^{t-1,t}}\vert^{\frac{2}{(2+\varepsilon_0)(1+\varepsilon_0)}}$ and \textbf{%
Lemma 6.2 $iii)$} is proved.
\end{proof}

\bigskip\bigskip

Then we can prove (\ref{app11}). By \textbf{Lemma 6.2 $i)$}, the Sobolev
norms of $F_{r-t+1}^{\mathcal{P},{M_{\mathcal{P}}}}(x)$ are bounded,
uniformly with respect to $x$. Using \textbf{Lemma 6.2 $ii)$}, the
covariance matrix of $F_{r-t+1}(x)$ is non-degenerated. Then we are able to
apply \textbf{Lemma 3.5} for $F=F_{r-t+1}^{\mathcal{P},{M_{\mathcal{P}}}}(x)$
and $Q=F_{r-t+1}(x)$ so (\ref{1*}) holds for $F=F_{r-t+1}^{\mathcal{P},{M_{%
\mathcal{P}}}}(x)$ and $Q=F_{r-t+1}(x)$. Thanks to the remark of \textbf{%
Lemma 3.5}, (\ref{1*}) also holds for $F=X_{t,r}\circ X_{t-1,t}^{\mathcal{P},%
{M_{\mathcal{P}}}}(x)=F_{r-t+1}^{\mathcal{P},{M_{\mathcal{P}}}}(x)+x,\
Q=X_{t,r}\circ X_{t-1,t}(x)=F_{r-t+1}(x)+x $ and get 
\begin{eqnarray}
&&\left\vert \mathbb{E}(f(X_{t,r}\circ X_{t-1,t}^{\mathcal{P},{M_{\mathcal{P}%
}}}(x)))-\mathbb{E}(f_{\delta }(X_{t,r}\circ X_{t-1,t}^{\mathcal{P},{M_{%
\mathcal{P}}}}(x)))\right\vert  \nonumber \\
&&\leq C\left\Vert f\right\Vert _{\infty }\times (\frac{\delta ^{q}}{%
\eta^{2q}}+\eta^{-p}{\mathbb{E}\vert\det\sigma_{X_{t,r}\circ X_{t-1,t}^{%
\mathcal{P},{M_{\mathcal{P}}}}(x)}-\det\sigma_{X_{t,r}\circ
X_{t-1,t}(x)}\vert}^p+\eta^\kappa),  \label{regul}
\end{eqnarray}%
where we have also used the fact that $\sup\limits_{{\mathcal{P}}%
}\sup\limits_{x}\mathbb{E}(1/ \det\sigma_{X_{t,r}\circ
X_{t-1,t}(x)})^\kappa\leq C_\kappa$ from \textbf{Lemma 6.2 $ii)$}.\newline
We take $p=\frac{2}{1+\varepsilon_0}$ for any small $\varepsilon_0$. Thanks
to \textbf{Lemma 6.2 $iii)$}, 
\[
\sup_x\mathbb{E}\vert\det\sigma_{X_{t,r}\circ X_{t-1,t}^{ \mathcal{P},{M_{%
\mathcal{P}}}}(x)}-\det\sigma_{X_{t,r}\circ X_{t-1,t}(x)}\vert ^{\frac{2}{%
1+\varepsilon_0}}=\sup_x\mathbb{E}\vert \det\sigma_{F_{r-t+1}^{\mathcal{P},{%
M_\mathcal{P}}}(x)}-\det\sigma_{F_{r-t+1}(x)}\vert^{\frac{2}{1+\varepsilon_0}%
}\leq C\vert{\mathcal{P}^{t-1,t}}\vert^{\frac{2}{(2+\varepsilon_0)(1+%
\varepsilon_0)}}.
\]%
This implies that
\[
\sup_x\mathbb{E}\vert\det\sigma_{X_{t,r}\circ X_{t-1,t}^{ \mathcal{P},{M_{%
\mathcal{P}}}}(x)}-\det\sigma_{X_{t,r}\circ X_{t-1,t}(x)}\vert ^{\frac{2}{%
1+\varepsilon_0}}\leq C\gamma_{N(t-1)}^{\frac{2}{(2+\varepsilon_0)(1+%
\varepsilon_0)}}.
\]
Substituting into (\ref{regul}), we obtain 
\[
\sup_x\left\vert \mathbb{E}(f(X_{t,r}\circ X_{t-1,t}^{\mathcal{P},{M_%
\mathcal{P}}}(x)))-\mathbb{E}(f_{\delta }(X_{t,r}\circ X_{t-1,t}^{\mathcal{P}%
,{M_\mathcal{P}}}(x)))\right\vert \leq C\left\Vert f\right\Vert _{\infty
}\times (\frac{\delta ^{q}}{\eta^{2q}}+\eta^{-\frac{2}{1+\varepsilon_0}}{%
\gamma_{N(t-1)}^{\frac{2}{(2+\varepsilon_0)(1+\varepsilon_0)}}}+\eta^\kappa).
\]

By a similar argument, we have 
\[
\sup_x\left\vert \mathbb{E}(f(X_{t,r}^{\mathcal{P},{M_\mathcal{P}}}\circ
X_{t-1,t}^{\mathcal{P},{M_\mathcal{P}}}(x)))-\mathbb{E}(f_{\delta }(X_{t,r}^{%
\mathcal{P},{M_\mathcal{P}}}\circ X_{t-1,t}^{\mathcal{P},{M_\mathcal{P}}%
}(x)))\right\vert \leq C\left\Vert f\right\Vert _{\infty }\times (\frac{%
\delta ^{q}}{\eta^{2q}}+\eta^{-\frac{2}{1+\varepsilon_0}}{\gamma_{N(t-1)}^{%
\frac{2}{(2+\varepsilon_0)(1+\varepsilon_0)}}}+\eta^\kappa).
\]

So (\ref{app11}) holds for $p=\frac{2}{1+\varepsilon_0}$ and $\beta=\frac{1}{%
2+\varepsilon_0}$.

\bigskip\bigskip

Finally, we can apply \textbf{Proposition 2.1.1} for $X^{\mathcal{P},{M_%
\mathcal{P}}}_{0,\Gamma_n} $ and $X_{0,\Gamma_n}$ with $\alpha=1,\ k_0=0,\ p=\frac{2%
}{1+\varepsilon_0},\ \beta=\frac{1}{2+\varepsilon_0}$ (for any small $%
\varepsilon_0$), and obtain the following result: for every $\varepsilon>0$,
there exists a constant $C$ such that 
\begin{eqnarray}
d_{TV}(X_{0,\Gamma_n}^{\mathcal{P},{M_\mathcal{P}}},X_{0,\Gamma_n})\leq
C\gamma_n^{\frac{2}{(2+\varepsilon_0)(1+\varepsilon_0)}-\varepsilon}=C%
\gamma_n^{1-\bar{\varepsilon}},  \label{TV1}
\end{eqnarray}%
with $\bar{\varepsilon}=1-\frac{2-\varepsilon(2+\varepsilon_0)(1+%
\varepsilon_0)}{(2+\varepsilon_0)(1+\varepsilon_0)}.$\newline
And moreover, we have 
\[
d_{TV}(\mathcal{L}(X^{\mathcal{P},{M_\mathcal{P}}}_{\Gamma_n}),\nu )\leq
C(\gamma _{n}^{1-\varepsilon }+\int_{\mathbb{R}^d} \left\vert x-y\right\vert
d\nu (y)e^{-\frac{\theta}{2} \Gamma _{n}}),
\]%
where $\nu$ is the unique invariant probability measure.

$\bigskip$

\section{Appendix}

\subsection{The numerical lemma}

In  Section 2, we need to use the following
numerical lemma.

\begin{lemma}
\textbf{(A)} Take an integer $n_{\ast}$. Let $(\gamma_n)_{n\in\mathbb{N}}$
be a non-increasing positive sequence such that for $n\geq n_{\ast}$, we
have 
\begin{equation}
\frac{\gamma _{n}-\gamma _{n+1}}{\gamma _{n+1}^{2}}\leq 2\overline{\omega}.
\label{steptech0}
\end{equation}%
We denote $\Gamma_n=\sum_{i=1}^n\gamma_i$. Then for every $n_{\ast}\leq
i\leq n$, we have 
\begin{equation}
\gamma_{i}\leq e^{2\overline{\omega}(\Gamma_n-\Gamma_i)}\times\gamma_n.
\label{steptech}
\end{equation}

\textbf{(B)} We assume that $(\gamma_n)_{n\in\mathbb{N}}$ is a
non-increasing positive sequence verifying 
\begin{equation}
\frac{\gamma _{n}-\gamma _{n+1}}{\gamma _{n+1}^{2}}\leq c_{\ast }<\frac{\rho 
}{\alpha }.  \label{App13}
\end{equation}%
We denote $\Gamma_n=\sum_{i=1}^n\gamma_i$. Then 
\begin{equation}
u_{n}:=\sum_{i=1}^{n}\gamma _{i}^{1+\alpha }e^{-\rho (\Gamma _{n}-\Gamma
_{i})}\leq C\gamma _{n}^{\alpha }.  \label{App14}
\end{equation}
\end{lemma}

\textbf{Proof of \textit{(A)} }Notice that (\ref{steptech0}) implies%
\[
\frac{\gamma _{n}}{\gamma _{n+1}}\leq 1+2\overline{\omega}\gamma _{n+1}\leq
e^{2\overline{\omega}\gamma _{n+1}}.
\]%
Then 
\[
\frac{\gamma_{i}}{\gamma_n}=\prod_{k=i}^{n-1}\frac{\gamma_k}{\gamma_{k+1}}%
\leq \prod_{k=i}^{n-1}e^{2\overline{\omega}(\gamma_{k+1})}\leq e^{2\overline{%
\omega}(\Gamma_n-\Gamma_i)}.
\]

\textbf{Proof of \textit{(B)} }Notice that (\ref{App13}) implies%
\[
\frac{\gamma _{n}}{\gamma _{n+1}}\leq 1+c_{\ast }\gamma _{n+1}\leq
e^{c_{\ast }\gamma _{n+1}}.
\]%
Then we define $v_{n}=u_{n}/\gamma _{n}^{\alpha }$ and we have the
recurrence relation%
\[
v_{n+1}=\theta _{n}v_{n}+\gamma _{n+1},\quad \theta _{n}=\frac{\gamma
_{n}^{\alpha }}{\gamma _{n+1}^{\alpha }}\times e^{-\rho \gamma _{n+1}}.
\]%
Using the previous inequality we obtain%
\[
v_{n+1}\leq e^{(\alpha c_{\ast }-\rho )\gamma _{n+1}}v_{n}+\gamma _{n+1}
\]%
and further 
\begin{eqnarray*}
e^{(\rho -\alpha c_{\ast })\Gamma _{n+1}}v_{n+1}&\leq& e^{(\rho -\alpha
c_{\ast })\Gamma _{n}}v_{n}+e^{(\rho -\alpha c_{\ast })\Gamma _{n+1}}\gamma
_{n+1} \\
&\leq& e^{(\rho -\alpha c_{\ast })\Gamma _{n}}v_{n}+C^\prime e^{(\rho
-\alpha c_{\ast })\Gamma _{n}}\gamma _{n+1},
\end{eqnarray*}%
with $C^\prime=\sup\limits_{k\geq1}e^{(\rho-\alpha
c_\ast)\gamma_k}=e^{(\rho-\alpha c_\ast)\gamma_1}$. We use recursively this
inequality and we obtain 
\begin{eqnarray*}
e^{(\rho -\alpha c_{\ast })\Gamma _{n+1}}v_{n+1} &\leq &e^{(\rho -\alpha
c_{\ast })\Gamma _{1}}v_{1}+C^\prime\sum_{i=1}^{n}e^{(\rho -\alpha c_{\ast
})\Gamma _{n}}\gamma _{n+1} \\
&\leq &e^{(\rho -\alpha c_{\ast })\Gamma _{1}}v_{1}+C^\prime\int_{0}^{\Gamma
_{n}}e^{(\rho -\alpha c_{\ast })s}ds \\
&\leq &e^{(\rho -\alpha c_{\ast })\Gamma _{1}}v_{1}+\frac{C^\prime}{\rho
-\alpha c_{\ast }}e^{(\rho -\alpha c_{\ast })\Gamma _{n+1}}.
\end{eqnarray*}%
That is 
\[
v_{n+1}\leq v_1+\frac{C^\prime}{\rho -\alpha c_{\ast }}\leq \gamma_1+\frac{%
C^\prime}{\rho -\alpha c_{\ast }}
\]%
which finally gives%
\[
u_{n+1}\leq (\gamma_1+\frac{C^\prime}{\rho -\alpha c_{\ast }})\gamma
_{n+1}^{\alpha }.
\]%
$\square $

\end{document}